\documentclass[12pt]{amsart}
\usepackage{amssymb}
\usepackage{amscd}
\usepackage{amsmath}

\newcommand{\gb}{{\mathfrak b}}
\newcommand{\gog}{{\mathfrak g}}
\newcommand{\gh}{{\mathfrak h}}
\newcommand{\gl}{{\mathfrak l}}
\newcommand{\gm}{{\mathfrak m}}
\newcommand{\gn}{{\mathfrak n}}
\newcommand{\gp}{{\mathfrak p}}
\newcommand{\gs}{{\mathfrak s}}
\newcommand{\bB}{{\bf B}}
\newcommand{\bG}{{\bf G}}
\newcommand{\bH}{{\bf H}}
\newcommand{\bL}{{\bf L}}
\newcommand{\bM}{{\bf M}}

\newcommand{\bP}{{\bf P}}
\newcommand{\bS}{{\bf S}}
\newcommand{\bT}{{\bf T}}
\newcommand{\bV}{{\bf V}}

\newcommand{\bD}{{\bf D}}
\newcommand{\bC}{{\bf C}}

\newcommand{\al}{\alpha}
\newcommand{\be}{\beta}

\newcommand{\om}{\omega}

\newcommand{\cal}{\mathcal} 
\newcommand{\Cscr}{{\cal C}}
\newcommand{\Nscr}{{\cal N}}
\newcommand{\Bscr}{{\cal B}}
\newcommand{\Oscr}{{\cal O}}
\newcommand{\Vscr}{{\cal V}}
\newcommand{\Uscr}{{\cal U}}

\newcommand{\Iscr}{{\cal I}}

\newcommand{\Wscr}{{\cal W}}
\newcommand{\Dscr}{{\cal D}}

\newcommand{\Yscr}{{\cal Y}}
\newcommand{\Lie}{{\rm Lie\,}}
\newcommand{\sh}{{\rm sh\,}}

\newcommand{\Id}{{\rm Id\,}}

\newcommand{\Spa}{{\rm span\,}}
\newcommand{\QED}{\par \hspace{15cm}$\blacksquare$ \par}
\newcommand{\pr}{^{\prime}}
\newcommand{\prpr}{^{\prime\prime}}
\newcommand{\st}{\subset}
\newcommand{\Co}{{\mathbb C}} 
\newcommand{\Na}{{\mathbb N}}  
\newcommand{\boT}{{\mathbb T}}
\newcommand{\boS}{{\mathbb S}}
\newcommand{\Pf}{\noindent{\bf Proof.}\par\noindent}

\newcommand{\sr}{\scriptscriptstyle}
\newcommand{\vb}{\vrule height 14pt depth 7pt} 
\newcommand{\ts}{\tabskip 4pt}	
\newcommand{\vsa}{\noalign{\vskip-7pt}}
\newcommand{\ssa}{\noalign{\vskip -1pt}}

\newcommand{\bk}{\break}
\newcommand{\ov}{\overline}
\newcommand{\rar}{\rightarrow}
\newcommand{\uar}{\uparrow}
\newcommand{\dar}{\downarrow}
\newcommand{\lar}{\leftarrow}

\newcommand{\Llrar}{\Longleftrightarrow}
\newcommand{\Da}{\Downarrow}
\newcommand{\Ra}{\Rightarrow}
\newcommand{\La}{\Leftarrow}
\newcommand{\Ua}{\Uparrow}

\newcommand{\parno}{\par\noindent}
\newcommand{\dor}{\stackrel{\rm D}{\leq}}     
\newcommand{\dg}{\stackrel{\rm D}{\geq}}
\newcommand{\dos}{\stackrel {\rm D}{<}}
\newcommand{\dgs}{\stackrel {\rm D}{>}}
\newcommand{\go}{\stackrel {\rm G}{\leq}}     

\newcommand{\gos}{\stackrel {\rm G}{<}}
\newtheorem*{theorem}{Theorem}
\newtheorem*{lemma}{Lemma}
\newtheorem*{defi}{Definition}
\newtheorem*{cor}{Corollary}
\newtheorem*{prop}{Proposition}

\begin{document}

\title[Duflo order]
{\bf On orbital variety closures in $\mathfrak{sl}_n$\\ 
I. Induced Duflo order} 
\author{Anna Melnikov}\thanks{This work
was partially supported by the EEC program TMR-grant ERB FMRX-C
T97-0100}
\address{Department of Mathematics,
University of Haifa,
Haifa 31905, Israel}
\email{melnikov@math.haifa.ac.il}

\begin{abstract}
For  a semisimple Lie algebra $\gog$ the orbit
method attempts to assign representations
of $\gog$ to (coadjoint) orbits in $\gog^*.$  Orbital
varieties are particular Lagrangian subvarieties of such orbits leading 
to highest weight representations of $\gog.$ 
In ${\gs\gl}_n$ orbital varieties are described by 
Young tableaux. Inclusion relation on orbital variety closures
defines a partial order on Young tableaux.
Our aim is to describe this order.
The paper is devoted to
the combinatorial description of induced Duflo order
on Young tableaux (the order generated by inclusion of
generating subspaces of orbital varieties). This is a very 
interesting and complex combinatorial
question. 
\par
This is the first paper in the series.
In Part II and Part III we use repeatedly the results of the paper
as a  basis for further study of orbital variety 
closures.
\end{abstract} 
\maketitle

\section {\bf Introduction}

This is the first paper in the series of three 
papers devoted to the study of orbital variety closures in 
$\mathfrak{sl}_n.$ They are referred  to as Part I, Part II and Part III 
respectively.
\subsection{}\label{1.1}
The orbital varieties derive from
the works of N. Spaltenstein ~\cite{Sp1} and ~\cite{Sp2}, 
and R. Steinberg ~\cite{St1} and \cite{St2} during their studies of unipotent 
variety of
a complex semi-simple group $\bf G.$ 
Let $\Bscr$ be the variety of Borel subgroups of $\bf G$ on which	
$\bf G$ acts by conjugation. For a fixed unipotent 
$u\in \bf G$
let ${\Bscr}_u$ be the subvariety of 
$\Bscr$ containing $u$
or equivalently the variety of flags in $\bf G/\bf B$ fixed by $u$ 
for some fixed Borel subgroup $\bf B.$ 
Spaltenstein and Steinberg studied the irreducible components of this 
variety. 
\par
Orbital varieties are 
the translation of these components from unipotent
variety of $\bf G$ to nilpotent cone of $\mathfrak g =\rm {Lie}(\bf G).$
We give their description in the next subsection.  
\par
\subsection{}\label{1.2}
Let $\bf G$ be a connected semisimple 
finite dimensional
complex algebraic group. Let $\mathfrak g$ be its Lie algebra  and
$U(\mathfrak g)$ be the enveloping algebra of $\mathfrak g.$
Consider the adjoint action of $\bf G$ on $\mathfrak g.$
Fix some triangular decomposition  $\mathfrak g=\mathfrak n\bigoplus
\mathfrak h\bigoplus\mathfrak n^-.$  A $\bf G$ orbit $\Oscr$ in $\mathfrak g$ is called
nilpotent if it consists of nilpotent elements, that is if
${\Oscr}={\bf G}\, x$ for some $x\in\mathfrak n.$ 
The intersection $\Oscr\cap\mathfrak n$ is reducible in general. 
Its irreducible components
are called orbital varieties associated to $\Oscr.$
Orbital varieties play a key role 
in the study of
primitive ideals in $U(\mathfrak g).$ They also play an important
role in Springer's Weyl group representations, 
described in terms of ${\Bscr}_u.$
\subsection{}\label{1.3}
The first role above can be detailed as follows.
Since $\mathfrak g$ is semisimple we can
identify $\mathfrak g^*$ with $\mathfrak g$ through the Killing
form. 
This identification gives an
adjoint orbit a  symplectic structure.
Let $\Vscr$ be an orbital
variety associated to $\Oscr.$
By ~\cite{Sp2} and ~\cite{St1}
one has $\dim {\Vscr}=\frac{1}{2}\dim{\Oscr}.$
Moreover as it was pointed out in ~\cite{J} this
implies that an orbital variety is a Lagrangian
subvariety of the nilpotent orbit it is associated to.
Following the orbit method one would like to attach an 
irreducible representation of $U(\mathfrak g)$ to $\Vscr.$ This
should be a simple highest weight module.
Combining the results of A. Joseph and T. A. Springer
one obtains  
a one to one correspondence
between the set of primitive ideals of $U(\mathfrak g)$ 
containing the augmentation ideal of its centre (thus corresponding to
integral weights) 
and the set of orbital varieties in $\mathfrak g$ corresponding
to Lusztig's special orbits (see for example ~\cite{B-B}).
The picture is especially beautiful for $\mathfrak g=\mathfrak s\mathfrak l_n.$
In this case  all orbits are special and by ~\cite{M} 
the associated variety of a simple highest (integral)
weight module is irreducible. By ~\cite{B-B} and ~\cite{J} in general
the orbital variety closures are just the irreducible 
components of an associated variety of a  simple 
highest weight module. Thus for $\mathfrak g=\mathfrak s\mathfrak l_n$
orbital variety closures are these associated varieties and
therefore give a natural geometric understanding of 
the classification of primitive ideals. This makes their study
especially interesting. 
\subsection{}\label{1.4}
Orbital varieties are very interesting objects
from algebro-geometric point of view as well. 
Given an orbital variety $\Vscr$ one can easily find 
$\mathfrak m_{\sr \Vscr}$ -- the nilradical of the smallest dimension 
containing $\Vscr$ (see \ref{2.1.7}). Consider an orbital variety 
closure as an algebraic variety in the affine linear space  $\mathfrak m_{\sr \Vscr}.$
Then vast majority of orbital varieties are not complete 
intersections.  
So orbital varieties are examples of algebraic varieties
which are both Lagrangian subvarieties and not 
complete intersections.
\subsection{}\label{1.5}
Orbital varieties still remain rather
mysterious objects. The only general description
was given by R. Steinberg ~\cite{St1}. It is explained in detail in \ref{2.1.2}
and \ref{2.1.3}. Briefly speaking let $\bf B$ be the standard
Borel subgroup of $\bf G$, i.e.  such that $\rm {Lie}({\bf B})=\mathfrak b=\mathfrak h\bigoplus\mathfrak n.$
$\bf B$ acts by conjugation on $\mathfrak n$ and its subsets.
Let $W$ be the Weyl group for the pair $(\mathfrak g,\ \mathfrak h).$ 
Then by ~\cite{St1}
there exists a surjection $\phi: w \mapsto \ov{\Vscr}_w$ from the 
Weyl group onto the set of orbital variety closures defined by 
$\phi(w)=\ov{{\bf B}(\mathfrak n\cap^w\mathfrak n)}=:\ov{\Vscr}_w.$ 
The fibers of this mapping, namely $\phi^{-1}(\ov {\Vscr})=\{w\in W\ :\ 
\ov{\Vscr}_w=\ov{\Vscr}\}$ are called geometric cells. 
\par
This description is not very satisfactory 
from the geometric 
point of view since a $\bf B$ invariant subvariety
generated by a linear space is a very complex object. 
\subsection{}\label{1.6}
On the other hand there exist a very nice combinatorial
characterization of orbital varieties in $\mathfrak s\mathfrak l_n$ in terms of
Young tableaux. The detailed description of this characterization
is given in \ref{2.4.5}, \ref{2.4.6}, \ref{2.4.17} and \ref{2.4.18}. It is defined by 
Robinson-Schensted procedure giving a bijection from the symmetric group 
${\bf S}_n$ onto the pairs of standard Young tableaux of the same 
shape $w\mapsto (T(w),Q(w))$. 
Let us identify $W$ with ${\bf S}_n$ (see \ref{2.2.2}). Then
geometric cells are given by Young tableaux as follows
$\phi^{\sr -1}({\Vscr}_w)=\{y\, :\, T(y)=T(w)\}.$ We will denote
$T_{{\Vscr}_w}:=T(w).$
\subsection{}\label{1.7}
The description of an orbital
variety closure has both geometric and 
combinatorial parts. The geometric component is whether an 
orbital variety closure is a union of orbital varieties.
The combinatorial part is to describe orbital variety closure in terms
of manipulations on Young tableaux.
\par
It is noted in \ref{4.1.1} that the projections
on the Levi factor of standard parabolic subalgebras of $\mathfrak g$ 
preserve orbital variety closures. Using this fact 
together with
computations in low rank cases we show in Part III
that if $\mathfrak g$ has factors not of type $A_n$ then
a closure of orbital variety is not necessarily a union of 
orbital varieties and includes some varieties of
smaller dimensions. 
\par
The same argument does not
work for $\mathfrak s\mathfrak l_n$ and there we conjecture
that an orbital variety closure is a union
of orbital varieties. This conjecture is supported by computations
for $n\leq 6.$ As well it is true for some special cases. In particular
this is true for orbital varieties of nilpotent order 2 as it is shown in ~\cite{M1}, and for
orbital varieties whose closure is a nilradical of some standard
parabolic subalgebras as it is shown in Part II.
\subsection{}\label{1.8}
To specify the combinatorial part let us define
the geometric order and the notion of a geometric descendant.
Given orbital varieties $\Vscr,\ \Wscr$ we define the
geometric order by ${\Vscr}\go\Wscr$ if $\Wscr\st\ov{\Vscr}.$
We say that $\Wscr\ne \Vscr$ is a geometric descendant of $\Vscr$
if $\Vscr\go\Wscr$ and for any $\Yscr$ such that 
${\Vscr}\go\Yscr\go\Wscr$ one has $\Yscr=\Vscr$ or $\Yscr=\Wscr.$
\parno
{\bf Remark.}\ \ The above definition of geometric order seems to be reverse to the obvious one
defined by the inclusion of closures.
Initially this order was introduced for inclusions of primitive ideals
described in \ref{1.3}. In $\gs\gl_n$ as it is shown in ~\cite{M} $I_w\subset I_y$
implies $\Vscr_{w^{-1}}\supset \Vscr_{y^{-1}}.$ So that the ``right'' definition
for primitive spectrum induces the ``reverse'' definition for orbital varieties. 
We will discuss the relationship between ordering of orbital varieties
and ordering of primitive spectrum in detail in Part III.
\par
We wish to describe the set of descendants of an orbital
varieties in terms of Young tableaux and in particular to determine whether
${\Vscr}\go\Wscr$ can be described by regarding  their Young tableaux.
\par
Defining the order on nilpotent orbits  (resp. a descendant of a 
nilpotent orbit) exactly in the same manner as for orbital
varieties we can ask the same questions about nilpotent orbits.
The construction of Gerstenhaber described in \ref{2.3.2} gives
very elegant combinatorial answers to both questions in terms
of Young diagrams.
\par
We would like to find similar answers for orbital variety 
closures. As we show this is much more complex.
For example given an orbital variety $\Vscr,$ let
${\Oscr}_{\Vscr}={\bf G}\Vscr$ be the nilpotent orbit, $\Vscr$ is associated to.
Then, as we show in Part II, $\Wscr$ being a geometric descendant 
of $\Vscr$ does not imply necessarily that ${\Oscr}_{\Wscr}$ is a
descendant of ${\Oscr}_{\Vscr}.$
\subsection{}\label{1.9}
Let us consider another order relation on orbital varieties,
which had long been thought to be the same as the geometric order.
\par
Let $R\st\mathfrak h^*$ denote the set of non-zero roots,  $R^+$ the set
of positive roots corresponding to $\mathfrak n$ in the triangular 
decomposition of $\mathfrak g$ and $\Pi\st R^+$ the resulting set of 
simple roots. 
Each $w$ in $W$ is a product of fundamental reflections $s_{\al}\ :\ \al\in\Pi.$ 
We denote by $\ell (w)$
the minimal  length  of any such expression for $w$. Consider
the order generated by the following preorder. 
For  $s_{\al},\ \al\in \Pi$ and $w\in W$ put 
$$w \dor ws_{\al}\quad {\rm if}\quad \ell(w s_{\al})=\ell(w)+1. \eqno{(*)}$$
We call it the (right) Duflo order. This is also known as the weak (right)
Bruhat order. We prefer the former nomenclature in
the present context since it was Duflo ~\cite{D} who first discovered 
the implication of the (left) Duflo order for the primitive spectrum.
\par
In \ref{2.1.6} we induce the Duflo order to orbital varieties,
geometric cells and Young tableaux in the obvious manner and 
call it the induced Duflo 
order. It is weaker than the geometric order, that is ${\Vscr}\dor \Wscr$
implies ${\Vscr}\go\Wscr.$
\par
As we show in Part III the induced Duflo order coincides with the 
geometric order up to $n=5$ and it is strictly weaker than the geometric
order for $n\geq 6.$
\subsection{}\label{1.10}
Part I is devoted to the combinatorial description
of the induced Duflo order on Young tableaux. This is a very 
interesting and complex combinatorial
question. In Part II and Part III we use repeatedly this 
description as a  basis for further study of orbital variety 
closures.
\par
All nilradicals of standard parabolic subgroups are orbital
variety closures. They are called also Richardson orbital varieties.
They are the simplest examples of orbital variety closures since
they are linear subspaces of $\mathfrak n.$ For Richardson orbital variety 
$\Vscr$  the question whether ${\Vscr}\go\Wscr$ has a very simple 
answer in combinatorics of Young tableaux. Consider the invariant 
$\tau(T)$ determined in \ref{2.4.14}. Then ${\Vscr}\go\Wscr$
if and only if $\tau(T_{\Wscr})\supset \tau(T_{\Vscr}).$ The description
of the set of descendants of Richardson orbital variety $\Vscr$ is
a much more delicate problem. We give the full solution to this 
problem in Part II.
\par
In Part III we study geometric properties of
orbital variety closures. In particular we use the Vogan 
$T_{\al,\be}$ operators to strengthen the induced Duflo order. 
This gives
a combinatorial description of the geometric order up to at least
$n=10.$ In Part III we also discuss in detail a 
connection between primitive ideals and orbital varieties mentioned 
briefly in \ref{1.8} and use it to study the properties  of orbital 
varieties as well as primitive ideals in $\mathfrak s\mathfrak l_n.$
\subsection{}\label{1.11}
Let us describe in more detail the results of
Part I. Since the induced Duflo order is generated by \ref{1.9} $(*)$
we call ${\Wscr}={{\Vscr}}_{{ws_{\al}}}$ (resp. Young tableau $S=T(ws_{\al})$)
a Duflo offspring of ${\Vscr}={\Vscr}_w$ (resp. of $T=T(w)$) if 
$\al\in \Pi$ and $\ell(ws_a)=\ell(w)+1$ (see \ref{2.5.2}). We define
the Duflo descendant of a given orbital variety or Young tableau 
with respect to the induced
Duflo order exactly in the same manner as we have defined the geometric 
descendant of a given orbital variety in \ref{1.8} (see \ref{2.5.1}). 
It is obvious that the set of Duflo descendants of a given
orbital variety or Young tableau is a subset of its set of 
Duflo offsprings.
But it is much easier to describe the set of Duflo offsprings
than the set of Duflo descendants since the first set has
a definitive combinatorial definition.
\par
We show  that the induced Duflo order can be 
completely described by the natural
ordering on ${\bf S}_2$
and the Robinson - Schensted insertion of an 
element $a$ into a Young tableau $T$ from the 
above $T\mapsto (T\Da a)$
and from the left $T\mapsto (a\Ra T)$   Precisely, in \ref{3.4.5} we show
that for any $S$ in the set of Duflo offsprings of $T$ there
exists $a$ such that $T=(T\pr\Da a)$ or $T=(a\Ra T\pr)$
and $S=(S\pr\Da a)$ or respectively $S=(a\Ra S\pr)$ where
$S\pr$ is an offspring of $T\pr.$ We also provide  in \ref{3.3.3} 
an exact way to compute the set of 
offsprings of a Young tableau of size $n$ from the knowledge of 
sets of offsprings of Young tableaux of size $n-1.$ It involves
the shuffling of numbered boxes in a manner 
prescribed by the Robinson - Schensted procedure
and, technically, is the most difficult theorem of 
this work.
\subsection{}\label{1.12}
Using the results described in \ref{1.11} we show
that the set of offsprings is preserved under projections
on a Levi factor (see \ref{4.1.3}) and under Robinson-Schensted 
embeddings $\Da$ and 
$\Ra$ (see \ref{4.1.4}). 
\par
In \ref{4.1.2} we show as well that the geometric order is preserved
under projections on a Levi factor. In Part III we show that
$\Da$ and $\Ra$ preserve the geometric order as well.
\par
However as we show in \S \ref{4.2} neither projections nor embeddings
preserve the set of Duflo or geometric descendants. These
facts again underline the difficulty of the constructing
the set of descendants.
\par
Using the results on embeddings and projections we
show in \ref{4.1.8} that the
induced Duflo order is compatible with the order 
on nilpotent
orbits in the following strong sense: 
if ${{\Oscr}}_1,\ {\Oscr}_2$  are nilpotent
orbits  and ${\Oscr}_1 \subset \ov {\Oscr}_2$ 
then for every
orbital variety
${\Vscr}_2 \st {\Oscr}_2$  there exist an orbital variety
${\Vscr}_1 \st {\Oscr}_1$  such that  ${\Vscr}_2 \dor  {\Vscr}_1.$ 
This property is important in
consideration of an orbital variety closure and
fails to be true outside	
of $\mathfrak s\mathfrak l_n$ as we show in Part III.	 
\subsection{}\label{1.13}
The body of the paper consists of three  sections.
\par
In section 2 we explain all the background in geometry
of orbital varieties and combinatorics of Young tableaux
essential in  the subsequent analysis. 
I hope this part makes the paper 
self-contained.
\par
In section 3  we work out the machinery for recursive construction
of the set of Duflo offsprings  with the main 
results stated in \ref{3.3.3} and \ref{3.4.5}.
\par
Finally section 4 is devoted to the study of the properties of the
induced Duflo order resulting immediately from section 3. Further
properties of the induced Duflo order are studied in Part II.
\par
In the end one can find the index of notation in which 
symbols appearing frequently are given with the subsection
where they are defined. We hope that this will help the reader
to find his way through the paper.  
\parno
{\bf Acknowledgments.}\ \ \ I would like to express my deep gratitude
to A. Joseph for introducing the world of orbital varieties
to me, for posing the problems, suggesting ideas underlying 
this research and many fruitful discussions through the various 
stages of this work. I would also like to thank
V. Hinich for fruitful discussions.

\section {\bf Combinatorics of Symmetric group}
\subsection{\bf Steinberg map and Duflo order on Weyl groups}

\subsubsection{}\label{2.1.1}
Recall the notation from \ref{1.2}. 
Let $\Nscr={\bG}(\gn)$ denote the nilpotent cone in $\gog.$ 
For $u\in\Nscr$ let ${\Oscr}_u$ be the 
nilpotent orbit it defines, that is the orbit of $u$
under the coadjoint action of $\bG.$
We define an order relation  on the set of nilpotent orbits  by   
$$ {\Oscr}_u \geq {\Oscr}_v \ \ {\rm if} \ \ {\Oscr}_u \st\ov{\Oscr}_v \ . $$
\par
Let ${\Oscr}_u\ne {\Oscr}_v$ and ${\Oscr}_u\geq {\Oscr}_v.$
We call ${\Oscr}\pr$ a descendant of $\Oscr$ if for any $\Oscr\prpr$ such that
$\Oscr\geq\Oscr\prpr \geq\Oscr\pr$  one has $\Oscr\prpr=\Oscr$  or
$\Oscr\prpr=\Oscr\pr.$
\subsubsection{}\label{2.1.2}
An orbital variety $\Vscr$ associated to a nilpotent
orbit $\Oscr$ is an irreducible component of $\Oscr\cap \gn.$
Let $\bV$ denote the set of all orbital varieties of $\gog$ and 
$W$ the Weyl group for the pair $(\gog,\ \gh).$ We describe
first the Steinberg map of $W$ onto $\bV.$
Let $R\st\gh^*$ denote the set of non-zero roots,  $R^+$ the set
of positive roots corresponding to $\gn$ in the triangular 
decomposition of $\gog$ and $\Pi\st R^+$ the resulting set of 
simple roots. Let $X_{\al}={\Co} x_{\al}$ denote the root subspace 
corresponding to $\al\in R.$ 

Then $\gn=\bigoplus\limits_{\al\in R^+}X_{\al}$
(resp.$\gn^-=\bigoplus\limits_{\al\in -R^+}X_{\al}$).

\subsubsection{}\label{2.1.3}
Given $S,T\st R$ and $w\in W$ we set 
$S\cap^w T:=S\cap w(T)=\{\al\in S\ :\ \al\in w(T)\}.$
Then set
$$\gn\cap^w\gn:=\bigoplus\limits_{\al\in R^+\cap^w R^+}X_{\al}.$$
This is a subspace of $\gn.$ For each closed, irreducible  subgroup  $\bH$ of $\bG$ let
${\bH}(\gn\cap^w\gn)$ be the set of $\bH$
conjugates of $\gn\cap^w\gn.$ It is an irreducible locally 
closed subvariety. Since there are only finitely many nilpotent 
orbits in $\gog$ it follows that there exists a unique 
nilpotent orbit which we denote by ${\Oscr}_w$ such that 
$\ov{{\bG}(\gn\cap^w\gn)}=\ov{\Oscr}_{w}.$ 
A result of Steinberg ~\cite{St1} asserts that 
${\Vscr}_w:=\ov{{\bB}(\gn\cap^w\gn)}\cap{\Oscr}_w$ is an orbital
variety and that the map
$\varphi:w\mapsto {\Vscr}_w$
is a surjection of $W$ onto $\bV.$

\subsubsection{}\label{2.1.4}
Recall the notion of geometric order and geometric
descendant from \ref{1.8}.
Set ${\Vscr}_2\gos {\Vscr}_1$ if ${\Vscr}_2\go{\Vscr}_1$  and ${\Vscr}_2\ne {\Vscr}_1.$
\par
The partial order relation on the corresponding Weyl group $W$  
defined by inclusion of orbital variety closures takes the form
$$ w\go y \ \ {\rm if} \ \ \ov{{\bB}(\gn \cap^y\gn)}
\st \ov{{\bB}(\gn \cap^w\gn)} \ . $$
\par
We decompose $W$ according to this relation into the cells:
$${\Cscr}_w=\{y\in W:{\Vscr}_y={\Vscr}_w\}\ . $$
\par
We call these the (right) geometric cells of $W$ following ~\cite{St1}, 
and 
~\cite[6.8]{B-B}.
\par
Set $ {\Cscr}_y\go {\Cscr}_w$  (resp. $ {\Cscr}_y\gos {\Cscr}_w$) if  ${\Vscr}_y\go {\Vscr}_w$
(resp. $ {\Vscr}_y\gos {\Vscr}_w$). 
We call ${\Cscr}_w$ a geometric descendant of ${\Cscr}_y$
if ${\Vscr}_w$ is a geometric descendant of ${\Vscr}_y.$
\par
As well we decompose $W$ into double geometric cells 
according to nilpotent orbits
as follows
$${\Cscr}^d_w=\{y\in W\, :\, {\Oscr}_y={\Oscr}_w\}.$$
One sees immediately that ${\Cscr}^d_w$ is the 
union of geometric cells whose orbital varieties
are attached to ${\Oscr}_w.$ 

\subsubsection{}\label{2.1.5}
Recall $\ell (w)$ from \ref{1.9}.
Set 
$S(w):=R^+ \cap^w R^-,$ then a classical result, described, for 
instance in ~\cite[\S 2.2]{Ca} provides
\begin{lemma} 
Take $x, y\in W$  and  set  $w=yx \ .$
  Then $\ell(w)=\ell(x) + \ell(y)$  if and only if
$$ S(y)\st S(w)  \eqno{(*)}  $$
 Moreover if $(*)$ holds there are exactly  $\ell(x)$  
roots $\be\in R^+$ such
that $y(\be) \in R^+ $ and $w(\be)\in R^-$.
\end{lemma}
It is known  that the map  $w\mapsto S(w)$ is injective. 
Hence we may define a partial order relation on $W$  
by $y \dor w$  when $(*)$ holds. 
It is called the Duflo order. It is 
generated by the preorder defined in \ref{1.9}.
\par
Note that if $y\dor w$ then $R^+ \cap^y R^+\supset R^+ \cap^w R^+$
thus  $\ov{\Vscr}_y\supset {\Vscr}_w$ just by the inclusion of generating subspaces.

\subsubsection{}\label{2.1.6}
We induce Duflo order on orbital varieties
(resp. on geometric cells) by the 
following. Set ${\Vscr}_y\dor {\Vscr}_w$ 
(resp. ${\Cscr}_y\dor {\Cscr}_w$) if there exist a chain 
$y=x_o^1,x_o^2, x_1^1,x_1^2,\cdots,x_k^1, x_k^2=w$ such that 
${\Cscr}_{x_i^1}={\Cscr}_{x_i^2}$ for $0\leq i\leq k$ and 
$x_i^2\dor x_{i+1}^1$ for $0\leq i\leq k-1.$
Set ${\Vscr}_y \dos {\Vscr}_w$ if  ${\Vscr}_y\dor {\Vscr}_w$ and 
${\Vscr}_y\ne {\Vscr}_w.$ In the same fashion we define 
${\Cscr}_y\dos {\Cscr}_w.$
\par 
By note \ref{2.1.5} we get that the induced Duflo order is weaker than the
geometric order, that is for cells ${\Cscr}_y,{\Cscr}_w$ the relation 
${\Cscr}_y\dor {\Cscr}_w$
implies ${\Cscr}_y\go {\Cscr}_w.$ 
\subsubsection{}
\label{2.1.7}
Recall the standard Borel subalgebra $\gb$ from 
\ref{1.5}.
For $\al\in\Pi$ let ${\bP}_{\al}$ be the standard 
parabolic subgroup with $\Lie({\bP}_{\al})=\gp_{\al}:=\gb\oplus X_{-\al}.$
\par
Given an orbital variety $\Vscr,$ let ${\bP}_{\Vscr}$
be its stabilizer in $\bG.$ This is a standard parabolic
(that is ${\bP}_{\Vscr}\supset\bB$) subgroup of $\bG.$ Let $\gm_{\sr \Vscr}$
be the maximal subalgebra of $\gog$ stabilized by ${\bP}_{\Vscr}.$
This is a nilradical and the linear subspace of $\gn$ of minimal
possible dimension containing $\Vscr.$
\par
Take $w\in W,\ {\Vscr}\in\bV,$
a standard parabolic subgroup $\bP$ and 
a standard parabolic subalgebra $\gp=\Lie(\bP).$
Define their $\tau$-invariants to be
$$ 
\begin{array}{rcl}
\tau(w)&:=&\Pi\cap S(w),\\
\tau({\bP})& :=&\{\al\in\Pi\ :{\bP}_{\al}\st {\bP}\},\\
\tau(\gp)& :=&\{\al\in \Pi \ :\gp_{\al}\st \gp\},\\
\tau(\Vscr)&:=&\{\al\in\Pi\ : {\bP}_{\al}(\Vscr)=\Vscr\}.\\
\end{array}
$$
Note that $\bP$ (resp. $\gp$) is uniquely determined by 
its $\tau$-invariant.
One has (see ~\cite[\S 9]{J}) 
\begin{lemma} 
$\tau(w)=\tau({\Vscr}_w)=\tau({\bP}_{{\Vscr}_w}).$
\end{lemma} 
Therefore $\tau({\Cscr}_w):=\tau(w)$ is well defined.
\subsubsection{}\label{2.1.8}
Given $\Iscr\st \Pi,$ let ${\bP}_{\Iscr}$ denote the unique 
standard parabolic subgroup of $\bG$ such that 
$\tau({\bP}_{\Iscr})=\Iscr.$ Let ${\bM}_{\Iscr}$ be the unipotent 
radical of ${\bP}_{\Iscr}$ and ${\bL}_{\Iscr}$ a Levi factor.
Let $\gp_{\sr \Iscr},\ \gm_{\sr \Iscr},\ \gl_{\sr \Iscr}$ denote 
the corresponding Lie algebras.
Set ${\bB}_{\Iscr}:={\bB}\cap{\bL}_{\Iscr}$ and
$\gn_{\sr \Iscr}:=\gn\cap\gl_{\sr \Iscr}.$
 We have decompositions ${\bB}={\bM}_{\Iscr}\ltimes{\bB}_\Iscr$ and 
$\gn=\gn_{\sr \Iscr}\oplus\gm_{\sr \Iscr}.$ They define projections 
${\bB}\rar {\bB}_\Iscr$ and $\gn\rar\gn_{\sr \Iscr}$ which we denote by 
$\pi_{\sr \Iscr}.$
\par
Set $W_\Iscr:=<s_{\al}\ :\ \al\in \Iscr>$ to be a 
parabolic subgroup of $W,$ set
$D_\Iscr:=\{w\in W\ :\ w(\al)\in R^+\ \forall \ \al\in \Iscr\}.$ 
Set $R_\Iscr^+=R^+\cap \Spa(\Iscr).$
A well-known result described
for example in ~\cite[2.5.8]{Ca}  gives
\begin{lemma} 
Each $w\in W$ has a unique expression of the 
form $w=w_{\sr \Iscr}d_{\sr \Iscr}$
where $d_{\sr \Iscr}\in D_\Iscr,\ w_{\sr \Iscr}\in W_\Iscr$ and $\ell(w)=\ell(w_{\sr \Iscr})+
\ell(d_{\sr \Iscr}).$  Moreover
$$R_\Iscr^+\cap^w R^+=R^+_\Iscr\cap^{w_{\sr \Iscr}} R^+_\Iscr.$$
\end{lemma}
The decomposition $W=W_\Iscr\times D_\Iscr$ explained in the lemma  
defines a projection $\pi_{\sr \Iscr}:W\rar W_\Iscr.$
\par
Set $D^{-1}_J:=\{f^{-1}\ :\ f\in D_J\}.$
Applying the lemma  to $w^{-1}$
we get that each $w\in W$ has a unique expression 
of the form $w=f_{\sr J}w_{\sr J}$ where 
$w_{\sr J}\in W_J,\ f_{\sr J}\in D^{-1}_J$ and $\ell(w)=\ell(w_{\sr J})+
\ell(f_{\sr J}).$
\subsection {\bf   ${\bS}_n$ as a Weyl group of $\gs\gl_n$}

\subsubsection{}\label{2.2.1}
From now and on we consider only $\gs\gl_n.$
It is convenient to replace $\gs\gl_n$ by $\gog=\gog\gl_n.$
This obviously makes no difference when we consider
nilpotent cone and adjoint action (conjugation) by ${\bG}={\bf GL}_n$ 
on it. Let $\gn$ be  the subalgebra
of strictly upper-triangular matrices and let $\gn^-$ 
be the subalgebra of strictly lower-triangular matrices.
Let $\bB$ be the (Borel) subgroup of upper-triangular
matrices in $\bG.$ All parabolic
subgroups we consider further are standard, that is contain
$\bB.$
Let $e_{i,j}$ be
the matrix having $1$ in the $ij-$th entry
and $0$ elsewhere. Set $B=\{e_{i,j}\}_{i,j=1}^n$ which is a basis 
of $\gog.$
\par
Take $i<j$ and let $\al_{i,j}$ be the root corresponding to
$e_{i,j}.$ Set $\al_{j,i}=-\al_{i,j}.$ We write $\al_{i,i+1}$
simply as $\al_i.$ Then $\Pi=\{\al_i\}_{i=1}^{n-1}.$ Moreover
$\al_{i,j}\in R^+$ exactly when $i<j.$ One has 
$$\al_{i,j}=\begin{cases} 
\sum\limits_{k=i}^{j-1}\al_k & {\rm if}\ i>j\\
                   -\sum\limits_{k=j}^{i-1}\al_k & {\rm if}\  i>j\\
\end{cases}
$$
For each $\al\in \Pi,$
let $s_{\al}\in W$ be the corresponding reflection and set 
$s_i=s_{\al_i}.$ 

\subsubsection{}\label{2.2.2}
We represent every element of the symmetric
group ${\bS}_n$  in word form 
$$ w =[a_{\sr 1}, a_{\sr 2},\ldots, a_n]\ ,\quad {\rm where}
                                    \ a_i=w(i).\eqno{(*)}$$
We identify $W$ with ${\bS}_n$ by taking $s_i$ to be the 
elementary permutation interchanging $i,\ i+1.$ We consider
multiplication from right to left that is given $w$ from $(*)$
one has $s_iw(j)=s_i(a_j).$
For example $ s_{\sr 1} s_{\sr 2} =[2,3,1].$
\parno
{\bf Remark.}\ \ In our notation if $w=[a_{\sr 1},\ldots,a_i,a_{i\sr + 1},\ldots,a_n]$ then
$ws_i$ is obtained from $w$ by interchanging $a_i$ and $a_{i\sr +1}$
that is $ws_i=[a_{\sr 1},\ldots,a_{i+\sr 1},a_i,\ldots,a_n].$
\par
\begin{defi} 
A word (of length $k$) is an ordered array 
$[a_1,\cdots,a_k],$ where all $a_i$ are distinct and
$\{a_i\}_{i=1}^k\subset \{j\}_{j=1}^n,\ n>k.$ If  
$\{a_i\}_{i=1}^k=\{j\}_{j=1}^k$ the word is called standard.
\end{defi}
Expression $(*)$ for $w$ is a standard word of length $n.$
\subsubsection{}\label{2.2.3}
Given $w=[a_1,\cdots,a_n].$ 
Set $p_w(i):=j$ if $a_j=i,$ that is
$p_w(i)$ is the place (index) of $i$ in the word form
of $w.$ One has $w(p_w(i))=w(j)=a_j=i,$ that is 
$p_w(i)=w^{-1}(i).$ As well one has
\begin{lemma} 
$w(\al_{i,j})=\al_{w(i),w(j)}.$
\end{lemma}
\Pf
It is enough to show only for $i<j.$
\par
One has
$$s_k(m)=\begin{cases} 
m & {\rm if}\ m\ne k,k+1\\
                 k+1 & {\rm if}\ m=k\\
                 k   & {\rm if}\ m=k+1\\.
\end{cases}$$
As well
$$s_k(\al_{m,m+1})=s_k(\al_m)=\begin{cases}
-\al_k & {\rm if}\ m=k\\
 \al_{k-1,k+1} & {\rm if}\ m=k-1\\
 \al_{k,k+2}   & {\rm if}\ m=k+1\\
 \al_m & {\rm if}\ m\ne k-1, k,k+1\\
\end{cases}$$ 
Thus $s_k(\al_{m,m+1})=\al_{s_k(m),s_k(m+1)}.$ Applying this to $\al_{i,j}$ we get 
$$s_k(\al_{i,j})=s_k(\sum\limits_{m=i}^{j-1}\al_m)=\al_{s_k(i),s_k(j)}.$$
\par
Now we can show lemma by induction on $\ell(w).$ Indeed let $w=s_k y$
where $\ell(w)=\ell(y)+1$ then
$$
\begin{array}{rll}
w(\al_{i,j})&=s_k y(\al_{i,j})=s_k(\al_{y(i),y(j)}) &{\rm by\ induction\ hypothesis}\\
               &=\al_{s_k(y(i)),s_k(y(j))}& \\
               &=\al_{w(i),w(j)} &{\rm by\ the\ previous\ computations}\\
\end{array}$$
\QED
\subsubsection{}\label{2.2.4}
Recall the notation $S(w)$ from \ref{2.1.5}. 
We get from \ref{2.2.3}  
\begin{cor}[{~\cite[2.3]{JM}}] 
Take $i<j.$ Then $\al_{i,j}\in S(w)$
if and only if $p_w(i)>p_w(j).$
\end{cor}
\Pf
Indeed $S(w)=\{\al\in R^+\ : \al\in w(R^-)\}=\{\al\in R^+\ : w^{-1}(\al)\in R^-\}$
and by \ref{2.2.3} for $i<j$ one has $w^{-1}(\al_{i,j})\in R^-$ if and only if $w^{-1}(i)>w^{-1}(j)$
which is equivalent to $p(i)>p(j).$
\QED
\parno
{\bf Remark 1.}\ \  In particular $\al_i\in \tau(w)$ exactly
when $i+1$ comes before $i$ in the word form of $w.$
This is of course well-known.
\parno
{\bf Remark 2.} As a corollary of \ref{2.2.2}, \ref{2.1.5} and \ref{2.2.4} we get that 
$w\dos ws_i$
if and only if $a_i<a_{i+\sr 1}.$ This is also well-known and given
for example in  ~\cite[pp. 73-74]{Ja}.
\parno
{\bf Remark 3.} As a corollary of \ref{2.2.2}, \ref{2.1.5} and \ref{2.2.4} we get that
${\bS}_n$ has the unique longest element $w_o=[n,n-1,\ldots,2,1].$
For any $w\in {\bS}_n$ one has $w\dor w_o.$ This is also well-known.  

\subsubsection{}\label{2.2.5}
Introduce the following useful notational conventions.
\par
Given a word $w=[a_1,\cdots, a_n]$ we denote by $<w>:=\{a_i\}_{i=1}^n$ the 
set of its entries. 
\begin{itemize}
\item[(i)] Given $m\in <w>,$ set $(w-m)$ to be a word obtained
from $w$ by deleting $m,$ that is if $m=a_i,$ then 
$(w-m):=[a_{\sr 1},\ldots,a_{i-\sr 1},a_{i+\sr 1}, \ldots,a_n].$
\item[(ii)] For $i<j$ set $w_{[i,j]}:=[a_i,a_{i+\sr 1}\cdots, a_j].$ 
\item[(iii)] Set $\ov w$ to be the word with  order reverse to the order of $w$, that is 
if $w=[a_1,\cdots,a_n]$ then $\ov w:=
[a_n,\cdots,a_1].$ We call $\ov w$ the reversal of $w.$
\item[(iv)]
Given words $x=[a_1,\cdots,a_k],\ y=[b_1,\cdots,b_l].$ If $<x>\cap <y>=
\emptyset$ we define a colligation $[x,y]=[a_1,\cdots,a_k,b_1,\cdots,b_l].$ 
\end{itemize}
Given a fixed set $E$ of $n$ distinct positive integers we 
let ${{\boS}}_n$ or ${{\boS}}_E$ denote the set of words $w$ such that the set 
of its entries $<w>=E.$
Let $E=\{a_i\}_{i=1}^n,$ where $a_i<a_{i+\sr 1}$
for all $i\ : 1\leq i<n.$ Taking $i$ instead of $a_i$ 
as a corresponding
entry in the word we get a bijection 
$\phi:{{\boS}}_E\rar {\bS}_n.$
This bijection is 
constructed as a composition of bijections $\phi^{\sr -1}_j$ 
where $\phi_j$ is induced from $\phi_j\ :E \rar E\pr$ defined by
$$\phi_j(a_i):=\begin{cases}
a_i & {\rm if} a_i<j\\
             a_i+1 & {\rm otherwise.}\\
\end{cases}$$
Thus  if $j\not\in E$ we define $\phi^{\sr -1}:E\rar E\pr$ by 
$$\phi_j^{-1}(a_i)=\begin{cases}
a_i & {\rm if}\ a_i<j\\
a_i-1 & {\rm if}\ a_i>j.\\
\end{cases}$$
\par
For $m<n$ we consider ${\bS}_m$ as a subgroup of ${\bS}_n$ by the following
identification
$${\bS}_m=\{w\in {\bS}_n\ :\ a_i=w(i)=i\ \forall i>m\}$$
\par
Set 
$$\begin{array}{rcl}

s_{i,j}^{\sr <}&:=&
\begin{cases}
s_is_{i+\sr 1}\cdots s_j & {\rm if}\ 1\leq i\leq j\\
                    \Id & {\rm otherwise}\\
\end{cases}\\
   s_{i,j}^{\sr >}&:
          =&\begin{cases}
s_is_{i-\sr 1}\cdots s_j & {\rm if}\ i\geq j\geq 1\\
                    \Id & {\rm otherwise}\\
\end{cases}\\
\end{array}
$$
\subsubsection{}\label{2.2.6}
One has
\begin{lemma}
\begin{itemize}
\item[(i)] Given $w \in {\bS}_n.$ Let $p_w(n)=i.$
  Then  $w=ys_{n-{\sr 1},i}^{\sr >} $ where  $y=(w-n),\ y\in {\bS}_{n-1}.$ 
In this $\ell(w) = n-i+\ell(y)$ and 
$$S(w)=S(y)\bigcup\{\al_{a_k,n}\}_{k=i+1}^n.$$
\item[(ii)] Given $w \in {\bS}_n.$ Let $w(n)=i.$
Then $w=s_{i,n-\sr 1}^{\sr <}y$ where $y=\phi^{-1}_i(w_{<{\sr 1},n-\sr 1>}),\ y\in {\bS}_{n-1}.$
In this $\ell(w) = n-i+\ell(y)$ and      
$$S(w)=\{\al_{\phi_i(j),\phi_i(k)}\, :\, \al_{j,k}\in S(y)\}\bigcup\{a_{i,j}\}_{j=i+1}^n.$$
\end{itemize}
\end{lemma}
\Pf
Applying remark \ref{2.2.2} consequently to multiplication of $y=[a_1,\ldots,a_{i-\sr 1}, a_{i+\sr 1},\ldots, a_n]$
by $s_{n-{\sr 1},i}^{\sr >}$ we get
$$\begin{array}{rcl} 
ys_{n-{\sr 1},i}^{\sr >}=[a_1,\ldots, a_n,n]s_{n-\sr 1}s_{n-{\sr 2},i}^{\sr >}=
[a_1,\ldots,n, a_n]s_{n-{\sr 2},i}^{\sr >}&=&\cdots\\
&=&[a_1,\ldots,a_{i-\sr 1},n, a_{i+\sr 1},\ldots, a_n].\\
\end{array}$$
By \ref{2.1.8} one has $\ell(w) = n-i+\ell(y)$ and by \ref{2.2.4} we get  
$$S(w)=S(y)\bigcup\{\al_{j,n}\}_{p_w(j)>p_w(n)}=S(y)\bigcup \{\al_{a_k,n}\}_{k=i+1}^n.$$
\par
(ii) is obtained by applying (i) to $w^{\sr -1}$ and again by applying  \ref{2.2.4}.
for the computation of $S(w).$
\QED
\subsubsection{}\label{2.2.7}
There are two obvious subgroups of ${\bS}_n$ isomorphic to 
${\bS}_{n-1}:$
${\bS}_{n-1}=W_J$ where 
$J=\{\al_i\}_{i=1}^{n-2}$ and ${\bS}\pr_{n-1}=W_{J\pr}$ 
where $J\pr=\{\al_i\}_{i=2}^{n-1}.$ 
In the notation of \ref{2.2.5} one has
${\bS}\pr_{n-1}=\phi_1(\bS_{n-1})$ 
where we add $1$ to the first place in all
the words of $\phi_1({\bS}_{n-1}).$ Note that we get in such a way 
that $\phi_{\sr 1}(s_i)=s_{i+\sr 1}$ and 
$\phi_{\sr 1}(\al_i)=\al_{i+\sr 1}.$
Lemma \ref{2.2.6} can be reformulated for ${\bS}\pr_{n-1}$ as follows
\begin{lemma}\begin{itemize} 
\item[(i)] Given $w \in {\bS}_n.$  Let $p_w(1)=i$
Then $w=ys_{{\sr 1},i-{\sr 1}}^{\sr <} $ where  $y=(w-1),\ y\in {\bS}\pr_{n-1}.$ 
In this $\ell(w) = i-1+\ell(y)$ and 
$S(w)=S(y)\bigcup\{\al_{{\sr 1},a_k}\}_{k=1}^{i-1}.$
\item[(ii)] Given $w \in {\bS}_n.$ Let $w(1)=i.$
Then $w=s_{i-\sr 1,1}^{\sr >}y$ where $y=\phi^{\sr -1}_{i+\sr 1} (\phi_{\sr 1}(w_{<{\sr 2},n>})).$
In this $\ell(w) = i-1+\ell(y)$ and
$S(w)=\{\al_{\phi_{i+\sr 1}(j),\phi_{i+\sr 1}(k)}\ :\ \al_{j,k}\in S(y)\}\bigcup\{a_{j,i+\sr 1}\}_{j=1}^i.$
\end{itemize}
\end{lemma}
\subsubsection{}\label{2.2.8}
We will also need.
\begin{lemma} 
Take $y_1,y_2\in {\bS}_{n-1}$ and $i\ :\ 1\leq i\leq n-1.$
If $y_1 \dor y_2$ then $s_{i,n-1}^{\sr <}y_1\dor s_{i,n-1}^{\sr <}y_2.$
\end{lemma}
\Pf
This follows from the addition of lengths criterion in \ref{2.1.8} 
using \ref{2.2.6} (ii).
\QED
\subsection{\bf Young Diagrams and nilpotent orbits}

\subsubsection{}\label{2.3.1}
Let us define Young diagrams corresponding to the partitions.
Let\ $\lambda=(~\lambda_1~\geq ~\lambda_2 ~\geq \cdots ~\geq ~\lambda_k~ >~ 0~)$
be a partition of $n$. Set
$\lambda^*:=\{\lambda^*_1\geq \lambda^*_2\geq \cdots\geq \lambda^*_l>0\}$
to be the dual partition, that is $\lambda^*_i=
\sharp\{j\ |\ \lambda_j\geq i\}.$ For example $k=\lambda^*_1.$
\par
We define the corresponding Young
diagram $D_\lambda$  of $\lambda$  to be an array of $k$  rows of boxes
starting on the left  with the $i$-th row containing $\lambda_i$ boxes.
\par
For example given $\lambda=(4,3,1)$ then $\lambda^*=(3,2,2,1)$ and
 the corresponding Young tableau is
$$D_\lambda=
\vcenter{
\halign{& \hfill#\hfill
\tabskip4pt\cr
\multispan{9}{\hrulefill}\cr
\ssa
\vb & \ &\vb  &  \ &\vb &  \ &\vb& \ &\ts\vb\cr
\vsa
\multispan{9}{\hrulefill}\cr
\ssa
\vb & \ &\vb & \ & \vb &\ &  \ts\vb\cr
\vsa
\multispan{7}{\hrulefill}\cr
\ssa
\vb & \ & \ts\vb\cr
\vsa
\multispan{3}{\hrulefill}\cr}}$$
The set of Young diagrams with $n$ boxes is denoted by ${\bD}_n$.
\par
Recall order relation on nilpotent orbits
defined in \ref{2.1.1}. In case $\gog=\gs\gl_n$ the  nilpotent 
orbits and  the above order relation on them have 
nice and simple combinatorial descriptions. 
In this case ${\bG=SL}_n$ acts on $\gog$ by conjugation.
For $u\in \gog$ its $\bG$ orbit is 
determined uniquely by the Jordan form of $u.$
If $u\in \Nscr$ all the eigenvalues of $u$ are zero
and the Jordan form of $u$ is determined only
by the length of its Jordan blocks.   
Let us write the Jordan blocks of $u$
in decreasing order and denote
the length of the $i$-th block by $\lambda_i$ .
The resulting partition
 $\lambda = (\lambda_{\sr 1},\cdots,\lambda_k)$  is denoted by $J(u).$  
For example 
$$   u=\left(\begin{array}{cccccccc}
0 & 1 & 0 & 0 & 0 & 0 & 0 & 0 \\
0 & 0 & 1 & 0 & 0 & 0 & 0 & 0 \\
0 & 0 & 0 & 1 & 0 & 0 & 0 & 0 \\
0 & 0 & 0 & 0 & 0 & 0 & 0 & 0 \\
0 & 0 & 0 & 0 & 0 & 1 & 0 & 0 \\
0 & 0 & 0 & 0 & 0 & 0 & 0 & 0 \\
0 & 0 & 0 & 0 & 0 & 0 & 0 & 0 \\
0 & 0 & 0 & 0 & 0 & 0 & 0 & 0 \\
\end{array}\right),
\quad J(u)=(4,2,1,1)$$
The map $u \mapsto J(u)$ gives a bijection of 
$\Nscr/G$  onto ${\bD}_n$. Given $J(u)=\lambda$  we also write 
${\Oscr}_{\lambda}:={\Oscr}_u$ and $D_u:=D_{\lambda}.$
\par
\subsubsection{}\label{2.3.2}
Define  an order relation on Young diagrams as 
follows. Let $\lambda=(\lambda_{\sr 1},\cdots,\lambda_k)$ and 
$\mu = (\mu_{\sr 1},\cdots,\mu_j) $
be partitions of $n$  with corresponding  diagrams
$ D_\lambda, \ D_\mu \in {\bD}_n$. If $j\ne k$ complete the partition
with the lesser number of parts by adding the appropriate number
of $0$'s. In this manner we can consider that both
partitions have $\max(j,k)$ elements.
Define
$D_\lambda \geq D_\mu$ if for each
$i\ :\ 1\leq i\leq \max(j,k)$ one has
$$ \sum^i_{m=1}\lambda_m \leq \sum^i_{m=1}\mu_m \ . $$
\par
The following
result of Gerstenhaber (see ~\cite[\S 3.10]{H}  for example)
describes the closure of a nilpotent orbit.
\begin{theorem} 
Let $\mu$ be a partition of $n$
and ${\Oscr}_\mu$ be the corresponding nilpotent orbit in
$\gs\gl_n.$ One has
$$ \ov{\Oscr}_\mu = \coprod_{\lambda\vert D_\lambda\geq D_\mu}
                                                        {\Oscr}_\lambda \ . $$
In particular 
$$ {\Oscr}_\lambda \geq {\Oscr}_\mu \Llrar D_\lambda\geq D_\mu \ . $$
\end{theorem}
This describes the order relation on nilpotent orbits through the combinatorics of
Young diagrams.
\subsubsection{}\label{2.3.3}
From theorem \ref{2.3.2} we easily obtain
\begin{cor} 
Let ${\Oscr}_\lambda$ be a  descendant 
of ${\Oscr}_\mu.$
Let $D_\lambda, \ D_\mu$ be the corresponding 
Young diagrams. \ Then $D
_\lambda$ is obtained from $D_\mu$ in one of two ways:
\begin{itemize}
\item[(i)] There exists $i$ such that $\mu_i-\mu_{i+1}\geq 2$. \ Then
$\lambda_j = \mu_j$ for $j\neq i, \ i+1$  and
$$ \lambda_i = \mu_i-1, \ \ \lambda_{i+1} = \mu_{i+1}+1 \ . $$
\item[(ii)] There exists $i$ such that
            $ \mu_{i+1}=\mu_{i+2}=\cdots=\mu_{i+k}= \mu_i-1$ 
for  some  $k\geq 1$  and
$\mu_{i+k+1} = \mu_i-2.$ 
Then $\lambda_j = \mu_j$ for $j\neq i, \ i+k+1$ \ and 
$$ \lambda_i = \mu_i-1, \ \ \lambda_{i+k+1} = \mu_i-1 \ . $$
\end{itemize}
\end{cor}
The above result can be described pictorially as follows.
\par
In the first case $D_\lambda$ is obtained from $D_\mu$ by pushing
one box down  one row (and possible across several columns). For example
$$ D_\mu =
\vcenter{
\halign{& \hfill#\hfill
\tabskip4pt\cr
\multispan{9}{\hrulefill}\cr
\ssa
\vb & \quad & \vb &  \quad & \vb  & \quad & \vb & \quad &\ts\vb\cr
\vsa
\multispan{9}{\hrulefill}\cr
\ssa
\vb & \ \ & \vb &  \ \ & \vb  & \ \  & \vb &X&\ts\vb\cr
\vsa
\multispan{9}{\hrulefill}\cr
\ssa
\vb & \ & \ts\vb \cr
\vsa
\multispan{3}{\hrulefill}\cr
\ssa
\vb  & \ &\ts\vb\cr
\vsa
\multispan{3}{\hrulefill}\cr}}\ ,
\qquad
 D_\lambda =
\vcenter{
\halign{& \hfill#\hfill
\tabskip4pt\cr
\multispan{9}{\hrulefill}\cr
\ssa
\vb & \quad & \vb &  \quad & \vb  & \quad & \vb & \quad &\ts\vb\cr
\vsa
\multispan{9}{\hrulefill}\cr
\ssa
\vb & \ & \vb &  \ &  \vb & \ &\ts\vb\cr
\vsa
\multispan{7}{\hrulefill}\cr
\ssa
\vb & \ & \vb &X&\ts\vb\cr
\vsa
\multispan{5}{\hrulefill}\cr
\ssa
\vb  & \ &\ts\vb\cr
\vsa
\multispan{3}{\hrulefill}\cr}}\ . $$
\par
In the second case diagram $D_\lambda$ is obtained from $D_\mu$ by
pushing one box across one column(and possible down several rows).
For example
$$ D_\mu =
\vcenter{
\halign{& \hfill#\hfill
\tabskip4pt\cr
\multispan{9}{\hrulefill}\cr
\ssa
\vb & \quad & \vb &  \quad & \vb  & \quad & \vb &X &\ts\vb\cr
\vsa
\multispan{9}{\hrulefill}\cr
\ssa
\vb & \ & \vb & \ & \vb  & \ &\ts\vb\cr
\vsa
\multispan{7}{\hrulefill}\cr
\ssa
\vb & \ & \vb & \ & \ts\vb \cr
\vsa
\multispan{5}{\hrulefill}\cr
\ssa
\vb  & \ &\ts\vb\cr
\vsa
\multispan{3}{\hrulefill}\cr}}\ ,
\qquad
 D_\lambda =
\vcenter{
\halign{& \hfill#\hfill
\tabskip4pt\cr
\multispan{7}{\hrulefill}\cr
\ssa
\vb & \quad & \vb &  \quad & \vb  & \quad & \ts\vb\cr
\vsa
\multispan{7}{\hrulefill}\cr
\ssa
\vb & \ & \vb &  \ &  \vb & \ & \ts\vb\cr
\vsa
\multispan{7}{\hrulefill}\cr
\ssa
\vb & \ & \vb &  \ &  \vb &X&  \ts\vb\cr
\vsa
\multispan{7}{\hrulefill}\cr
\ssa
\vb  & \ &\ts\vb\cr
\vsa
\multispan{3}{\hrulefill}\cr}}\ .$$
\par
In these cases we say that $D_\lambda$ is a descendant of $D_\mu$.
\subsection{\bf  Young Tableaux and orbital varieties}
\parno
\subsubsection{}\label{2.4.1}
Fill the boxes of Young diagram
$D_\lambda$ with $n$ distinct positive
integers. \ If the entries increase in rows from left to right and in
columns from top to bottom we call such an array a Young tableau or
simply a tableau.
If the numbers in a
Young tableau form the set of integers from 1 to $n$,
then the tableau is called a standard Young tableau.
For example
$$ T =
\vcenter{
\halign{& \hfill#\hfill
\tabskip4pt\cr
\multispan{9}{\hrulefill}\cr
\ssa
\vb & 1 &  &  2 & &  4 & &  8 & \ts\vb \cr
\vsa
&&&&&&\multispan{3}{\hrulefill}\cr
\ssa
\vb & 3 &  & 5  & & 9 & \ts\vb \cr
\vsa
&&\multispan{5}{\hrulefill}\cr
\ssa
\vb & 6 & \ts\vb \cr
\vsa
\multispan{0}{\hrulefill}\cr
\ssa
\vb & 7 & \ts\vb \cr
\vsa
\multispan{3}{\hrulefill}\cr}} $$
is a standard Young tableau.
\par
Let ${\bT}_n$ denote the set 
of all standard Young tableaux of size $n.$
The shape of a Young tableau $T$  is defined to be the 
Young diagram, denoted  $\sh(T)$ , from which it was built.
\subsubsection{}\label{2.4.2}
We will use the following notation for Young tableaux.
Let $T$ be a Young tableau and let $T^i_j$ for $i,j\in \Na$ denote the entry on the intersection of
$i$-th row and $j$-th column. It is sometimes convenient (for example in defining the 
Robinson - Schensted insertion) to assume that each
row and
column is completed to semi-infinite
length by insertion of $\infty$ i.e. if the length
of $i$-th row or $j$-th column in the corresponding 
diagram is $k$ then we set
$$T^i_{k+1}=T^i_{k+2}=\cdots=\infty \ ,\qquad 
  T^{k+1}_j=T^{k+2}_j=\cdots=\infty $$
Given $u$ an entry of $T$ set $r_{\sr T}(u)$ to be the
number of the row, $u$ belongs to and $c_{\sr T}(u)$ to
be the number of the column, $u$ belongs to.  That is if $u=T_j^i$ then 
$r_{\sr T}(u)=i$ and $c_{\sr T}(u)=j.$ The hook number 
of the $ij-$th entry of $T$ is 
defined by $h(T^i_j):= 1+(\lambda_j^*-j)+(\lambda_i-i).$
\par
Let $T^i$ (resp. $T_i$) denote the $i$-th row (resp. column) of $T,$ that is
$$ T^i : = (T^i_1,\ldots),\qquad
T_i :=\left( \begin{array}{c}
T^1_i \\ 
\vdots \\ 
\end{array}  \right), $$
We let $|T^i|$ (resp. $|T_i|$) denote the number of  finite elements in the row 
$T^i$ (resp. column $T_i$) and  $\om^i(T)$ (resp. $\om_i(T)$)  denote the largest finite entry 
of $T^i$ (resp. $T_i$). 
\par

We  consider a tableau as a  matrix $T:=(T_i^j)$  
and  write $T$ by rows or by columns :
$$ T = \left( \begin{array}{c}
T^1 \\ 
\vdots \\ 
T^m \\
\end{array} \right)=
       (T_1, \cdots, T_l) \ . $$
\par
We set $T^{i,j},\ i<j$ to be the subtableau of $T$ consisting of rows from
$i$ to $j,$ that is 
$$ T^{i,j} = \left( \begin{array}{c}
T^i \\ 
\vdots \\ 
T^j \cr
\end{array} \right)\ , $$
and $T^{i,\infty}$ to be the subtableau of $T$ consisting of all $T^j$ such that $j\geq i.$
\par
We set $T_{i,j},\ i<j $ to be the subtableau of $T$ consisting of columns from $i$ 
to $j,$ that is
$$ T_{i,j}=(T_i,\cdots, T_j) $$
and $T_{i,\infty}$ to be a subtableau consisting of all $T_j$ such that $j\geq i.$
\par
Writing $T=T_{1,l}^{1,k}$ designates that $T$ has $l$ columns 
and $k$ rows.
\par
For each tableau $T$ let $T^{\dagger}$ denote the transposed tableau, that is
$${\rm if}\quad T = \left( \begin{array}{c} 
T^1 \\ 
\vdots \\ 
T^m \\
\end{array} \right)=
       (T_1, \cdots, T_l), \quad {\rm then}\quad  T^{\dagger}= (T^1,\cdots, T^m)=
        \left( \begin{array}{c}
T_1 \\ 
\vdots \\ 
T_l \\
\end{array}\right) $$
Note that $\sh(T^\dagger)=\sh(T)^*.$
\par
Given a fixed set $E$ of $n$ distinct positive integers, let 
${\boT}_n$ or ${\boT}_E$ be the set of Young tableaux $T$ such that its set 
of entries  $<T>=E.$
Then ${\boT}_E$ is in bijection with ${\bT}_n$ exactly
the same way as ${{\mathbb S}}_E$ with ${\bS}_n.$ 
We define $\phi_j,\ \phi_j^{-1},\ \phi$  to be the maps induced on
${\bT}_n$ (or ${\boT}_E$) by their action on ${\bS}_n$ (or ${{\boS}}_E$) defined in \ref{2.2.5}. 
\par
To each row $T^i=(T^i_1,\cdots,T^i_k,\cdots)$ with $|T^i|=k$ and column
$T_j=(T_j^1,\cdots,T_j^m,\cdots)$ with $|T_j|=m$ we associate
words, denoted by $[T^i]$ and $[T_j],$ given by 
$$[T^i]=[T^i_1, \cdots , T^i_k]\ ,\qquad
              [T_j]=[T_j^1, \cdots, T_j^m]\ . $$
For simplicity of notation we will omit square brackets inside a colligation. 
For example $[T^i,T_j]$ means $[[T^i],[T_j]].$
\subsubsection{}\label{2.4.3}
A row (resp. column) of $T$ is determined by the set of its entries since
these must increase from left to right (resp. from top to bottom).
\par
Let $T=(T_1,T_2,\ldots,T_l),\ S=(S_1,S_2,\ldots,S_{l\pr})$
be Young tableaux given by their columns. Assume that $T,S$ have no
common entries. Then we define $(T,S)$ to be the array
whose rows  are the same as the rows
of $(T_1,T_2,\ldots,T_l,S_1,S_2,\ldots,S_{l\pr}),$ that is
$r_{\sr(T,S)}(T^i_j)=i$ and $r_{\sr (T,S)}(S_j^i)=i,$ and
ordered in the increasing order.
Of course this involves the shuffling of numbered boxes within a 
row.  
\begin{lemma} [{~\cite[2.7]{JM}}] 
$(T,S)$ is a Young tableau. 
\end{lemma}
If the entries 
of $S$ all exceed those of $T$ then
one only needs to shift numbered boxes (to the left).
\par
Taking $T$ and $S$ by rows instead of columns, that is
$$T=\left( \begin{array}{c} 
T^1 \\ 
\vdots \\ T^k \\
\end{array} \right),\quad 
S=\left(\begin{array}{c}
S^1\cr\vdots\cr S^{k\pr}\cr\end{array}\right)$$
we define $\left(\begin{array}{c}T\\ S\\ \end{array}\right)$ in a  fashion similar  to $(T,S).$ One has 
$\left(\begin{array}{c}T\\ S\\ \end{array}\right)^{\dagger}=(T^\dagger,S^\dagger).$
\par
Note that $T^{i,i+1}=\left(\begin{array}{c}T^i\\ T^{i+1}\\ \end{array}\right),$ for example.
\subsubsection{}\label{2.4.4}
Given $D_\lambda \in {\bD}_n$ with $\lambda=(\lambda_1,\cdots,\lambda_j)$
we define a corner box (or simply, a corner) of the Young diagram 
to be a box with
no neighbours to right and below.
\par
For example in $D$ below all the corner boxes are labeled by $X$.
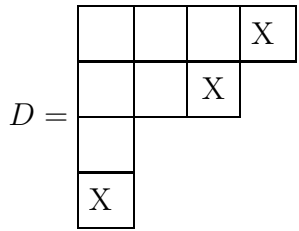
$$ D =
\vcenter{
\halign{& \hfill#\hfill
\tabskip4pt\cr
\multispan{9}{\hrulefill}\cr
\ssa
\vb & \quad & \vb &  \quad & \vb  & \quad & \vb & X &\ts\vb\cr
\vsa
\multispan{9}{\hrulefill}\cr
\ssa
\vb & \ & \vb &  \ & \vb  &X&\ts\vb\cr
\vsa
\multispan{7}{\hrulefill}\cr
\ssa
\vb & \ & \ts\vb \cr
\vsa
\multispan{3}{\hrulefill}\cr
\ssa
\vb  &X &\ts\vb\cr
\vsa
\multispan{3}{\hrulefill}\cr}} $$
The entry of a Young tableau in a corner is called a corner entry.
Take  $D_\lambda$ with  $\lambda=(\lambda_1,\cdots,\lambda_k).$ Then  there is a 
corner  
entry $\om^i(T)$ at the corner  $c=c(i,\lambda_i)$ with coordinates $(i,\lambda_i)$ iff 
$\lambda_{i+1}<\lambda_i.$
We order corners by ordering their first coordinate, i.e.
for $c=c(i,\lambda_i),\ c\pr=c(j,\lambda_j)$ one has $c<c\pr$ iff $i<j.$   
\subsubsection{}\label{2.4.5} 
We now define the insertion algorithm. 
Consider a row  $R = (a_1,a_2 , \cdots )$ completed by $\infty.$
Given $j \in {\Na}^+, \ \ j \not \in <R>$ or $j=\infty.$ Let $a_i$  be the first smallest
integer greater or equal to $j$ (possibly $\infty$) in $<R>$  and set
$$(R \dar j):= \begin{cases}
(a_1,\cdots,a_{i-1},j,a_{i+1},\cdots) \ ,\ \ j_{\sr R}=a_i & {\rm if} j\ne\infty\\
                      R & {\rm otherwise}\\
\end{cases}$$
The inductive extension of this operation 
 to a Young tableau
$T$  with $m$  rows for $j\not\in <T>$ given by
$$ (T \Da j)=\left(\begin{array}{l}(T^1\dar j) \cr
                  (T^{2,m} \Da j_{_{\sr T^1}}) \cr\end{array}\right) $$
is called the  insertion algorithm.
\par
Note that the shape of $(T\Da j)$ is the shape of $T$ obtained by adding one new corner. 
The entry of this corner is denoted by $j_T.$ 
\par
\subsubsection{}\label{2.4.6}
Let $w=[a_{\sr 1},a_{\sr 2},\ldots,a_n].$ 
According to Robinson-Schensted procedure we associate an ordered pair
of Young tableaux $(T(w),Q(w))$ to $w.$ The procedure is fully explained 
in  ~\cite[5.1.4]{Kn}  or in ~\cite[2.5]{Re}. As well we explain it in \ref{2.4.18}. In what follows 
we call it RS procedure. Here we explain only the inductive procedure
of constructing the first tableau $T(w).$ 
\parno
\begin{itemize}
\item[(1)] Set $ {}_1T(w) = (a_1).$
\item[(2)] Set ${}_{j+1}T(w)=({}_jT(w) \Da a_{j+1}) \ .$
\item[(3)] Set $T(w)={}_nT(w) \ .$
\end{itemize}
\par
As an example we take
$$ w =  [ 2 , 5,  1, 4, 3 ] $$
This gives
$${}_1T(w) =
    \vcenter{
\halign{& \hfill#\hfill
\tabskip4pt\cr
\multispan{3}{\hrulefill}\cr
\ssa
\vb & 2  & \ts\vb\cr
\vsa
\multispan{3}{\hrulefill}\cr}} \quad
{}_2T(w) =
\vcenter{
\halign{& \hfill#\hfill
\tabskip4pt\cr
\multispan{5}{\hrulefill}\cr
\ssa
\vb &2 & & 5 & \ts\vb\cr
\vsa
\multispan{5}{\hrulefill}\cr}}\quad
{}_3T(w) =
\vcenter{
\halign{& \hfill#\hfill
\tabskip4pt\cr
\multispan{5}{\hrulefill}\cr
\ssa
\vb&1 & & 5  & \ts\vb\cr
\vsa
&&\multispan{3}{\hrulefill}\cr
\ssa
\vb & 2 &  \ts\vb\cr
\vsa
\multispan{3}{\hrulefill}\cr}}$$
$${}_4T(w) =
\vcenter{
\halign{& \hfill#\hfill
\tabskip4pt\cr
\multispan{5}{\hrulefill}\cr
\ssa
\vb & 1 && 4  & \ts\vb\cr
\vsa
&&&&\multispan{0}{\hrulefill}\cr
\ssa
\vb&2 && 5 & \ts\vb\cr
\vsa
\multispan{5}{\hrulefill}\cr}}\qquad 
{}_5T(w)=
\vcenter{
\halign{& \hfill#\hfill
\tabskip4pt\cr
\multispan{5}{\hrulefill}\cr
\ssa
\vb&1 &&3  & \ts\vb\cr
\vsa
&&&&\multispan{0}{\hrulefill}\cr
\ssa
\vb& 2 & &4 & \ts\vb\cr
\vsa
&&\multispan{3}{\hrulefill}\cr
\ssa
\vb&5 & \ts\vb\cr
\vsa
\multispan{3}{\hrulefill}\cr}} $$
\par
The result due to Robinson and Schensted (see \ref{2.4.18})  implies the map
$\varphi: w \mapsto T(w)$ is a surjection of  ${\bS}_n$ (resp. ${{\boS}}_E$) onto 
${\bT}_n$ (resp. ${\boT}_E$). 
\par
Recall the notion of geometric cells from \ref{2.1.4}.
By R. Steinberg ~\cite{St1} the fibres of $\varphi$ provide a partition of ${\bS}_n$  into
geometric cells, that is
\begin{theorem} 
For any $w\in{\bS}_n$ one has ${\Cscr}_w=\{y\in {\bS}_n\ :\ T(y)=T(w)\}.$ 
\end{theorem}

We generalize this description to ${{\boS}}_E$
by taking ${\Cscr}_w:=\{y\in  {{\boS}}_E\ :\ T(y)=T(w)\}.$
As well we define 
${\Cscr}_T:=\{y\in {\bS}_n\ ({\rm or\ } {{\boS}}_E)
\ :\ T(y)=T\}.$
We also denote by $T_{\Cscr}$ the tableau corresponding to 
the cell $\Cscr.$  
\subsubsection{}\label{2.4.7} 
RS procedure together with 
lemma \ref{2.2.6} (ii)  implies
\begin{prop} 
For $y\in {\bS}_n$ one has
$$T(s_{i,n}^{\sr <}y)=(\phi_i(T(y))\Da i)\qquad where\quad 1\leq i\leq n+1 $$
\end{prop}
\Pf
Set $w=s_{i,n}^{\sr<}y.$ By \ref{2.2.6} (ii) $w(n+1)=i$ and $y=\phi_i^{-1}(w_{<{\sr 1},n>}).$
Hence by RS procedure 
$$T(w)=({}_nT(w)\Da i)=(T(\phi_i(y))\Da i)
=(\phi_i(T(y))\Da i).$$
\QED 
\subsubsection{}\label{2.4.8}
Let us describe a few algorithms connected to RS
procedure which we use for proofs and constructions. 
\par
First let us describe some operations for rows and tableaux.
Consider a row $R=(a_1,\cdots).$ 
\begin{itemize}
\item[(i)] Set $(R-a_i):=(a_1,\cdots,a_{i-1},a_{i+1},\cdots).$
\item[(ii)] For $j \in {\Na},\ j\not\in <R>$ let $a_i$ be the greatest
element of $<R>$ smaller than $j$ and set 
$$(R+j):
=(a_1,\cdots,a_i,j,a_{i+1},\cdots).$$
\item[(iii)] We define a pushing up operation. Again let 
$j \in {\Na},\ j\not\in <R>$ and $j>a_1.$ 
Let $a_i$
be the greatest entry of $R$ smaller than $j$ and set :
$$(R \uar j):=(a_1,\cdots,a_{i-1},j,a_{i+1},\cdots),\qquad j^{\sr R}:
                                                                     =a_i.$$
The last operation is extended to a Young tableau $T$ by induction
on the number of rows. Let $T^m$ be the last row of $T$ and assume $T^m_1<j.$ Then
$$(T\uar j)=\left(\begin{array}{l}(T^{1,m-1} \uar j^{\sr T^m})\cr
                  (T^m \uar j) \cr\end{array}\right)$$
We denote by $j^{\sr T}$ the element pushed out 
from the first row of the tableau in the last step.
\end{itemize}
\subsubsection{}\label{2.4.9}
The pushing up operation gives us a procedure of deleting a corner inverse to the
insertion algorithm. This is also described in ~\cite[5.1.4]{Kn}  or 
~\cite[2.8]{Re}. 
\par
 As a result of insertion we get a new tableau with a
shape obtained from the old one just by adding one corner. As a result of 
deletion we  get a new tableau with the shape obtained from the old one
by extracting one corner.
\par
Let $T$ be a tableau. Recall the definition  of $\om^i(T)$ from \ref{2.4.2}. 
Assume $\lambda_i>\lambda_{i+1}$ and let $c=c(i,\lambda_i)$ be the corner of $T$  on the $i$-th row.
To delete the  corner $c$ we
delete $\om^i(T)$ from the row $T^i$ and
push it up through the tableau $T^{1,i-1}.$ 
The position of the displaced boxes can be joined to form a segment $s_c(T)$
in $\sh(T).$
The element pushed out from
the tableau is denoted by $c^{\sr T}.$ This is written
$$(T\Ua c):=\left(\begin{array}{l}(T^{1,i-1} \uar \om^i(T))\cr
                          (T^i-\om^i(T))\cr
                          T^{i+1,\infty}\cr
\end{array}\right)$$
For example
$$\left (\ \vcenter{
\halign{& \hfill#\hfill
\tabskip4pt\cr
\multispan{5}{\hrulefill}\cr
\ssa
\vb&1 &&3  & \ts\vb\cr
\vsa
&&&&\multispan{0}{\hrulefill}\cr
\ssa
\vb& 2 & &4 & \ts\vb\cr
\vsa
&&\multispan{3}{\hrulefill}\cr
\ssa
\vb&5 & \ts\vb\cr
\vsa
\multispan{3}{\hrulefill}\cr}}\Ua c(3,1)\right )=
\vcenter{
\halign{& \hfill#\hfill
\tabskip4pt\cr
\multispan{5}{\hrulefill}\cr
\ssa
\vb & 1 && 4  & \ts\vb\cr
\vsa
&&&&\multispan{0}{\hrulefill}\cr
\ssa
\vb&2 && 5 & \ts\vb\cr
\vsa
\multispan{5}{\hrulefill}\cr}},\qquad c^{\sr T}=3,\qquad 
s_c(T)=\vcenter{
\halign{& \hfill#\hfill
\tabskip4pt\cr
\multispan{5}{\hrulefill}\cr
\ssa
\vb& \ &\vb&\ $\uparrow$  & \ts\vb\cr
\vsa
\multispan{5}{\hrulefill}\cr
\ssa
\vb& \ &\vb &\ $\uar$ & \ts\vb\cr
\vsa
\multispan{5}{\hrulefill}\cr
\ssa
\vb& $\nearrow$& \ts\vb\cr
\vsa
\multispan{3}{\hrulefill}\cr}}.$$ 
\par
Note that insertion and deletion are  indeed inverse since 
 for any $T\in {\bT}_n$ 
\bk
$$((T\Ua c)\Da c^{\sr T})=T\quad {\rm and} 
\quad((T\Da j)\Ua j_{\sr T})=T\quad ({\rm for}\ j\not\in <T>).\eqno{(*)}$$
Note that sometimes we will write $(T\Ua a)$ where $a$ is a corner entry 
just as we have written above.
\subsubsection{}\label{2.4.10}
Inserting, deleting, adding and  pushing up operations are defined for rows. 
For our further analysis we need to define these operations also for columns. 
The most simple way to define them is by the operation of 
transposition ($\dagger$). 
Given a column $C;$ given $a\in<C>$ and $j\in \Na,\ j\not\in <C>,$ or $j=\infty$   set
\begin{itemize}
\item[(i)] $(C-a):=(C^\dagger-a)^\dagger;$
\item[(ii)] $(C+j):=(C^\dagger+j)^\dagger;$
\item[(iii)] $(j\rar C):=(C^\dagger\dar j)^\dagger$ and ${}_{\sr C}j:=j_{\sr C^\dagger};$
\item[(iv)] if $j>C^1$ then $(C\lar j):=(C^\dagger\uar j)^\dagger$ 
and ${}^{\sr C}j:=j^{\sr C^\dagger};$
\end{itemize}
We also extend the operations to tableaux:
$$(j\Ra T):=(T^{\dagger}\Da j)^{\dagger}\ ,\quad {}_{\sr T}j:=j_{\sr T^{\dagger}};\qquad\qquad 
(T\La c):=(T^{\dagger}\Ua c)^{\dagger}\ ,\quad {}^{\sr T}c:=c^{\sr T^{\dagger}}\ .\eqno{(*)}$$

\subsubsection{}\label{2.4.11}
Recall the notion of $h(T^i_j)$ from \ref{2.4.2}. Note that if 
$h(T^i_j)>1$ and $T^{i+1}_j$ (resp.$T^i_{j+1}$) is defined that is $T^{i+1}_j\ne \infty$ 
(resp. $T^i_{j+1}\ne \infty$) then $h(T^{i+1}_j)<h(T^i_j)$
(resp. $T^i_{j+1}<h(T^i_j)$).
\par
Let us describe the jeu de taquin (see ~\cite{Sch}) which removes $T^i_j$
from $T.$ The resulting tableau is denoted by $T-T^i_j$ and is 
obtained by the following  procedure recursive on $h(T^i_j).$
\begin{itemize}
\item[(1)] If $ h(T^i_j)=1$ then 
          $$(T- T^i_j)=\left(\begin{array}{c}T^{1,i-1}\cr
                          (T^i- T^i_j)\cr
                             T^{i+1,\infty}\cr\end{array}\right)$$
\item[(2)]  If $h(T^i_j)>1,$ then
\begin{itemize}
\item[(i)] If $T^i_{j+1}> T^{i+1}_j$ or $\lambda_i=j,$ set
$$(T-T^i_j)=\left(\begin{array}{l}T^{1,i-1}\cr
                            (T^i\uar T^{i+1}_j)\cr
                            (T^{i+1,m}-T^{i+1}_j)\cr\end{array}\right)$$
\item[(ii)] If $T^i_{j+1}< T^{i+1}_j$ or $\lambda_j^*=i,$ set
$$(T-T^i_j)=(T_{1,j-1},(T_j\lar T^i_{j+1}),
(T_{j+1,\lambda_{\sr 1}}-T^i_{j+1}))$$
\end{itemize}
\end{itemize}
\parno
The result due to M. P. Sch\"utzenberger ~\cite{Sch} gives
\begin{theorem} 
If $T$ is a Young tableau then $(T- T^i_j)$
is a Young tableau.
\end{theorem}
As an example we take 
$$ T =
\vcenter{
\halign{& \hfill#\hfill
\tabskip4pt\cr
\multispan{7}{\hrulefill}\cr
\ssa
\vb & 1 &  &  2 & &  5 & \ts\vb\cr
\vsa
&&&&\multispan{3}{\hrulefill}\cr
\ssa
\vb & 3 &  &  4 & \ts\vb\cr
\vsa
&&\multispan{3}{\hrulefill}\cr
\ssa
\vb & 6 & \ts\vb\cr
\vsa
\multispan{3}{\hrulefill}\cr}}$$
Then
$$(T-6) =
\vcenter{
\halign{& \hfill#\hfill
\tabskip4pt\cr
\multispan{7}{\hrulefill}\cr
\ssa
\vb & 1 &  &  2 & &  5 & \ts\vb\cr
\vsa
&&&&\multispan{3}{\hrulefill}\cr
\ssa
\vb & 3 &  &  4 & \ts\vb\cr
\vsa
\multispan{5}{\hrulefill}\cr}},\quad
(T-5)=
\vcenter{
\halign{& \hfill#\hfill
\tabskip4pt\cr
\multispan{5}{\hrulefill}\cr
\ssa
\vb & 1 &  &  2 &  \ts\vb\cr
\vsa
&&&&\multispan{0}{\hrulefill}\cr
\ssa
\vb & 3 &  &  4 & \ts\vb\cr
\vsa
&&\multispan{3}{\hrulefill}\cr
\ssa
\vb & 6 & \ts\vb\cr
\vsa
\multispan{3}{\hrulefill}\cr}},\quad
(T-4) =
\vcenter{
\halign{& \hfill#\hfill
\tabskip4pt\cr
\multispan{7}{\hrulefill}\cr
\ssa
\vb & 1 &  &  2 & &  5 & \ts\vb\cr
\vsa
&&\multispan{5}{\hrulefill}\cr
\ssa
\vb & 3 & \ts\vb\cr
\vsa
&&\cr
\ssa
\vb & 6 & \ts\vb\cr
\vsa
\multispan{3}{\hrulefill}\cr}}.$$
$$ (T-3) =
\vcenter{
\halign{& \hfill#\hfill
\tabskip4pt\cr
\multispan{7}{\hrulefill}\cr
\ssa
\vb & 1 &  &  2 & &  5 & \ts\vb\cr
\vsa
&&\multispan{5}{\hrulefill}\cr
\ssa
\vb & 4 & \ts\vb\cr
\vsa
&&\cr
\ssa
\vb & 6 & \ts\vb\cr
\vsa
\multispan{3}{\hrulefill}\cr}},\quad
(T-2) =
\vcenter{
\halign{& \hfill#\hfill
\tabskip4pt\cr
\multispan{7}{\hrulefill}\cr
\ssa
\vb & 1 &  &  4 & &  5 & \ts\vb\cr
\vsa
&&\multispan{5}{\hrulefill}\cr
\ssa
\vb & 3 & \ts\vb\cr
\vsa
&&\cr
\ssa
\vb & 6 & \ts\vb\cr
\vsa
\multispan{3}{\hrulefill}\cr}},\quad
 (T-1) =
\vcenter{
\halign{& \hfill#\hfill
\tabskip4pt\cr
\multispan{7}{\hrulefill}\cr
\ssa
\vb & 2 &  &  4 & &  5 & \ts\vb\cr
\vsa
&&\multispan{5}{\hrulefill}\cr
\ssa
\vb & 3 & \ts\vb\cr
\vsa
&&\cr
\ssa
\vb & 6 & \ts\vb\cr
\vsa
\multispan{3}{\hrulefill}\cr}}.$$
\subsubsection{}\label{2.4.12}
The importance of  removing  the largest entry and the smallest one 
is  provided  by the following  
\begin{prop}[{~\cite[p.~60]{Kn}}]
Consider $w=[a_1,\cdots,a_n]$ and set $M:=
\max<w>.$ Then $T(w-M)=(T(w)-M).$ 
\end{prop}
and by the 
following theorem due to 
M. P. Sch\"utzenberger (see ~\cite[5.1.4]{Kn}).
\begin{theorem} 
Consider $w=[a_1,\cdots,a_n]$ and set $m:=min<w>.$ 
Then $T(w-m)=(T(w)-m).$
\end{theorem}
\subsubsection{}\label{2.4.13}
\ \begin{cor}\begin{itemize}
\item[(i)] Let $w\in {\bS}_n$ have a decomposition
        $w=ys_{n-1,i}^{\sr >}$ where $y \in {\bS}_{n-1}.$ 
        If $T(w)=T$  then  $T(y)=(T-n).$  
\item[(ii)] Let $w\in {\bS}_n$ have a decomposition
        $w=ys_{{\sr 1},i}^{\sr <}$ where 
$y \in {\bS}\pr_{n-1}.$ If $T(w)=T$ then $T(y)=(T-1).$  
\end{itemize}
\end{cor}
\Pf
Indeed, this is straightforward from \ref{2.4.12} and lemmas \ref{2.2.6} (i),
\ref{2.2.7} (i). 
\QED
\subsubsection{}\label{2.4.14}
Recall the notion of $\tau-$invariant from \ref{2.1.7}.
Given $T\in{\bT}_n$ set 
$$\tau(T):=\{\al_i\ : r_{\sr T}(i+1)>r_{\sr T}(i)\}.$$
By \ref{2.2.4} Remark 1, \ref{2.1.7} and \ref{2.4.6} one has  
\begin{lemma} 
For $w\in{\bS}_n$ set $T=T(w).$ Then
$\tau(w)=\tau(T)=\tau({\Vscr}_T)=\tau({\Cscr}_T).$
\end{lemma}
\subsubsection{}\label{2.4.15}
After Schensted - Sch\"utzenberger (see ~\cite[5.4.1]{Kn} ) one has
\begin{theorem} 
$T^{\dagger}(w)=T(\ov w), \ \forall w\in {\bS}_n.$
\end{theorem}
\subsubsection{}\label{2.4.16}
Recall projections $\pi_{\sr \Iscr}$ from \ref{2.1.8}. We consider
these projections for $\gog=\gs\gl_n$ and $\Iscr\st \Pi$ defined as follows.
For $1\leq i< j\leq n$ set 
$$<i,j>:=\{r\in{\Na}\ | i\leq r\leq j\},\
\Pi_{i,j}:=\{\al_r\ |\ r,r+1\in <i,j>\}\ {\rm and}\ \pi_{i,j}:=
\pi_{\Pi_{i,j}}.$$
\par
Define $\pi_{i,j}:{\bT}_n\rar {\boT}_{<i,j>},$ through the jeu de taquin applied
to the entries of $T$ not lying in $<i,j>$ (taken in any order).
Note that the order of elimination is indeed not important. Denote 
$T^{<i,j>}:=\pi_{i,j}(T).$
Set $D^{<i,j>}_T:=\sh(T^{<i,j>}).$ 
\par
We can represent each Young tableau as a chain of Young diagrams by
the following.
Take
$D_\lambda \in {\bD}_n,\ \ D_\mu \in {\bD}_m,\ \ n<m$.
We set $D_\lambda \subset D_\mu$  if $\lambda_i \leq \mu_i,\ \
\forall i.$ \  
Let ${\bC}_n$ denote the set of all decreasing chains
of Young diagrams :
$${\bC}_n:=\{(D_n,\cdots,D_1)\ :\ D_i \in {\bD}_i\ {\rm such\ that}\ D_{i-1}\subset D_i\}$$
We define a map
$\psi:{\bT}_n\longrightarrow {\bC}_n$ by
$\psi(T)=( D^{<1,n>}_T,\ D^{<1,n-1>}_T,\ \cdots, D^{<1,1>}_T).$
\par
For example, if
$$ T =
\vcenter{
\halign{& \hfill#\hfill
\tabskip4pt\cr
\multispan{7}{\hrulefill}\cr
\ssa
\vb & 1 &  &  2 & &  5 & \ts\vb\cr
\vsa
&&&&\multispan{3}{\hrulefill}\cr
\ssa
\vb & 3 &  &  4 & \ts\vb\cr
\vsa
&&\multispan{3}{\hrulefill}\cr
\ssa
\vb & 6 & \ts\vb\cr
\vsa
\multispan{3}{\hrulefill}\cr}}$$
Then
$$\psi(T)= \left(\ 
\vcenter{
\halign{& \hfill#\hfill
\tabskip4pt\cr
\multispan{7}{\hrulefill}\cr
\ssa
\vb & \ &\vb& \ &\vb &  \ & \ts\vb\cr
\vsa
\multispan{7}{\hrulefill}\cr
\ssa
\vb & \ &\vb&  \ & \ts\vb\cr
\vsa
\multispan{5}{\hrulefill}\cr
\ssa
\vb & \ & \ts\vb\cr
\vsa
\multispan{3}{\hrulefill}\cr}},\ 
\vcenter{
\halign{& \hfill#\hfill
\tabskip4pt\cr
\multispan{7}{\hrulefill}\cr
\ssa
\vb & \ &\vb &  \ &\vb &  \ & \ts\vb\cr
\vsa
\multispan{7}{\hrulefill}\cr
\ssa
\vb & \ &\vb&  \ & \ts\vb\cr
\vsa
\multispan{5}{\hrulefill}\cr}},\ 
\vcenter{
\halign{& \hfill#\hfill
\tabskip4pt\cr
\multispan{5}{\hrulefill}\cr
\ssa
\vb & \ & \vb&  \  & \ts\vb\cr
\vsa
\multispan{5}{\hrulefill}\cr
\ssa
\vb & \ &\vb  &  \ & \ts\vb\cr
\vsa
\multispan{5}{\hrulefill}\cr}},\ 
\vcenter{
\halign{& \hfill#\hfill
\tabskip4pt\cr
\multispan{5}{\hrulefill}\cr
\ssa
\vb & \ &\vb &  \  & \ts\vb\cr
\vsa
\multispan{5}{\hrulefill}\cr
\ssa
\vb & \ &  \ts\vb\cr
\vsa
\multispan{3}{\hrulefill}\cr}},\ 
\vcenter{
\halign{& \hfill#\hfill
\tabskip4pt\cr
\multispan{5}{\hrulefill}\cr
\ssa
\vb & \ &\vb  &  \  & \ts\vb\cr
\vsa
\multispan{5}{\hrulefill}\cr}},\
\vcenter{
\halign{& \hfill#\hfill
\tabskip4pt\cr
\multispan{3}{\hrulefill}\cr
\ssa 
\vb & \  & \ts\vb\cr
\vsa
\multispan{3}{\hrulefill}\cr}}\ \right).$$
We  reconstruct $T$ from this chain by inserting $i$ into the box
deleted on going from $D^{<1,i>}_T$  to $D^{<1,i-1>}_T.$  It follows easily
that $\psi$ is a bijection.
\subsubsection{}\label{2.4.17}
The interpretation of $\psi$ from \ref{2.4.16}  given by
Spaltenstein in ~\cite{Sp1} provides us a description of the
connection between orbital varieties and Young tableaux somewhat different from 
Steinberg's construction described in \ref{2.1.3} and \ref{2.4.6}.
\par
Let $0=V_0\subset V_1\subset\cdots\subset V_n={\Co}^n$ be the standard flag
invariant under $\gn.$ 
Consider $x\in \gn$ and let $x_i=x\vert_{V_i}$ be the restriction of 
$x$ to $V_i.$ Then $C=(D_{x_n}, D_{x_{n-\sr 1}},\cdots,D_{x_{\sr 1}})$
is  a decreasing chain of Young diagrams and so defines a 
standard Young tableau $T(x)=\psi^{-1}(C)$ .  
Set $\eta:\gn \rar {\bT}_n$ by $\eta(x)= T(x).$ 
Consider $T\in{\bT}_n$ and 
denote the  fiber of $\eta$ in $\gn$ by 
$$\gn(T)=\{x\in \gn\ :\ \eta(x)=T\}.$$
\par
By ~\cite{Sp1} one has   
\begin{theorem} 
Every $\gn(T)$ is dense in  a unique 
orbital variety $\Vscr\subset \gs\gl_n$, and
every orbital variety is obtained in that way.
\end{theorem}
\subsubsection{}\label{2.4.18}
Recall RS procedure from \ref{2.4.6}. 
Given $w\in {\bS}_n$ let 
$C(w):=(\sh({}_nT(w)),\ldots,\sh({}_1T(w)))$
and set $Q(w):=\psi^{\sr -1}(C(w)).$ 
Define $\theta(w):=(T(w),Q(w)).$ Set 
$$\widetilde {\bT}_n:=\{(T,S)\in {\bT}_n\times{\bT}_n\ :\ \sh(T)=\sh(S)\}$$ 
It was shown by Robinson and Schensted 
(see ~\cite[5.1.4]{Kn}  or in ~\cite[2.5]{Re} )
\begin{theorem} 
$\theta:{\bS}_n\rar \widetilde{\bT}_n$ is a bijection.
Moreover $\theta(w^{\sr -1})=(Q(w),T(w)).$
\end{theorem}
In particular one has $|{\Cscr}_w|$ is the number of standard tableaux
of $\sh(T(w)).$ Thus for any $y\in{\Cscr}^d_w$ one has $|{\Cscr}_y|=|{\Cscr}_w|.$
\subsection{\bf  Duflo descendant and Duflo offspring of a Young cell }

\subsubsection{}\label{2.5.1}
Recall the induced Duflo order from \ref{2.1.6}. 
We can speak about  ordering of cells in  ${{\boS}}_n$
as well as about ordering of cells in ${\bS}_n.$ We set 
$\Cscr\dor \Cscr\pr$ if $\phi(\Cscr) \dor \phi(\Cscr\pr).$
We define the induced Duflo order on tableaux as well by $T\dor S$
if ${\Cscr}_T\dor{\Cscr}_S.$
\par
We define a Duflo descendant
of an orbital variety exactly as we have defined a geometric descendant 
of an orbital variety in \ref{1.8} only with respect to Duflo order.
\par
Given orbital varieties ${\Vscr}_1$ and ${\Vscr}_2$ 
such that ${\Vscr}_1\dos {\Vscr}_2.$ 
\begin{defi} 
We call ${\Vscr}_2$ a Duflo descendant
of ${\Vscr}_1$ if for any orbital variety $\Wscr$ one has,
${\Vscr}_1\dor \Wscr\dor {\Vscr}_2$ implies $\Wscr={\Vscr}_1$ or $\Wscr={\Vscr}_2.$
\end{defi}
In that case we call ${\Cscr}_{{\Vscr}_2}$ (resp. $T_{{\Vscr}_2}$)  
a Duflo descendant of ${\Cscr}_{{\Vscr}_1}$ (resp. $T_{{\Vscr}_1}$).
\subsubsection{}\label{2.5.2}
Given $y\in {\bS}_n$ 
we call $ys_i$ an offspring of $y$ if $y\dos ys_i.$ We induce this definition
to the cells.
\begin{defi}
${\Cscr}_2$ is called a Duflo offspring of ${\Cscr}_1$   
if there exist $y \in {\Cscr}_1$ and $i\ :\ 1\leq i\leq n-1$ such that $y \dos ys_i$ and $ys_i\in {\Cscr}_2.$  
\end{defi}
In that case we call ${\Vscr}_{{\Cscr}_2}$ (resp. $T_{{\Cscr}_2}$)  
a Duflo offspring of ${\Vscr}_{{\Cscr}_1}$ (resp. $T_{{\Cscr}_1}$).
\par
Let $\Dscr(T)$ (resp. $\Dscr(\Cscr)$) denote the set of all offsprings of a given tableau $T$
(resp. of a  given cell $\Cscr$).
\par
Note that each non-trivial cell (i.e. $\Cscr\ne\{[1,2,\cdots,n]\},
\ \{[n,n-1,\cdots,1]\}$) is an offspring of itself. For uniqueness we define
$T([1,2,\cdots,n])\in\Dscr(T([1,2,\cdots,n]))$ and 
$\Dscr(T([n,n-1,\cdots,1])):=\{T([n,n-1,\cdots,1])\}.$
\par
Note also that some  offsprings of $\Cscr$ can be an offspring  
of another offspring of $\Cscr$ as it is shown by the
\parno
{\bf Example.}\ \ \ Consider the following cells in ${\bS}_4$
$$\begin{array}{rcl}
{\Cscr}_1&=&\{s_2,\ s_2 s_1,\  s_2 s_3\},\\
 {\Cscr}_2&=&\{s_2 s_1 s_3,\ s_2 s_1 s_3 s_2,\},\\
 {\Cscr}_3&=&\{s_2 s_1 s_2,\  s_2 s_1 s_2 s_3,\ s_2 s_1 s_2 s_3 s_2,\}.\\
\end{array}$$
Then ${\Cscr}_2,{\Cscr}_3\in\Dscr({\Cscr}_1)$ since $s_2 s_1,\dos s_2 s_1 s_3$ and 
$s_2 s_1\dos s_2 s_1 s_2.$ As well ${\Cscr}_3\in \Dscr({\Cscr}_2)$ since
$ s_2 s_1 s_3 s_2\dos s_2 s_1 s_3 s_2 s_3=s_2 s_1 s_2 s_3 s_2.$ 
\par
Obviously the set of Duflo descendants of $\Cscr$ is a subset of Duflo offsprings
of $\Cscr.$

\section {\bf Description of Induced Duflo Order}

\subsection{\bf Decomposition of a cell}
\parno
\subsubsection{}\label{3.1.1} 
Recall the notation from \ref{2.2.5} and \ref{2.4.4}. 
Consider a standard Young tableau $T \in \bT_n$ with $m$ corners  
$\{c_i\}_{i=1}^m$ where $c_{\sr 1}<c_{\sr 2}<\cdots<c_m.$  
The deletion of a corner $c_i$ by the  algorithm
described in \ref{2.4.9} gives rise to the tableau 
$(T\Ua c_i)$ and the pushed out element $p_i=c_i^{\sr T}$ where $<(T\Ua c_i)>=\{j\}_{j=1,\ j\ne p_i}^n.$
By \ref{2.4.2} there exist a bijection 
$\phi:{\boT}_{n-1}\rar \bT_{n-1}$ such that $\phi(T\Ua c_i)$
is a standard Young tableau  (in our case $\phi=\phi_{p_i}^{-1}$ in notation
of \ref{2.2.5}). In what follows we denote ${\Cscr}_{(T\Ua c_i)}:={\Cscr}_{\phi(T\Ua c_i)}$
just for simplicity of notation. This is a cell in $\bS_{n-1}.$ 
\begin{prop} For each standard Young tableau $T \in \bT_n$ one
has that ${\Cscr}_T$ is a disjoint union:
$${\Cscr}_T=\coprod_{i=1}^m  s^{\sr <}_{p_i,n-\sr 1}{\Cscr}_{(T\Ua c_i)}$$
\end{prop} 
\Pf
Take $w\in {\Cscr}_T.$ Let $w(n)=p.$
Then by  \ref{2.2.6}(ii) $w=s^{\sr <}_{p,n-\sr 1}y$ with 
$y\in \bS_{n- {\sr 1}}.$ 
By RS procedure $T=({}_{n- {\sr 1}}T(w) \Da p).$ Moreover
by procedure \ref{2.4.9} one has
$(T\Ua p_{\sr T})={}_{n- {\sr 1}}T(w)=\phi_p(T(y)).$ 
Consequently $y\in {\Cscr}_{(T\Ua p_T)}.$ 
Writing $p_{\sr T}=c_j$ 
we have $c_j^T=
(p_{\sr T})^T=p$ and this establishes the inclusion $\st.$
\par 
Conversely to any corner
$c$ of $T$ and $y\in {\Cscr}_{(T\Ua c)}$ we obtain that
$T(s^{\sr <}_{c^T,n-\sr 1}y)
=((T\Ua c)\Da c^T)=T$ and this establishes the reverse inclusion.
\QED
For example take 
$$T=\vcenter{
\halign{& \hfill#\hfill
\tabskip4pt\cr
\multispan{7}{\hrulefill}\cr
\ssa
\vb & 1 &  &  2 & &  5 & \ts\vb\cr
\vsa
&&&&\multispan{3}{\hrulefill}\cr
\ssa
\vb & 3 &  &  4 & \ts\vb\cr
\vsa
&&\multispan{3}{\hrulefill}\cr
\ssa
\vb & 6 & \ts\vb\cr
\vsa
\multispan{3}{\hrulefill}\cr}}$$
Then
$$\phi_{\sr 5}(T\Ua c_{\sr 1})=\vcenter{
\halign{& \hfill#\hfill
\tabskip4pt\cr 
\multispan{5}{\hrulefill}\cr
\ssa
\vb & 1 &  &  2 &  \ts\vb\cr
\vsa
&&&&\multispan{0}{\hrulefill}\cr
\ssa
\vb & 3 &  &  4 & \ts\vb\cr
\vsa
&&\multispan{3}{\hrulefill}\cr
\ssa
\vb & 5 & \ts\vb\cr
\vsa
\multispan{3}{\hrulefill}\cr}},\ 
\phi_{\sr 2}(T\Ua c_{\sr 2})=\vcenter{
\halign{&\hfill#\hfill
\tabskip4pt\cr
\multispan{7}{\hrulefill}\cr
\ssa
\vb & 1 &  &  3 & & 4 & \ts\vb\cr
\vsa
&&\multispan{5}{\hrulefill}\cr
\ssa
\vb & 2 &  \ts\vb\cr
\vsa
&&\multispan{0}{\hrulefill}\cr
\ssa
\vb & 5 & \ts\vb\cr
\vsa
\multispan{3}{\hrulefill}\cr}},\ 
\phi_{\sr 2}(T\Ua c_{\sr 3})=\vcenter{
\halign{&\hfill#\hfill
\tabskip4pt\cr
\multispan{7}{\hrulefill}\cr
\ssa
\vb & 1 &  &  3 & &  4 & \ts\vb\cr
\vsa
&&&&\multispan{3}{\hrulefill}\cr
\ssa
\vb & 2 &  &  5 & \ts\vb\cr
\vsa
\multispan{5}{\hrulefill}\cr}}.$$
Hence ${\Cscr}_T=s^{\sr <}_{\sr 5,5}{\Cscr}_{(T\Ua c_1)}\coprod 
s^{\sr <}_{\sr 2,5}{\Cscr}_{(T\Ua c_2)} 
\coprod s^{\sr <}_{\sr 2,5}{\Cscr}_{(T\Ua c_3)}.$
\subsubsection{}\label{3.1.2}
Let $\Cscr$ be a cell of $\bS_{n- {\sr 1}}$ and
put $T=T_\Cscr.$  By \ref{2.4.7} 
$s^{\sr <}_{ i,n-\sr 1}\Cscr$ lies in a unique cell ${\Cscr}_{(\phi_i(T)\Da i)}$
of $S_n.$  When it does not cause an ambiguity 
we shall denote this cell also by $s^{\sr <}_{i,n-\sr 1} \Cscr.$
\subsubsection{}\label{3.1.3}
Recall the notion of $\om_i(T),\ \om^i(T)$ from \ref{2.4.2}, 
$j^T$ from \ref{2.4.8} and $c^T,\ s_c(T)$
from \ref{2.4.9}.
\begin{lemma} For a given tableau $T$ and two numbers 
$j_{\sr 1},j_{\sr 2}\in {\Na}$ such that $j_{\sr 1},j_{\sr 2}\not 
\in<T>$ 
and $j_{\sr 1},j_{\sr 2}>\om_{\sr 1}(T)$ one has that
$j_{\sr 1}\leq j_{\sr 2}$ implies 
$$j_{\sr 1}^{\sr T}\leq j_{\sr 2}^{\sr T} \eqno{(*)}$$
In particular given $c,c\pr$ two corners of $T$ one has:
\begin{itemize}
\item[(i)]    $c<c\pr$ implies 
$c^{\sr T}\geq {c\pr}^{\sr T}.$
\item[(ii)]   The segments $s_{c}(T),s_{c\pr}(T)$ cannot cross, but 
              coalesce exactly when equality holds in (i).
\item[(iii)]  Let $c<c\pr$ and $c^{\sr T}={c\pr}^{\sr T}$ then 
              for all the corners 
              $c\prpr\ :\ c<c\prpr<c\pr$ one has ${c\prpr}^T=c^T.$
\end{itemize}
\end{lemma}
\Pf
For a given row $R=(a_{\sr 1},\cdots)$ and two numbers 
$j_{\sr 1},j_{\sr 2} \in{\Na}\ \setminus <R>$ such that
$j_{\sr 1},j_{\sr 2}>a_{\sr 1}$ one has that $(R\uar j_i)$ is defined 
for $i=1,2$ and
$j_{\sr 1}\leq j_{\sr 2}$
implies 
$$j_{\sr 1}^R\leq j_{\sr 2}^R.$$
Using this result inductively on rows we obtain $(*).$ 
\par
(i).\ \ \ Set $c=c(r,s)$ and $p_r={c\pr}^{\sr T^{r,\infty}}.$
Since $p_r\in <T^r>$ one has $p_r\leq \om^r(T).$
Hence by $(*)$ one gets ${c\pr}^{\sr T}={p_r}^{\sr T^{1,r- {\sr 1}}}\leq 
({\om^r(T)})^{\sr T^{1,r- {\sr 1}}}=c^{\sr T}.$ 
\par
(ii) and (iii) follow from (i).
\QED  
For example consider $T$ from \ref{3.1.1}. One has $c_{\sr 1}^{\sr T}=5,\ 
c_{\sr 2}^{\sr T}=c_{\sr 3}^{\sr T}=2$ and $s_{c_3}(T)$ coalesces with
$s_{c_2}(T).$
\subsubsection{}\label{3.1.4}
We can change the order of deletions of corners if their segments do not coalesce. 
This may be precisely expressed as follows.
\begin{lemma} Given a tableau $T$ and two distinct numbers 
$j_{\sr 1},j_{\sr 2}\in {\Na}$ such that $j_{\sr 1}, j_{\sr 2}\not\in <T>$ 
and $\om_{\sr 1}(T) < j_{\sr 1}<j_{\sr 2}.$
If $j_{\sr 1}^{\sr T}\ne j_{\sr 2}^{\sr T}$ then 
$$((T\uar j_{\sr 1})\uar j_{\sr 2})=((T\uar j_{\sr 2})\uar j_{\sr 1}).\eqno{(*)}$$
In particular

\begin{itemize}
\item[(i)] if $c,c\pr$ are two corners of $T$ such that 
$c^{\sr T}\ne {c\pr}^{\sr T}$
            then 
$$((T\Ua c)\Ua c\pr)=((T\Ua c\pr)\Ua c).$$
\end{itemize}
\end{lemma}
\Pf
Given a row $R=(a_{\sr 1},\cdots)$ and two numbers 
$j_{\sr 1},j_{\sr 2} \not\in <R>$ such that $a_{\sr 1}<j_{\sr 1}<j_{\sr 2}.$ 
Set $a_{i_k}=j_k^{\sr R}$ for $k=1,2.$ If $i_{\sr 1}<i_{\sr 2}$ then 
$((R\uar j_{\sr 1})\uar j_{\sr 2})=
((R\uar j_{\sr 2})\uar j_{\sr 1}).$ Using this result inductively on rows we 
obtain $(*).$
\par
To show (i) let $c<c\pr, \ c=c(r,s).$ Consider the procedure of pushing out
from $T^{r,\infty}.$ If $c^{\sr T}>{c\pr}^{\sr T}$ then 
in particular  
$j_r={c\pr}^{\sr T^{r,\infty}}<\om_r(T)$
and $(T\Ua c)^{r,\infty}=(T^{r,\infty}-\om^r(T)).$ 
Hence 
\goodbreak
$$((T\Ua c)\Ua c\pr)=
\left(\begin{array}{c}(T^{1,r- {\sr 1}}\uar \om^r(T))\uar j_r\cr
                            (T^{r,\infty}-\om^r(T))\Ua c\pr\cr\end{array}\right)
{\buildrel {(*)} \over =}\left(\begin{array}{c}(T^{1,r- {\sr 1}}\uar j_r) \uar\om^r(T)\cr
(T^{r,\infty}\Ua c\pr)-\om^r(T)\cr\end{array}\right)
=((T\Ua c\pr)\Ua c).$$
\QED 
As an example let us consider once more $T$ from \ref{3.1.1}. One has
$$\left(\left(\ 
\vcenter{
\halign{&\hfill#\hfill
\tabskip4pt\cr
\multispan{7}{\hrulefill}\cr
\ssa
\vb & 1 &  &  2 & &  5 & \ts\vb\cr
\vsa
&&&&\multispan{3}{\hrulefill}\cr
\ssa
\vb & 3 &  &  4 & \ts\vb\cr
\vsa
&&\multispan{3}{\hrulefill}\cr
\ssa
\vb & 6 & \ts\vb\cr
\vsa
\multispan{3}{\hrulefill}\cr}}
\Ua c(1,3)\right)\Ua c(2,2)\right)=\left(\ 
\vcenter{
\halign{&\hfill#\hfill
\tabskip4pt\cr
\multispan{5}{\hrulefill}\cr
\ssa
\vb & 1 &  &  2 &  \ts\vb\cr
\vsa
&&&&\multispan{0}{\hrulefill}\cr
\ssa
\vb & 3 &  &  4 & \ts\vb\cr
\vsa
&&\multispan{3}{\hrulefill}\cr
\ssa
\vb & 6 & \ts\vb\cr
\vsa
\multispan{3}{\hrulefill}\cr}}\Ua c(2,2)\ \right)=
\vcenter{
\halign{&\hfill#\hfill
\tabskip4pt\cr
\multispan{5}{\hrulefill}\cr
\ssa
\vb & 1 &  &  4&  \ts\vb\cr
\vsa
&&\multispan{3}{\hrulefill}\cr
\ssa
\vb & 3 &  \ts\vb\cr
\vsa
&&\multispan{0}{\hrulefill}\cr
\ssa
\vb & 6 & \ts\vb\cr
\vsa
\multispan{3}{\hrulefill}\cr}}.$$
On the other hand 
$$\left(\left(\ 
\vcenter{
\halign{&\hfill#\hfill
\tabskip4pt\cr
\multispan{7}{\hrulefill}\cr
\ssa
\vb & 1 &  &  2 & &  5 & \ts\vb\cr
\vsa
&&&&\multispan{3}{\hrulefill}\cr
\ssa
\vb & 3 &  &  4 & \ts\vb\cr
\vsa
&&\multispan{3}{\hrulefill}\cr
\ssa
\vb & 6 & \ts\vb\cr
\vsa
\multispan{3}{\hrulefill}\cr}}
\Ua c(2,2)\right)\Ua c(1,3)\right)
=\left(\ \vcenter{
\halign{&\hfill#\hfill
\tabskip4pt\cr
\multispan{7}{\hrulefill}\cr
\ssa
\vb & 1 &  &  4 & & 5 & \ts\vb\cr
\vsa
&&\multispan{5}{\hrulefill}\cr
\ssa
\vb & 3 &  \ts\vb\cr
\vsa
&&\multispan{0}{\hrulefill}\cr
\ssa
\vb & 6 & \ts\vb\cr
\vsa
\multispan{3}{\hrulefill}\cr}}\Ua c(1,3)\ \right)=\vcenter{
\halign{&\hfill#\hfill
\tabskip4pt\cr
\multispan{5}{\hrulefill}\cr
\ssa
\vb & 1 &  &  4&  \ts\vb\cr
\vsa
&&\multispan{3}{\hrulefill}\cr
\ssa
\vb & 3 &  \ts\vb\cr
\vsa
&&\multispan{0}{\hrulefill}\cr
\ssa
\vb & 6 & \ts\vb\cr
\vsa
\multispan{3}{\hrulefill}\cr}}.$$
\subsubsection{}\label{3.1.5}
Let $c=c(r,s)$ be some corner of $T$ and let $\dot T=(T\Ua c).$ 
\begin{lemma} Let $c\pr=c(p,t)$ be a corner of $\dot T.$ One has
$${c\pr}^{\sr \dot T}\begin{cases} >c^{\sr T},& {\rm if}\ p<r\\
                            <c^{\sr T},& {\rm otherwise}\\ \end{cases}.$$
\end{lemma}
\Pf
The proofs for $p<r$ and $p\geq r$ are the same so we prove the lemma only
for $p<r.$
\par
Let us show this by induction on $r.\ r>p$ implies that $r\geq 2.$ If $r=2$ then one has $c(r,s)=\om^{\sr 2}(T)$
so that $c^{\sr T}<\om^{\sr 2}(T).$ On the other hand by pushing up process
$\om^{\sr 2}(T)\in \dot T^1,$ thus 
$${c\pr}^{\sr \dot T}=\om^{\sr 1}(\dot T)\geq \om^{\sr 2}(T)>c^{\sr T}.$$
Assume that this is true for $r<q$ and show it for $r=q.$
If $p=1$ then exactly as in previous case
$${c\pr}^{\sr \dot T}=\om^{\sr 1}(\dot T)\geq c^{\sr T^{2,\infty}}>c^{\sr T}.$$
If $p>1$ then by induction assumption one has $c^{\sr T^{2,\infty}}<{c\pr}^{\sr \dot T^{2,\infty}}$
and by pushing up process $c^{\sr T^{2,\infty}}\in \dot T^1.$ By pushing up
process one has 
$${c\pr}^{\sr \dot T}=({c\pr}^{\sr \dot T^{2,\infty}})^{\dot T^1}\geq c^{\sr T^{2,\infty}}>c^{\sr T}.$$
\QED
\subsubsection{}\label{3.1.6}
Let $T=T^{1,k}_{1,l}$ then  $c_{\sr 1}=c(r,l)$ for some $r\geq 1.$
In particular one always has ${c_{\sr 1}}^{\sr T}=\om^{\sr 1}(T).$
\subsection
{\bf Standard Young tableaux and canonical elements}

\subsubsection{}\label{3.2.1}
We need to construct some ``canonical'' representatives of a given cell for
simplification of proofs and calculations. Recall the 
notion of $\ov w$ and of colligation from \ref{2.2.5}.
We need the following very simple but important construction. 
\begin{lemma} Consider $T=T_{1,l}^{1,k}.$ 

\begin{itemize}
\item[(i)]  Take $w\pr\in{\Cscr}_{T^{2,\infty}}$ and set $w=[w\pr,T^1].$ 
            Then $T(w)=T.$
\item[(ii)] Take $y\pr\in{\Cscr}_{T_{1,l-1}}$ and set 
            $y=[y\pr,\ov T_l].$ Then $T(y)=T.$
\end{itemize}
\end{lemma}
\Pf
(i)\ \ \ By  RS algorithm $T(w)=((\cdots(T^{2,\infty} \Da T^{\sr 1}_{\sr 1})\Da \cdots)\Da T^{\sr 1}_l).$
In that case roughly speaking the insertion of $T^{\sr 1}_i$ just pushes down the
$i$-th column $T_i.$ Formally 
 $T^j_{\sr 1}<T^{j+1}_{\sr 1}$ and $T^j_{i-\sr 1}<T^j_i<T^{j+1}_i$ 
for all $j\geq 1,\ i\geq 2$ gives $T(w)=T.$
\parno
(ii)\ \ \ Again using RS algorithm we get that
$T(y)=((\cdots(T_{1,l-1}\Da \om_l(T))\Da \cdots)\Da T^{\sr 1}_l).$
Note that $T^j_l>T^i_{l-1}$ for all $j:\ j\geq 1$ and all $i:\ 1\leq i\leq j,$
and also $T^j_l>T^{j-1}_l$ for all $j:\ j\geq 2.$ Hence $T(y)=T.$
\QED  
\subsubsection{}\label{3.2.2}
In what follows we frequently use the two ``canonical'' words defined below so 
we introduce  special 
notations for them. For  a given $T=T_{1,l}^{1,k}$ 
set
 $w_r(T):=[T^k,\cdots,T^1],\  w_c(T):=[\ov T_1,\cdots,\ov T_l].$ 
\par
For example for $T$ from \ref{3.1.1} one has 
$w_r(T)=[6,3,4,1,2,5]$ and $w_c(T)=[6,3,1,4,2,5].$
\par 
The inductive use of lemma \ref{3.2.1} provides the
\begin{cor} $T(w_r(T))=T(w_c(T))=T.$
\end{cor}
\Pf
Let us prove the statement for $w_r(T)$ by the induction on the number of
 rows in~$T.$
\par
If $T$ has only one row that is $T=T^1=(T_1^1,\cdots,T^1_l)$ 
then $w_r(T)=[T^1_1,\cdots,T^1_l]$
        and by RS procedure 
$T(w_r(T))=T$.
\par
Now assume that the statement is true for a 
     tableau with $k-1$ rows
     and show that it is true for a tableau with $k$ 
     rows. In that case 
     $w_r(T)=[T^k,\cdots,T^2,T^1].$ Consider the 
     word  $w_r( T^{2,\infty})=
     [T^k,\cdots,T^2].$ By the induction hypothesis 
     we obtain that
     $T(w_r(T^{2,\infty}))=T^{2,\infty}.$ Then by 
     lemma \ref{3.2.1} $T(w_r(T))=T.$     
\parno
Exactly the same way we obtain $T(w_c( T))=T.$
\QED
\subsubsection{}\label{3.2.3}
Lemma \ref{3.2.1} provides us also the following
\begin{cor} Given a Young tableau $T=T_{1,l}^{1,k}$     
\begin{itemize}
\item[(i)]   For $i\ :\ 2\leq i\leq k$ let $w\pr\in {\Cscr}_{T^{i,\infty}}$ and 
             set $w=[w\pr,w_r(T^{1,i- {\sr 1}})].$ Then $T(w)=T.$
\item[(ii)]  For $i\ :\ 1\leq i\leq l-1$ let $x\pr\in{\Cscr}_{T_{1,i}}$ and set 
             $x=[x\pr,w_c(T_{i+1,l})].$
             Then $T(x)=T.$
\item[(iii)] For $i\ :\ 2\leq i\leq l$ let $y\pr\in{\Cscr}_{T_{i,\infty}}$ and 
             set 
             $y=[w_c(T_{1,i-1}),y\pr].$ Then $T(y)=T.$
\item[(iv)]  For $i\ :\ 1\leq i\leq k-1$ let $z\pr\in{\Cscr}_{T^{1,i}}$ and set
             $z=[w_r(T^{i+1,\infty})\, ,\,z\pr].$ Then $T(z)=T.$
\item[(v)]   For $i\ :\ 2\leq i\leq l$ let $w\pr\in {\Cscr}_{T^{i,\infty}}$ and $z\pr\in{\Cscr}_{T^{1,i-1}}.$ 
             Set $w=[w\pr,z\pr].$ Then $T(w)=T.$ 
\item[(vi)]  For $i\ :\ 1\leq i\leq k-1$ let $x\pr\in{\Cscr}_{T_{1,i}}$ and $y\pr\in{\Cscr}_{T_{i+1,\infty}}.$ 
             Set $x=[x\pr,y\pr].$ Then $T(x)=T.$ 
\end{itemize}
\end{cor}
\Pf
(i), (ii) generalize corollary \ref{3.2.2} and are proved similarly.
Part (iii)  is obtained from (i) by  \ref{2.4.15}. Indeed
if $T(y\pr)=T_{i,\infty}$, then by \ref{2.4.15}
$T(\ov{y\pr})=(T_{i,\infty})^{\dagger}=
({T^{\dagger}})^{i,\infty}.$ Now by (i) 
one has $T([\ov{y\pr},T_{i- {\sr 1}},\cdots,T_1])=T^{\dagger}.$ 
Using again \ref{2.4.15}
we obtain 
$$T([\ov T_1,\cdots,\ov T_{i- {\sr 1}},y\pr])=
      T^{\dagger}({[ \ov{y\pr},T_{i- {\sr 1}},\cdots,T_1]})=T.$$
Part (iv) is similarly  obtained from part (ii).
\par
To show (v) we use the following 
easy property of RS procedure. If $s,t$ are two words such that 
$T(s)=T(t)$ and $q$ is a word 
such that $<q>\cap <s>=
\emptyset$ then $T([s,q])=T([t,q]).$
\parno  
By \ref{3.2.2} one has $T(w\pr)=T^{i,\infty}=T([w_r(T^{i,m})]).$
Thus $T(w)=T([w\pr,z\pr])=T([w_r(T^{i,m}),z\pr])=T$ just by (iv).
\par
Part (vi) is similarly obtained from part (iii).
\QED
For example consider
$$T=\vcenter{
\halign{&\hfill#\hfill
\tabskip4pt\cr
\multispan{11}{\hrulefill}\cr
\ssa
\vb &1 & & 2  && 3 & & 7 & & 11& \ts\vb\cr
\vsa
&&&&&&\multispan{5}{\hrulefill}\cr
\ssa
\vb & 4 && 5 && 8  &  \ts\vb\cr
\vsa
&&&&\multispan{3}{\hrulefill}\cr
\ssa
\vb & 6 && 10 &  \ts\vb\cr
\vsa
&&\multispan{3}{\hrulefill}\cr
\ssa
\vb & 9 &  \ts\vb\cr
\vsa
\multispan{3}{\hrulefill}\cr}}$$
One has
$$T^{3,\infty}=\vcenter{
\halign{&\hfill#\hfill
\tabskip4pt\cr
\multispan{5}{\hrulefill}\cr
\ssa
\vb & 6 && 10 &  \ts\vb\cr
\vsa
&&\multispan{3}{\hrulefill}\cr
\ssa
\vb & 9 &  \ts\vb\cr
\vsa
\multispan{3}{\hrulefill}\cr}}
\qquad {\rm and}\qquad
T^{1,2}=\vcenter{
\halign{&\hfill#\hfill
\tabskip4pt\cr
\multispan{11}{\hrulefill}\cr
\ssa
\vb &1 & & 2  && 3 & & 7 & & 11& \ts\vb\cr
\vsa
&&&&&&\multispan{5}{\hrulefill}\cr
\ssa
\vb & 4 && 5 && 8  &  \ts\vb\cr
\vsa
\multispan{7}{\hrulefill}\cr}}$$
so that $w\pr=[9,10,6]\in {\Cscr}_{T^{3,\infty}}$
and $z\pr=[4,1,5,2,3,8,7,11]\in{\Cscr}_{T^{1,2}}.$ Hence
$w=[w\pr,w_r(T^{1,2})]=[9,10,6,4,5,8,1,2,3,7,11]$ as well as
$z=[w_r(T^{3,\infty}),z\pr]=
[9,6,10,4,1,5,2,3,8,7,11]$
and $v=[w\pr,z\pr]=[ 9,10,6,4,1,5,2,3,8,7,11]$
belong to ${\Cscr}_T.$
\subsubsection{}\label{3.2.4}
Let us note that the converse to lemma \ref{3.2.1}(i) is also true, moreover
\begin{lemma} Given $T=T_{1,l}^{1,k}.$
Let $x$ (resp. $y$) be a word  such that $T=T([\ov T_1, x])$ (resp. $T=T([y,T^1]).)$
Then $T(x)=T_{2,\infty}$ (resp. $T( y)=T^{2,\infty}.$)
\end{lemma}
\Pf
By \ref{2.4.15} it is enough to show one of the 
assertions.
The proof is based on a counting principle. Let us show it for $x.$
\par
One has $|T_1|=k.$ Recall the notion of $Q(w)$ from \ref{2.4.18}.
Consider the following subsets of ${\Cscr}_T$
\begin{itemize}
\item[(i)]     $\Uscr=\{w\in {\Cscr}_T\ :\ w=[\ov T_1,w\pr]\}.$
\item[(ii)]    $\Vscr=\{[\ov T_1,w]\ :\ T(w)=T_{2,\infty}\}.$ 
\item[(iii)]   $\Wscr=\{w\in {\Cscr}_T\ :\  (Q(w))_1=(1,\cdots,m)\}.$
\end{itemize}
\parno
Applying  corollary \ref{3.2.3}(iii) to $i=2$ we obtain $\Vscr\st\Uscr.$ 
On the other hand by RS procedure 
$$(Q([\ov T_1,x]))_1 =\left(\begin{array}{c}1\cr\vdots\cr l\cr\end{array}\right)$$ 
hence $\Uscr \st \Wscr.$
\par 
Yet by \ref{2.4.18}
$$\begin{array}{rcl}
\vert \Vscr\vert&
           =&\vert\{Q\in \bT_{n-m}\ :\ \sh(Q)=\sh(T_{2,\infty})\}\vert,\cr
            \vert \Wscr\vert&
           =&\vert\{Q\in {\boT}_{\{l+1,\cdots,n\}}\ :\ \sh(Q)=\sh(T_{2,\infty})\}\vert,\cr
\end{array}$$
which provides $\Uscr=\Vscr.$
\QED
The same argument gives us that the assertions converse to corollary 
\ref{3.2.3} (i), (iii) are also true.
\par 
The converse of lemma \ref{3.2.1}(ii) is false in general as it is shown by the 
\par
\goodbreak
{\bf Example.}
\par
Consider $$T=
\vcenter{
\halign{&\hfill#\hfill
\tabskip4pt\cr
\multispan{5}{\hrulefill}\cr
\ssa
\vb&1 & & 2  & \ts\vb\cr
\vsa
&&\multispan{3}{\hrulefill}\cr
\ssa
\vb & 3 &  \ts\vb\cr
\vsa
\multispan{3}{\hrulefill}\cr}}.$$ Note that  
$w=[1,3,2]=[1,3,T_2]\in {\Cscr}_T$, 
and 
$$T([1,3])=\vcenter{
\halign{&\hfill#\hfill
\tabskip4pt\cr
\multispan{5}{\hrulefill}\cr
\ssa
\vb&1 & & 3  & \ts\vb\cr
\vsa
\multispan{5}{\hrulefill}\cr}}\ne 
\vcenter{
\halign{&\hfill#\hfill
\tabskip4pt\cr
\multispan{3}{\hrulefill}\cr
\ssa
\vb &1&\ts\vb\cr
\vsa
\multispan{0}{\hrulefill}\cr
\ssa
\vb & 3 &  \ts\vb\cr
\vsa
\multispan{3}{\hrulefill}\cr}}.$$
Respectively the assertions converse to corollary \ref{3.2.3} (ii), (iv), (v), (vi)
are false in general.
\subsubsection{}\label{3.2.5}
Now we continue the description of representatives for a given cell.
Given $T=T_{1,l}\in {\boT}_E$ and $h\in {\Na}\ :\ h\not\in E.$ 
We want to describe special words connected to an insertion of $h$ into $T.$
\begin{lemma} Take $T,h$ as above. Assume that $T^{\sr 1}_{i- {\sr 1}}<h< T^{\sr 1}_i$ for some
$i\ :\ 1\leq i\leq l.$ Then
\begin{itemize}
\item[(i)]   $(T\Da h)=T([w,T^{\sr 1}_i,(T^1\dar h)])$ for any word 
             $w\in {\Cscr}_{T^{2,\infty}}.$
\item[(ii)]  For all $j\geq i$ one has $(T\Da h)=T([y_j,h,z_{j+\sr 1}])$
             for any words \bk 
             $y_j\in{\Cscr}_{T_{1,j}},\ z_{j+\sr 1}\in{\Cscr}_{T_{j+1,\infty}}$
\item[(iii)] In particular $(T\Da h)=T([w_c( T_{1,j}),h,w_c(T_{j+1,\infty})])$ for all 
             $j\ :\ i\leq j\leq l.$
\end{itemize}
\end{lemma}
\Pf
Set $U:=(T\Da h).$
\begin{itemize} 
\item[(i)] To show (i)
we recall that by the insertion algorithm
$$U=(T \Da h)=\left(\begin{array}{c}T^1 \dar h\cr T^{2,\infty} \Da T_i^{\sr 1}\cr\end{array}\right).$$
Set $w\pr=[w,T_i^{\sr 1}].$ Since $T(w)=T^{2,\infty}$ by hypothesis, the
insertion algorithm gives $T(w\pr)=(T^{2,\infty}\Da T^{\sr 1}_i)=U^{2,\infty},$ 
hence lemma \ref{3.2.1} gives $T([w,T^{\sr 1}_i,(T^1\dar h)])=U.$
\item[(ii)] Since $h<T^{\sr 1}_i<T^{\sr 2}_i<\cdots$ the insertion of $h$ into $T$ does
not affect $T_{i+1,\infty}$ so  $U_{i+1,\infty}=T_{i+1,\infty}.$
That means that
$$U_{1,j}=(T_{1,j} \Da h) \qquad {\rm for\ every}\ j\geq i$$ 
Hence  for every $j\ :\ i\le j\le l$ the word
$y\pr_j=[y_j,h]$ satisfies
$T(y\pr_j)=(T_{1,j} \Da h)=U_{1,j}.$ Then corollary \ref{3.2.3}(vi) gives  
$T([y_j,h,z_{j+\sr 1}])=U.$
\end{itemize}
\QED
\subsection
{\bf Inductive definition of the set of offsprings.  }
\subsubsection{}\label{3.3.1}
Consider $w\in \bS_n.$ Write $w=[a_{\sr 1},\cdots,a_n]$ and recall that 
by \ref{2.2.6} (ii)  $w=s^{\sr <}_{a_n,n-\sr 1}y,$
 where $y=\phi^{- {\sr 1}}_{a_n}([a_{\sr 1},\cdots,a_{n-\sr 1}]),\ y\in \bS_{n- {\sr 1}}.$
By  \ref{3.1.2} given a cell $\Cscr\st \bS_{n- {\sr 1}}$
we can define a cell $s^{\sr <}_{i,n-\sr 1}\Cscr$ in  $\bS_n.$ 
\par
By \ref{2.2.2} Remark and \ref{2.2.4} Remark 2   for $i\ :\ 1\le i\le n-2$ one has
$ w \dos ws_i$ iff $y \dos ys_i.$ Thus 
\begin{lemma} If ${\Cscr}_1,{\Cscr}_2$ are cells of $\bS_{n- {\sr 1}}$ satisfying 
${\Cscr}_1\dor {\Cscr}_2$ then
$s^{\sr <}_{i,n-\sr 1}{\Cscr}_1\dor 
s^{\sr <}_{i,n-\sr 1}{\Cscr}_2$ for all $i\ :\ 1\leq i\leq n-1.$ 
\end{lemma}
\Pf
By definition of induced Duflo order it is enough to show this only
for the case when ${\Cscr}_2$ is an offspring of ${\Cscr}_1$ which is straightforward
corollary of \ref{2.2.8}.
\QED
It can occur that ${\Cscr}_1\dos {\Cscr}_2$ in $\bS_{n- {\sr 1}}$
and yet $s^{\sr <}_{i,n-\sr 1}{\Cscr}_1= s^{\sr <}_{i,n-\sr 1}{\Cscr}_2.$
This is shown by the
\par
{\bf Example.}
\par
Consider $\bS_2=\{{\Cscr}_1=\{[1,2]\},{\Cscr}_2=\{[2,1]\}\},$
it is obvious that ${\Cscr}_1\dos {\Cscr}_2.$ Now consider sets  $s_{\sr 2}{\Cscr}_1=
\{[1,3,2]\},\ s_{\sr 2} {\Cscr}_2=\{[3,1,2]\}.$ One has 
$$T([1,3,2])=T([3,1,2])=
\vcenter{
\halign{&\hfill#\hfill
\tabskip4pt\cr
\multispan{5}{\hrulefill}\cr
\ssa
\vb&1 & & 2  & \ts\vb\cr
\vsa
&&\multispan{3}{\hrulefill}\cr
\ssa
\vb & 3 &  \ts\vb\cr
\vsa
\multispan{3}{\hrulefill}\cr}}.$$
Thus $s_{\sr 2}{\Cscr}_1=s_{\sr 2} {\Cscr}_2.$
\subsubsection{}\label{3.3.2}
Let ${\Cscr}_T$ be a cell in $\bS_n.$ By proposition \ref{3.1.1}
it can be considered as a disjoint union 
$${\Cscr}_T=\coprod_{i=1}^m  s^{\sr <}_{p_i,n-\sr 1}{\Cscr}_{(T\Ua c_i)}$$
where $p_i={c_i}^{\sr T}.$
\par
Consider right multiplication by $s_j \in \Pi_{n- {\sr 2}}$ on some element $w\in {\Cscr}_T.$ 
Since $w=s^{\sr <}_{p_i,n-\sr 1}y$ for $y\in {\Cscr}_{(T\Ua c_i)}$ one has
${\Cscr}_{ws_j}=s^{\sr <}_{p_i,n-\sr 1}{\Cscr}_{ys_j}.$ 
Recall the notion of $\Dscr({\Cscr}_T)$ from \ref{2.5.2}. Set for $T\in \bT_{n- {\sr 1}}$  
$$s^{\sr <}_{i,n-\sr 1}\Dscr({\Cscr}_T):=\{\ s^{\sr <}_{i,n-\sr 1}
\Cscr \ \mid \ \Cscr \in\Dscr({\Cscr}_T) \, \}=\{{\Cscr}_{(\phi_i(T\pr)\Da i)}\ \mid\ T\pr\in\Dscr(T)\, \}.$$
Set
$${\Dscr}_o({\Cscr}_T):=\bigcup_{i=1}^m s^{\sr <}_{p_i,n-\sr 1}\Dscr( {\Cscr}_{(T\Ua c_i)})\quad {\rm and}\quad
 {\Dscr}_o(T):=\bigcup_{i=1}^m\{(S\Da p_i)\ \mid\ S\in \Dscr(T\Ua c_i)\}.$$
these are subsets of respectively $\Dscr({\Cscr}_T)$ and  $\Dscr(T)$. 
\par
Set 
$${\Dscr}_n({\Cscr}_T):=
\Dscr({\Cscr}_T)\setminus {\Dscr}_o({\Cscr}_T)\quad {\rm and}\quad{\Dscr}_n(T):=\Dscr(T)\setminus {\Dscr}_o(T).$$
Clearly if $\Cscr\pr\in {\Dscr}_n({\Cscr}_T)$ 
then $\Cscr\pr={\Cscr}_{ws_{n- {\sr 1}}}$ where $w\in{\Cscr}_T$ and 
$ws_{n- {\sr 1}}\dgs w.$ Summarizing
\begin{cor} $\Dscr({\Cscr}_T)={\Dscr}_o({\Cscr}_T)\coprod {\Dscr}_n({\Cscr}_T).$
\end{cor}
\subsubsection{}\label{3.3.3}
Recall the notion of $\om^i(T)$ from \ref{2.4.2} and of $c^T$ from 
\ref{2.4.9}. Take some $T\in\bT_n.$ Let $\{c_i\}_{i=1}^m$ be the set of its
corners. 
\parno
If $|T^1|>|T^2|$ then $c_1$ is in the first row. In that case 
consider an array
$$S:=\left(\begin{array}{c}T^1-\om^{\sr 1}(T)\cr
                  T^{2,\infty}\Da\om^{\sr 1}(T)\cr \end{array}\right)\eqno{(*)}$$
and define the following conditions
\begin{itemize}
\item[(i)]  $|T^1|\geq |T^2|+2;$
\item[(ii)] $|T^1|= |T^2|+1$ and $\om^{\sr 1}(T)<\om^{\sr 2}(T);$
\end{itemize}
\parno
Set
$$S_T(c_1):=\begin{cases}S,& {\rm if}\ T\ {\rm satisfies\ (i)\ or\ (ii)};\\
                   \emptyset, & {\rm otherwise}.\\ \end{cases}$$
\par
If there exist $c_i$ such that $c_i^T=\om_{\sr 1}(T)$ then consider an array 
$$S:=\left(\begin{array}{l}(T\Ua c_i)^1\cr
           (T\Ua c_i)^{2,\infty}\Da \om^{\sr 1}(T)\cr\end{array}\right)\eqno{(**)}.$$
Set 
$$S_T(c_i):=\begin{cases}S,& {\rm if}\ c_i^T=\om_{\sr 1}(T)\ {\rm and}\ S\ {\rm is\ standard\ with}\ 
\sh(S)>\sh(T)\\
                  \emptyset & {\rm otherwise}.\end{cases}$$ 
Set $\Dscr\pr_n(T):=\{S_T(c_i)\}_{i=1}^m.$
\begin{theorem} $\Dscr(T)={\Dscr}_o(T)\bigcup \Dscr\pr_n(T).$ 
\end{theorem}
This theorem provides a justification to the rough idea that induced Duflo
order is given by lowering numbered boxes in the sense of RS procedure.
\par 
Indeed if $c_{\sr 1}$ is in the first row
then $S$ defined in $(*)$ is a tableau if and only if $T$ satisfies (i) or (ii).
In that case $S_T(c_1)$ is obtained by  lowering the (numbered) box from the first
row of $T$ to the second row by RS procedure. 
\par
If $c_i=c(j,|T^j|)$ where $j>1$ then the procedure described
by $(**)$ is a more complicated shuffle in which a corner of $T$ 
is deleted pushing out the last (numbered) box of the first row of $T$ which is then 
inserted into $(T\Ua c_i)^{2,\infty}.$ 
\par
We further note that any 
$T\pr\in {\Dscr}_o(T)$ is obtained by first deleting some corner $c$ of $T$ then computing
$T\prpr \in \Dscr(T\Ua c)$ and finally computing
$T\pr=(T\prpr \Da c^T).$ 
Thus \ref{3.3.2} and \ref{3.3.3} together describe the offsprings of ${\Cscr}_T$ by a precisely determined 
process of deleting corners and lowering (possibly with shuffling) numbered boxes.
\par
For example let us compute $\Dscr(T)$ for 
$$T=\vcenter{
\halign{&\hfill#\hfill
\tabskip4pt\cr
\multispan{9}{\hrulefill}\cr
\ssa
\vb & 1 &  &  2 & &  6 && 7 & \ts\vb\cr
\vsa
&&&&\multispan{5}{\hrulefill}\cr
\ssa
\vb & 3 &  &  5 & \ts\vb\cr
\vsa
&&\multispan{3}{\hrulefill}\cr
\ssa
\vb & 4 & \ts\vb\cr
\vsa
\multispan{3}{\hrulefill}\cr}}$$
Note that $S_T(c)\ne\emptyset$ only for $c=c_{\sr 1}=c(1,4).$ 
Hence $\Dscr\pr_n(T)=\{S_T(c_1)\}$ where 
$$S_T(c_1)=\vcenter{
\halign{&\hfill#\hfill
\tabskip4pt\cr
\multispan{7}{\hrulefill}\cr
\ssa
\vb & 1 &  &  2 & &  6 &\ts\vb\cr
\vsa
&&&&&&\multispan{0}{\hrulefill}\cr
\ssa
\vb & 3 &  &  5 && 7& \ts\vb\cr
\vsa
&&\multispan{5}{\hrulefill}\cr
\ssa
\vb & 4 & \ts\vb\cr
\vsa
\multispan{3}{\hrulefill}\cr}}$$
Now one has
$$(T\Ua c_{\sr 1})=\vcenter{
\halign{&\hfill#\hfill
\tabskip4pt\cr
\multispan{7}{\hrulefill}\cr
\ssa
\vb & 1 &  &  2 & &  6 &\ts\vb\cr
\vsa
&&&&\multispan{3}{\hrulefill}\cr
\ssa
\vb & 3 &  &  5 & \ts\vb\cr
\vsa
&&\multispan{3}{\hrulefill}\cr
\ssa
\vb & 4 & \ts\vb\cr
\vsa
\multispan{3}{\hrulefill}\cr}},\ (T\Ua c_{\sr 2})=
\vcenter{
\halign{&\hfill#\hfill
\tabskip4pt\cr
\multispan{9}{\hrulefill}\cr
\ssa
\vb & 1 &  &  5 & &  6 && 7 &\ts\vb\cr
\vsa
&&\multispan{7}{\hrulefill}\cr
\ssa
\vb & 3 &  \ts\vb\cr
\vsa
&&\multispan{0}{\hrulefill}\cr
\ssa
\vb & 4 & \ts\vb\cr
\vsa
\multispan{3}{\hrulefill}\cr}},\ 
(T\Ua c_{\sr 3})=\vcenter{
\halign{&\hfill#\hfill
\tabskip4pt\cr
\multispan{9}{\hrulefill}\cr
\ssa
\vb & 1 &  &  3 & &  6 && 7 &\ts\vb\cr
\vsa
 &&&&\multispan{5}{\hrulefill}\cr
\ssa
\vb & 4 &  &  5 & \ts\vb\cr
\vsa
\multispan{5}{\hrulefill}\cr}}$$
Applying further our algorithm to $(T\Ua c_{\sr 1}),\ (T\Ua c_{\sr 2}),\ 
(T\Ua c_{\sr 3})$ we get that ${\Dscr}_o({\Cscr}_T)$ is described by
the following
set of tableaux:
$$\vcenter{
\halign{&\hfill#\hfill
\tabskip4pt\cr
\multispan{7}{\hrulefill}\cr
\ssa
\vb & 1 &  &  2 & &  6 &\ts\vb\cr
\vsa
&&&&&&\multispan{0}{\hrulefill}\cr
\ssa
\vb & 3 &  &  5 && 7& \ts\vb\cr
\vsa
&&\multispan{5}{\hrulefill}\cr
\ssa
\vb & 4 & \ts\vb\cr
\vsa
\multispan{3}{\hrulefill}\cr}},\qquad
\vcenter{
\halign{&\hfill#\hfill
\tabskip4pt\cr
\multispan{7}{\hrulefill}\cr
\ssa
\vb & 1 &  &  2 & &  6 &\ts\vb\cr
\vsa
&&&&\multispan{3}{\hrulefill}\cr
\ssa
\vb & 3 &  &  5 & \ts\vb\cr
\vsa
&&&&\multispan{0}{\hrulefill}\cr
\ssa
\vb & 4 && 7& \ts\vb\cr
\vsa
\multispan{5}{\hrulefill}\cr}},\qquad
\vcenter{
\halign{&\hfill#\hfill
\tabskip4pt\cr
\multispan{7}{\hrulefill}\cr
\ssa
\vb & 1 &  &  2 & &  7 &\ts\vb\cr
\vsa
&&&&\multispan{3}{\hrulefill}\cr
\ssa
\vb & 3 &  &  5 & \ts\vb\cr
\vsa
&&&&\multispan{0}{\hrulefill}\cr
\ssa
\vb & 4 && 6& \ts\vb\cr
\vsa
\multispan{5}{\hrulefill}\cr}},$$
\vskip 0.5truecm
$$\vcenter{
\halign{&\hfill#\hfill
\tabskip4pt\cr
\multispan{9}{\hrulefill}\cr
\ssa
\vb & 1 &  &  2 & &  6 && 7 &\ts\vb\cr
\vsa
&&\multispan{7}{\hrulefill}\cr
\ssa
\vb & 3 & \ts\vb\cr
\vsa
&&\multispan{0}{\hrulefill}\cr
\ssa
\vb & 4 &\ts\vb\cr
\vsa
&&\multispan{0}{\hrulefill}\cr
\ssa
\vb & 5 &\ts\vb\cr
\vsa
\multispan{3}{\hrulefill}\cr}},\ 
\vcenter{
\halign{&\hfill#\hfill
\tabskip4pt\cr
\multispan{9}{\hrulefill}\cr
\ssa
\vb & 1 &  &  5 & &  6 && 7 &\ts\vb\cr
\vsa
&&\multispan{7}{\hrulefill}\cr
\ssa
\vb & 2 & \ts\vb\cr
\vsa
&&\multispan{0}{\hrulefill}\cr
\ssa
\vb & 3 &\ts\vb\cr
\vsa
&&\multispan{0}{\hrulefill}\cr
\ssa
\vb & 4 &\ts\vb\cr
\vsa
\multispan{3}{\hrulefill}\cr}},\ 
\vcenter{
\halign{&\hfill#\hfill
\tabskip4pt\cr
\multispan{7}{\hrulefill}\cr
\ssa
\vb & 1 &  &  2 & &  6 &\ts\vb\cr
\vsa
&&&&\multispan{3}{\hrulefill}\cr
\ssa
\vb & 3 && 5 & \ts\vb\cr
\vsa
&&\multispan{3}{\hrulefill}\cr
\ssa
\vb & 4 &\ts\vb\cr
\vsa
&&\multispan{0}{\hrulefill}\cr
\ssa
\vb & 7 &\ts\vb\cr
\vsa
\multispan{3}{\hrulefill}\cr}},\ 
\vcenter{
\halign{&\hfill#\hfill
\tabskip4pt\cr
\multispan{7}{\hrulefill}\cr
\ssa
\vb & 1 &  &  2 & &  7 &\ts\vb\cr
\vsa
&&&&\multispan{3}{\hrulefill}\cr
\ssa
\vb & 3 && 5 & \ts\vb\cr
\vsa
&&\multispan{3}{\hrulefill}\cr
\ssa
\vb & 4 &\ts\vb\cr
\vsa
&&\multispan{0}{\hrulefill}\cr
\ssa
\vb & 6 &\ts\vb\cr
\vsa
\multispan{3}{\hrulefill}\cr}},\ 
\vcenter{
\halign{&\hfill#\hfill
\tabskip4pt\cr
\multispan{7}{\hrulefill}\cr
\ssa
\vb & 1 &  &  2 & &  7 &\ts\vb\cr
\vsa
&&&&\multispan{3}{\hrulefill}\cr
\ssa
\vb & 3 && 6 & \ts\vb\cr
\vsa
&&\multispan{3}{\hrulefill}\cr
\ssa
\vb & 4 &\ts\vb\cr
\vsa
&&\multispan{0}{\hrulefill}\cr
\ssa
\vb & 5 &\ts\vb\cr
\vsa
\multispan{3}{\hrulefill}\cr}}.$$
Note that $S_T(c_{\sr 1})\in{\Dscr}_o(T),$ hence ${\Dscr}_n(T)=\emptyset.$
\par
To prove the theorem we need to consider
the effect 
of right multiplication by $s_{n-\sr 1}$  
on the elements of a given cell ${\Cscr}_T.$
 We do this in \ref{3.3.4} and then give the proof 
of \ref{3.3.3} in \ref{3.3.5}.
\subsubsection{}\label{3.3.4}
Take $w\in {\Cscr}_T$ and write 
$w=[a_{\sr 1},\cdots,a_n].$ If $a_n<a_{n-\sr 1}$ 
then by  \ref{2.2.4} Remark 2, $ws_{n-\sr 1}\dor w.$ We therefore  need to 
consider only  those $w$ with $a_n>a_{n-\sr 1}.$ 
Assume that $w$
is such a word. 
One has $ws_{n-\sr 1}=
[a_{\sr 1},a_{\sr 2},\cdots,a_n,a_{n-\sr 1}].$ 
\par
Let $T=T(w)$ be a tableau with $l$ columns (that is $|T^1|=l$)
and with $m$ corners. It determines the set $\{ (T\Ua c_i)\}_{i=1}^m.$
Each one of the $(T\Ua c_i)$ determines in turn the set of 
$\{((T\Ua c_i)\Ua c\pr_j)\}_{j=1}^{m_i}.$ It is obvious that 
 one has ${}_{n-2}T(w)=((T\Ua c_i)\Ua c\pr_j)$ 
for some $c_i,\,c\pr_j,$ 
and $a_{n-\sr 1}={c\pr_j}^{(T\Ua c_i)},\ a_n={c_i}^T.$ Hence the 
offsprings 
obtained by right multiplication by $s_{n-\sr 1}$ can be computed from considering
two-fold deletions of corners. We show that  only rather few two-fold 
deletions intervene. 
We  describe the results in terms of Young tableaux.
Set $\ddot T={}_{n-2}T(w)$ which we saw was some $((T\Ua c_i)\Ua c\pr_j).$
Set 
$\dot T={}_{n- {\sr 1}}T(w).$ Then $\dot T=(\ddot T\Da a_{n-\sr 1})
=(\ddot T \Da {c\pr_j}^{(T\Ua c_i)} )=(T\Ua c_i).$ 
The resulting offsprings are described in the following
\begin{lemma} Let $T$ be a tableau described above, such that $|T^1|=l.$  
Consider $w \in {\Cscr}_T,\ w=[a_{\sr 1},\cdots,a_n]$ with 
$a_n>a_{n- {\sr 1}}.$ Set $\ddot T={}_{n-2}T(w)$ and 
$\dot T={}_{n- {\sr 1}}T(w).$ Set $l\pr=|\ddot T^1|.$
Then $T\dor T( ws_{n-\sr 1})$ and one of   the following holds: 
\begin{itemize}
\item[(i)] If there exist a pair $i,j$ with 
$1<i<j\leq l\pr$ or $1< i=j\leq l\pr-2$ such that $\ddot T^{\sr 1}_{i-\sr 1}<
           a_{n-\sr 1}<\ddot T^{\sr 1}_i$ and
           $\ddot T^{\sr 1}_{j-\sr 1}<a_n<\ddot T^{\sr 1}_j,$ then there exist $y \in {\Cscr}_T,\quad 
           \al_e \in \Pi_{n-2}$ such that $T( ws_{n-\sr 1})=T(ys_e).$
\item[(ii)]  If $a_{n-\sr 1}<\ddot T^{\sr 1}_{l\pr}<a_n$ then  $T( ws_{n-\sr 1})=T.$ 
\item[(iii)] \begin{itemize}\item[(a)]
If $a_{n-\sr 1}>\ddot T_{l\pr}^{\sr 1}$ then $l=l\pr+2;$ 
\item[(b)] If $\ddot T^{\sr 1}_{l\pr-1}<a_{n-\sr 1}$ and $a_n<\ddot T^{\sr 1}_{l\pr}$ then $l=l\pr+1;$
\item[(c)] If $\ddot T^{\sr 1}_{l\pr-2}<a_{n-\sr 1}$ and $a_n<\ddot T^{\sr 1}_{l\pr-1},$ then $l=l\pr.$
\end{itemize}   
\end{itemize}
\parno
\quad\quad In all the three cases {\rm (a), (b)} and {\rm (c)} 
$T^{\sr 1}_{l-1}=a_{n-\sr 1},\  \om^{\sr 1}(T)=T^{\sr 1}_l=a_n$ 
         and 
$$T( ws_{n-\sr 1})=\left(\begin{array}{c}\dot T^1 \cr (\dot T^{2,\infty} 
\Da \om^{\sr 1}(T))\cr\end{array}\right).$$
\end{lemma}
\Pf
Set $S=T(ws_{n-\sr 1}).$ Let us choose ``good'' elements in ${\Cscr}_T$ and ${\Cscr}_S.$
\par
Assume that there exist $j<l\pr$ such that $\ddot T^{\sr 1}_{j-\sr 1}<a_n<\ddot T^{\sr 1}_j.$
Then since $a_{n-\sr 1}<a_n$ there must exist $i\le j$ such that $\ddot T^{\sr 1}_{i-\sr 1}<a_{n-\sr 1}
<\ddot T^{\sr 1}_i.$ By lemma \ref{3.2.5}(iii) we have $\dot T=T(x)$ where 
$x=([w_c(\ddot T_{1,j}), 
a_{n-\sr 1}, w_c(\ddot T_{j+1,\infty})]).$
Then $T=(T(x) \Da a_n).$ Note that $\dot T_{1,j}=(\ddot T_{1,j}\dar a_{n- {\sr 1}})=
T([w_c(\ddot T_{1,j}),a_{n-\sr 1}]).$ Then  by lemma \ref{3.2.5}(ii) we can write
$T =(T(x) \Da a_n) =T(y)$ where
$$y=\begin{cases}
[w_c(\ddot T_{1,j}),a_{n-\sr 1},a_n, w_c(\ddot T_{j+1,\infty})],& {\rm if}\ i<j, \cr
       [w_c(\ddot T_{1,j}),a_{n-\sr 1},\ov {\ddot T}_{j+1},a_n,w_c(\ddot T_{j+2,\infty})],
& {\rm if}\ i=j. \cr \end{cases}$$
Note that in case $i=j$ we obtain that 
$$\dot T^1=(\ddot T^1 \Da a_{n-\sr 1})=
(\ddot T_{\sr 1}^{\sr 1},\cdots,\ddot T^{\sr 1}_{j-\sr 1},
a_{n-\sr 1},\ddot T^{\sr 1}_{j+\sr 1},\cdots),$$ 
so that $a_{n- {\sr 1}}<a_n<\ddot T_{j+\sr 1}^{\sr 1}.$ That is why we can insert $a_n$ only after 
inserting the column $\ddot T_{j+1}.$ 
\par
Using once more lemma \ref{3.2.5}(ii) for the case $i=j$ in definition of $y$ we 
can also write $T=T(z)$ where 
$$z=\begin{cases}[w_c(\ddot T_{1,j}),a_{n-\sr 1},a_n, w_c(\ddot T_{j+1,\infty})],&{\rm if}\ i<j, \cr
           [w_c(\ddot T_{1,j+1}),a_{n-\sr 1},a_n, w_c(\ddot T_{j+2,\infty})],&{\rm if}\ i=j. \cr 
\end{cases} 
                                                                                       \eqno(*)$$
In exactly the same way we obtain 
$S=((\ddot T \Da a_n)\Da a_{n-\sr 1})=T(u)$ where
$$   u=\begin{cases}[w_c(\ddot T_{1,j}),a_n,a_{n-\sr 1}, w_c(\ddot T_{j+1,\infty})],& {\rm if}\ i<j,\cr
          [w_c(\ddot T_{1,j+1}),a_n,a_{n-\sr 1}, w_c(\ddot T_{j+2,\infty})],& {\rm if}\ i=j. \cr
\end{cases}
                                                                                       \eqno(**)$$
\par
We now describe how cases (i)-(iii) result.
\begin{itemize}
\item[(i)]\begin{itemize}\item[(a)]  If $i<j<l\pr$ or $i=j<l\pr-2$ then the hypotheses 
                of (i) are fulfilled.  
                In this case the above computation shows that the conclusion of (i)
                is obtained. Indeed  denote by $e$ the place of $a_{n- {\sr 1}}$ 
                in $z$ and note that $e<n-1.$
                Then $(*)$ and $(**)$ gives $S=T(zs_e).$
\item[(b)] To complete case (i) it remains to study the case $j=l\pr$ and $i<l\pr.$ 
                By the insertion algorithm $T^1=(\ddot T^{\sr 1}_{\sr 1},\cdots\ddot T^{\sr 1}_{i-\sr 1},
                a_{n-\sr 1},\cdots,\ddot T^{\sr 1}_{l\pr-1},a_n)$ and
                $$T(w) = \left(\begin{array}{c}T^1\cr (\ddot T^{2,\infty}\Da \ddot T^{\sr 1}_i)\Da \ddot T^{\sr 1}_{l\pr}\cr\end{array}\right).$$
                On the other hand by the insertion algorithm  $S^1=T^1$ and
              $$S =\left(\begin{array}{c}T^1\cr (\ddot T^{2,\infty}\Da \ddot T^{\sr 1}_{l\pr})\Da \ddot T^{\sr 1}_i\cr\end{array}\right).$$
              Set $y=[w_r(\ddot T^{2,\infty}),\ddot T^{\sr 1}_i,\ddot T^{\sr 1}_{l\pr},T^1]$ and denote by $e$ the place of 
              $\ddot T^{\sr 1}_i$ in $y.$ Then $S=T(ys_e)$ which is conclusion of (i). 
\end{itemize}
\item[(ii)]  If $a_{n-\sr 1}<\ddot T^{\sr 1}_{l\pr}<a_n$ then the hypotheses of (ii) are
              fulfilled. In this case
$$\begin{array}{rcl}T =((\ddot T\Da a_{n- {\sr 1}})\Da a_n)&=&\left(\begin{array}{c}\dot T^1+ a_n\cr 
                                               \dot T^{2,\infty}\cr\end{array}\right)\cr
           S =((\ddot T\Da a_n)\Da a_{n- {\sr 1}})&=&\left(
        \left(\begin{array}{c}\ddot T^1+ a_n\cr \ddot T^{2,\infty}\cr\end{array}\right)\Da a_{n- {\sr 1}}\right )\cr
                                           &=&
        \left(\begin{array}{c}\dot T^1+ a_n\cr \dot T^{2,\infty}\cr\end{array}\right)=T.\cr
\end{array}$$
             This is just the conclusion of  (ii).
 
\item[(iii)]\begin{itemize}\item[(a)]  If $i=j=l\pr+1$ then $\ddot T^{\sr 1}_{l\pr}<a_{n-\sr 1}$  which are 
          hypotheses of iii (a).
          In that case $T^1=((\ddot T^1\dar a_{n-\sr 1})\dar a_n)=(\ddot T^1 +a_{n-\sr 1}+ a_n)$ so that
          $$T=\left(\begin{array}{c}T^1\cr \ddot T^{2,\infty}\cr\end{array}\right).$$
          Hence $l=l\pr+2,\ T^{\sr 1}_{l-1}=a_{n-\sr 1},\ \ \dot T^1=(\ddot T^1+a_{n-\sr 1})=(T^1-a_n)$ and
          $\om^{\sr 1}(T)=a_n.$  Thus 
$$S=((\ddot T\Da \om^{\sr 1}(T))\Da T^{\sr 1}_{l-1})=\left(\begin{array}{c}\dot T^1\cr \dot T^{2,\infty}\Da \om^{\sr 1}(T)\cr\end{array}\right)$$
exactly as in conclusion of (iii).
\item[(b)]  If $j=l\pr$ and $i=l\pr$ then $\ddot T^{\sr 1}_{l\pr-1}<a_{n-\sr 1},a_n<\ddot T^{\sr 1}_{l\pr}$ 
              which are the hypotheses
           of iii (b). In that case we get
$$T=((\ddot T\Da a_{n-\sr 1})\Da a_n)=\left(\begin{array}{c}\dot T^1+a_n\cr \dot T^{2,\infty}\cr
\end{array}\right)$$
             where $\dot T^1=(\ddot T^1\dar a_{n-\sr 1}),\ \dot T^{2,\infty}=
             (\ddot T^{2,\infty}\dar \ddot T^{\sr 1}_{l\pr}).$ Hence
             $T^1=(\ddot T^{\sr 1}_{\sr 1},\cdots,\ddot T^{\sr 1}_{l\pr-1},a_{n-\sr 1},a_n)$ and 
             $T^{2,\infty}=\dot T^{2,\infty}.$ In particular $l=l\pr+1$ and  
             $T^{\sr 1}_{l-1}=a_{n-\sr 1},\  \om^{\sr 1}(T)=a_n.$ Then
$${}_{n- {\sr 1}}T(ws_{n-\sr 1})=(\ddot T\Da \om^{\sr 1}(T))=
\left(\begin{array}{c}\ddot T^1\dar \om^{\sr 1}(T)\cr 
                                                             \dot T^{2,\infty}\cr\end{array}\right)$$
which provides us the conclusion of (iii), namely,
$$S=({}_{n- {\sr 1}}T(ws_{n-\sr 1})\Da T^{\sr 1}_{l-\sr 1})= \left(\begin{array}{c}
\ddot T^1\dar T^{\sr 1}_{l-1}\cr
\dot T^{2,\infty}\Da \om^{\sr 1}(T)\cr\end{array}\right)=
                       \left(\begin{array}{c}\dot T^1\cr \dot T^{2,\infty}\Da \om^{\sr 1}(T)\cr
\end{array}\right).$$             
\item[(c)] It remains to consider $i=j=l\pr-1.$ Here $\ddot T^{\sr 1}_{l\pr-2}<a_{n-\sr 1},a_n< \ddot T^{\sr 1}_{l\pr-1}$
            which are the hypotheses of iii(c). In that case 
$$\dot T=(\ddot T\Da a_{n-\sr 1})=\left(\begin{array}{c}\ddot T^1 \dar a_{n-\sr 1}\cr
                                \ddot T^{2,\infty}\Da \ddot T^{\sr 1}_{l\pr-\sr 1}\cr\end{array}\right)$$
      and so $\dot T^1=(\ddot T^1 \dar a_{n-\sr 1})=
     (\ddot T^{\sr 1}_{\sr 1},\cdots,\ddot T^{\sr 1}_{l\pr-2},a_{n-\sr 1},\ddot T^{\sr 1}_{l\pr})$
     Note that $a_{n-\sr 1}<a_n<\ddot T^{\sr 1}_{l\pr}$ which gives $a_n=\om^{\sr 1}(T)$ and $l=l\pr.$
     Now consider $S.$ By RS algorithm
   $$ {}_{n-\sr 1}T(ws_{n-\sr 1})=
                  \left(\begin{array}{c}\ddot T^1 \dar \om^{\sr 1}(T)\cr 
\ddot T^{2,\infty}\Da \ddot T^{\sr 1}_{l\pr-1}\cr\end{array}\right)
                             =\left(\begin{array}{c}\ddot T^1 \dar \om^{\sr 1}(T)\cr 
\dot T^{2,\infty}\cr\end{array}\right)$$
      where $(\ddot T^1 \dar \om^{\sr 1}(T))=(\ddot T^{\sr 1}_{\sr 1},\cdots,\ddot T^{\sr 1}_{l\pr-2},\om^{\sr 1}(T),
\ddot T^{\sr 1}_{l\pr},
      \infty \cdots).$ This again provides us the conclusion of (iii), namely,
   $$S=\left(\begin{array}{c}\ddot T^1 \dar a_{n-\sr 1}\cr \dot T^{2,\infty}\Da \om^{\sr 1}(T)\cr
\end{array}\right)
                       =\left(\begin{array}{c}\dot T^1\cr 
                                              \dot T^{2,\infty}\Da \om^{\sr 1}(T)\cr\end{array}\right).$$
\end{itemize}
\end{itemize}
\QED
\subsubsection{}\label{3.3.5} 
Now we are ready to prove theorem \ref{3.3.3}
\par
\Pf
Let $\Dscr\prpr(T)=\{T(ws_{n-\sr 1})\}$ for $w$ satisfying 
conditions of \ref{3.3.4}(iii). It is immediate from lemma \ref{3.3.4} 
${\Dscr}_n(T)\st \Dscr\prpr(T)\st \Dscr(T).$ Thus we must show
that $\Dscr\prpr(T)={\Dscr}_n\pr(T).$ 
\par
First we show that condition \ref{3.3.3}(i) is satisfied
if and only if there exists $w\in {\Cscr}_T$ satisfying
\ref{3.3.4} iii(a) and  in this case $T(ws_{n-\sr 1})=S_T(c_{\sr 1}).$
Then we show that if condition \ref{3.3.3}(ii) is satisfied
then there exists $w\in {\Cscr}_T$ satisfying
\ref{3.3.4} iii(b) and if there exists $w\in {\Cscr}_T$ satisfying
\ref{3.3.4} iii(b) then condition \ref{3.3.3} (i) or condition \ref{3.3.3} (ii)
is satisfied and in both cases 
$T(ws_{n-\sr 1})=S_T(c_{\sr 1}).$ Finally we show
that for any $c_i\ne c(1,l)$ one has $S_T(c_i)\ne \emptyset$
if and only if there exists $w\in {\Cscr}_T$ satisfying
\ref{3.3.4} iii (c) and in this case $T(ws_{n-\sr 1})=S_T(c_i).$ This
completes the proof of \ref{3.3.3}.  
\begin{itemize}
\item[(i)] If there exists $w\in {\Cscr}_T$ satisfying \ref{3.3.4} iii(a) then $\dot T^1=(\ddot T^1,a_{n-1}),\ 
T^1=(\ddot T^1,a_{n-1},a_n)$ and $T^{2,\infty}=\dot T^{2, \infty}=\ddot T^{2, \infty}.$
In particular $|T^1|=|\ddot T^1|+2\geq |\ddot T^2|+2=|T^2|+2.$ Hence  \ref{3.3.3} (i) is satisfied.
As well in this case
$$S=T(ws_{n-\sr 1})=\left(\begin{array}{c}\dot T^1\cr \dot T^{2,\infty}\Da \om^{\sr 1}(T)\cr
\end{array}\right)
=\left(\begin{array}{c}T^1-\om^{\sr 1}(T)\cr
T^{2,\infty}\Da \om^{\sr 1}(T)\cr\end{array}\right)=S_T(c_{\sr 1})$$ 
On the other hand if hypotheses of \ref{3.3.3}(i) is satisfied then  $w_r(T)$ satisfies \ref{3.3.4} iii (a).
\item[(ii)] If there exists $w\in {\Cscr}_T$ satisfying \ref{3.3.4} iii(b) then 
$\dot T^1=(\ddot T^1_1,\ldots, \ddot T^1_{l\pr-1},a_{n-1}),\ 
T^1=(\dot T^1,a_n)$ and $T^{2,\infty}=\dot T^{2, \infty}=(\ddot T^{2, \infty}\Da \ddot T^{\sr 1}_{l\pr}).$
In particular $|T^1|=|\ddot T^1|+1=l\pr+1$ and $\om^{\sr 2}(T)\geq \ddot T^{\sr 1}_{l\pr}>a_n=\om^{\sr 1}(T).$
Note that $|T^2|=|\dot T^2|\leq |\dot T^1|=l\pr.$ 
If  $|T^2|=l\pr$ then $|T^1|=|T^2|+1$ so that \ref{3.3.3}(ii) is satisfied.
If $ |T^2|<l\pr$ then $|T^1|\geq |T^2|+2$ so that \ref{3.3.3}(i) is satisfied.
In both cases
$$S=T(ws_{n-\sr 1})=\left(\begin{array}{c}\dot T^1\cr \dot T^{2,\infty}\Da \om^{\sr 1}(T)\cr
\end{array}\right)
=\left(\begin{array}{c}T^1-\om^{\sr 1}(T)\cr
T^{2,\infty}\Da \om^{\sr 1}(T)\cr\end{array}\right)=S_T(c_{\sr 1})$$
On the other hand if \ref{3.3.3}(ii) is satisfied then $w_r(T)$ satisfies \ref{3.3.4} iii (b).
Indeed in that case by RS procedure $\ddot T^1_{l\pr -1}=T^1_{l\pr-1}<T^1_{l\pr}=a_{n-\sr 1}$
and $\ddot T^1_{l\pr}=\om^{\sr 2}(T)>\om^{\sr 1}(T)=a_n.$
\item[(iii)]  If there exists $w\in {\Cscr}_T$ satisfying \ref{3.3.4} iii(c) then $a_n=\om^{\sr 1}(T)$ and
$\dot T=(T\Ua c_i)$ where $c_i$ is the corner of $(\om^{\sr 1})(\dot T).$
In that case
$$S=T(ws_{n-\sr 1})=\left(\begin{array}{c}\dot T^1\cr \dot T^{2,\infty}\Da \om^{\sr 1}(T)\cr
\end{array}\right)
=\left(\begin{array}{c}(T\Ua c_i)^1\cr
 (T\Ua c_i)\Da \om^{\sr 1}(T)\cr\end{array}\right)=S_T(c_i).$$
On the other hand let us show that if there exists $c_i$ such that $(c_i)^T=\om_{\sr 1}(T)$
for some $c_i\ne c(1,l)$ and \ref{3.3.4} iii (c) is not satisfied then  $\sh(T)\not< \sh(S)$
where 
$$S=\left(\begin{array}{c}(T\Ua c_i)^1\cr
 (T\Ua c_i)\Da \om^{\sr 1}(T)\cr\end{array}\right).$$ 
Set $c_i=c(r,s),\ \dot T=(T\Ua c_i).$
Let us show that we can reformulate conditions 
\ref{3.3.4} iii (c) as follows: 
$$c_i^{\sr T}=\om^{\sr 1}(T)\qquad {\rm and }\qquad 
\om^{\sr 1}(T)<\begin{cases}\om^{\sr 2}(\dot T) & {\rm if}\ r=2,\cr
                 (\om^r(\dot T))^{\dot T^{2,r-1}} & {\rm if}\ r>2.\cr
\end{cases} \eqno{(*)}$$
\end{itemize} 
Indeed, conditions \ref{3.3.4} iii(c) mean that $c_i^{\sr T}=\om^{\sr 1}(T)$ and 
there exist a corner $c\pr$ of $\dot T=(T\Ua c_i)$ such that 
${c\pr}^{\dot T}=a_{n-\sr 1}=T^{\sr 1}_{l-1}.$ 
Put $\ddot T=(\dot T\Ua c\pr)$. Note that $\ddot T^{\sr 1}_{l- 1}={c\pr}^{\dot T^{2,\infty}}.$ 
By \ref{3.3.4} iii(c) one has
that $\ddot T^{\sr 1}_{l- 1}>\om^{\sr 1}(T),T^{\sr 1}_{l- 1}.$ 
i.e. 
$${c\pr}^{\dot T}<\om^{\sr 1}(\dot T)\quad {\rm and}\quad {c\pr}^{\dot T^{2,\infty}}>\om^{\sr 1}(T). 
\eqno{(**)}$$
Note that by \ref{3.1.3} (i), (iii) for every corner $c\pr=c\pr(r\pr,s\pr)$
$${c\pr}^{\dot T}=\begin{cases}\om^{\sr 1}(\dot T) & {\rm if}\ r\pr<r,\cr
                         \leq T^{\sr 1}_{l- 1} & {\rm if}\ r\pr\geq r.\cr
\end{cases}$$
Furthermore applying \ref{3.1.3} (i) to two corners of $\dot T$ which are as well corners
of $\dot T^{2,\infty}$ we get
$c\pr>c\prpr$ implies ${c\pr}^{\dot T^{2,\infty}}\leq {c\prpr}^{\dot T^{2,\infty}}.$
Hence we can choose the smallest $c\pr=c\pr(r\pr,s\pr)$ such that $r\pr\geq r.$
Note also that $|\dot T^r|=s-1$ hence applying \ref{3.1.6} to $\dot T^{r,\infty}$ one has
${c\pr}^{\dot T^{r,\infty}}=\om^r(\dot T)$ thus by \ref{3.1.3} (ii) 
${c\pr}^{\dot T^{2,\infty}}={\om^r(\dot T)}^{T^{2,r-1}}.$  This
together with $(**)$ gives $(*)$.
\par
Hypotheses $(*)$ are not  satisfied only in two cases: if  $s=1$ or   
$s>1$ and 
$$\om^{\sr 1}(T)>\begin{cases}\om^{\sr 2}(\dot T^2) & {\rm if}\ r=2,\cr
                            {\om^r(\dot T)}^{T^{2,r-1}} & {\rm if}\ r>2.\cr
\end{cases}$$
\parno
If $s=1$ then $\sh(T)=(l,\cdots,|T^{k-1}|,1)$
and $\sh(S)$ is obtained by adding one box to 
$\sh(\dot T)=(l,\cdots,|T^{k-1}|)$ so that $\sh(S)=(l,\cdots,
|T^i|+1,\cdots)$ for some $i\ :\  1\leq i\leq k.$ Hence $\sh(S)\leq \sh(T)$ so that $ \sh(T)\not<\sh(S)$ as required.
\parno
Assume that $s>1$ and $r=2.$ Then 
$$S=\left(\begin{array}{c}\dot T^1\cr(\dot T^2+\om^{\sr 1}(T))\cr T^{3,\infty}\cr\end{array}\right)$$
and $\sh(S)=\sh(T),$ so that $ \sh(T)\not<\sh(S)$ as required.
\parno 
Finally consider the case  $s>1$ and $r>2$ and 
$\om^{\sr 1}(T)>{\om^r( \dot T)}^{T^{2,r-1}}\ \ (***).$    
Assume that  $\sh(S)>\sh(T).$ Then $\sh(S^{1,r})=\sh(\dot T^{1,r}).$ 
In particular putting $S\pr=(\dot T^{2,r}\Da \om^{\sr 1}(T))$
one gets $\sh(S\pr)=(|\dot T^2|,\cdots,|\dot T^r|,1).$ Hence $c(r+1,1)^{S\pr}=\om^{\sr 1}(T)$
and $\dot T^{2,r}=(S\pr\Ua c(r+1,1)).$ 
One has 
$$\om^r(\dot T)^{\sr T^{2,r-1}}=c\pr(r,s-1)^{\sr T^{2,r}}=
c\pr(r,s-1)^{(S\pr\Ua c(r+1,1))}>c(r+1,1)^{S\pr}=\om^{\sr 1}(T)$$
where inequality is obtained by \ref{3.1.5}. 
This contradicts to $(***).$ Hence $\sh(S)\not >\sh(T)$ as required.
\QED  
Note that the set $\Dscr\pr_n(T)$ may contain a few offsprings of type
$S_T(c_i)$ for $i>1.$ Let us illustrate this by 
\par
{\bf Example.}
\par
Let us regard $\bS_{10}$ and consider the following tableau:
$$T=\vcenter{
\halign{&\hfill#\hfill
\tabskip4pt\cr
\multispan{9}{\hrulefill}\cr
\ssa
\vb &1 & & 2  && 3 & & 4 & \ts\vb\cr
\vsa
&&&&\multispan{0}{\hrulefill}\cr
\ssa
\vb & 5 && 6 && 9 && 10 &  \ts\vb\cr
\vsa
&&&&\multispan{5}{\hrulefill}\cr
\ssa
\vb & 7 && 8 &  \ts\vb\cr
\vsa
\multispan{5}{\hrulefill}\cr}}$$
Deletion of the corners provides us $(T\Ua c_{\sr 1})=\dot T\pr,\ c_{\sr 1}^T=4,\ 
(T\Ua c_{\sr 2})=\dot T\prpr,\ c_{\sr 2}^T=4,$
where
$$\dot T\pr=\vcenter{
\halign{&\hfill#\hfill
\tabskip4pt\cr
\multispan{9}{\hrulefill}\cr
\ssa
\vb &1 & & 2  && 3 & & 10 & \ts\vb\cr
\vsa
&&&&&&\multispan{3}{\hrulefill}\cr
\ssa
\vb & 5 && 6 && 9  &  \ts\vb\cr
\vsa
&&&&\multispan{3}{\hrulefill}\cr
\ssa
\vb & 7 && 8 &  \ts\vb\cr
\vsa
\multispan{5}{\hrulefill}\cr}},\qquad
\dot T\prpr=\vcenter{
\halign{&\hfill#\hfill
\tabskip4pt\cr
\multispan{9}{\hrulefill}\cr
\ssa
\vb &1 & & 2  && 3 & & 6 & \ts\vb\cr
\vsa
&&&&\multispan{0}{\hrulefill}\cr
\ssa
\vb & 5 && 8 && 9 && 10 &  \ts\vb\cr
\vsa
&&\multispan{7}{\hrulefill}\cr
\ssa
\vb & 7  &  \ts\vb\cr
\vsa
\multispan{3}{\hrulefill}\cr}}.$$
A straightforward checking provides ${\Dscr}_n\pr(T)=\{S_T(c_{\sr 1}),\ S_T(c_{\sr 2})\}$ where
$$S_{\sr T}(c_{\sr 1})=\left(\begin{array}{c}\dot T^{\prime^1}\cr 
\dot T^{\prime^{2,\infty}}\Da 4\cr\end{array}\right)=
\vcenter{
\halign{&\hfill#\hfill
\tabskip4pt\cr
\multispan{9}{\hrulefill}\cr
\ssa
\vb&1 & & 2  && 3 & & 10 & \ts\vb\cr
\vsa
&&&&&&\multispan{3}{\hrulefill}\cr
\ssa
\vb & 4 && 6 && 9  &  \ts\vb\cr
\vsa
&&&&\multispan{3}{\hrulefill}\cr
\ssa
\vb & 5 && 8 &  \ts\vb\cr
\vsa
&&\multispan{3}{\hrulefill}\cr
\ssa
\vb & 7 &  \ts\vb\cr
\vsa
\multispan{3}{\hrulefill}\cr}},\, 
S_{\sr T}(c_{\sr 2})=\left(\begin{array}{c}\dot T^{\prime\prime^1}\cr 
\dot T^{\prime\prime^{2,\infty}}\Da 4\cr\end{array}\right)=
\vcenter{
\halign{&\hfill#\hfill
\tabskip4pt\cr
\multispan{9}{\hrulefill}\cr
\ssa
\vb&1 & & 2  && 3 & & 6 & \ts\vb\cr
\vsa
&&&&\multispan{1}{\hrulefill}\cr
\ssa
\vb & 4 && 8 && 9  && 10 &  \ts\vb\cr
\vsa
&&\multispan{7}{\hrulefill}\cr
\ssa
\vb & 5 &  \ts\vb\cr
\vsa
&&\multispan{1}{\hrulefill}\cr
\ssa 
\vb & 7 &  \ts\vb\cr
\vsa
\multispan{3}{\hrulefill}\cr}}.$$
Note also that  any
$S\in {\Dscr}_o(T)$ is obtained by $S=(S\pr \Da 4),$ where $S\pr\in\Dscr(T\pr)\cup\Dscr(T\prpr),$
in particular  $4\in<S^1>.$ Hence in that case ${\Dscr}_n(T)={\Dscr}_n\pr(T).$ 
\subsection
{\bf The set of offsprings as an induction of two subgroups}

\subsubsection{}\label{3.4.1}
In the previous section  we have described the set of offsprings  
of a cell in $\bS_n$ as formed from induced sets of offsprings of certain cells in
$\bS_{n- {\sr 1}}$ together with some new offsprings obtained by action of $s_{n-\sr 1}.$ 
Recall the  notation $\bS\pr_{n- {\sr 1}}$ from  \ref{2.2.7}. Here
we show that a set of offsprings can be represented as a union of 
induced sets of offsprings from $\bS_{n- {\sr 1}}$ and $\bS\pr_{n- {\sr 1}}$ which are two 
isomorphic subgroups of $\bS_n.$ 
\subsubsection{}\label{3.4.2}
Recall $w_o$ from \ref{2.2.4} Remark 3. 
Let $\dot w_o$
denote the unique largest element of the subgroup $\bS_{n- {\sr 1}}$ of $\bS_n.$ 
Recall the notation $\ov y$ and $s^{\sr >}_{i,j}$ from \ref{2.2.5}.   
One has
\begin{lemma}\begin{itemize} 
\item[(i)]   For any $y\in \bS_n$ one has $\ov y= yw_o.$
\item[(ii)]  $w_o$ can be decomposed as $w_o=y^{\sr -1}\ov y$ for any $y\in\bS_n.$          
\item[(iii)] In particular $w_o=\dot w_o s^{\sr >}_{n-\sr 1, 1}.$
\item[(iv)]  For any $y\in \bS_{n- {\sr 1}}$ set $\dot y=y\dot w_o.$ Then $w_o$ can be decomposed as
             $w_o=y^{\sr -1}s^{\sr >}_{n-\sr 1,1}\phi_{\sr 1}({\dot y}).$
\end{itemize}
\end{lemma}
\Pf 
Let $y=[a_{\sr 1},\cdots,a_n]$ then
    $$\ov y =[a_n,\cdots,a_{\sr 1}]=[a_{\sr 1},\cdots,a_n][n,\cdots,1]$$
which provides (i).
\par
(ii) follows from (i).
\par
To show (iii) note that 
$\dot w_o={\dot w}^{\sr -1}=[n-1,\cdots,1,n]$ and $\ov{\dot w_o}=[n,1,2,\ldots,n-1]=s^{\sr >}_{n-\sr 1, 1}.$ 
Thus (iii) is obtained by applying  (ii) to $y=\dot w_o.$
\par
To show (iv) note that by (ii) $\dot w_o=y^{\sr -1}\dot y.$ Then by (iii) 
$w_o=\dot w_o s^{\sr >}_{n-\sr 1, 1}=y^{\sr -1}\dot ys^{\sr >}_{n-\sr 1, 1}.$ 
Further note that for any $j\ :\ 1\leq j\leq n-2$ one has
$$s_js^{\sr >}_{n-\sr 1, 1}=s_{n-\sr 1}\cdots s_j s_{j+\sr 1} s_j
  \cdots s_{\sr 1}=s^{\sr >}_{n-\sr 1, 1}s_{j+1}=s^{\sr >}_{n-\sr 1, 1}\phi_{\sr 1}(s_j),$$
which implies $\dot ys^{\sr >}_{n-\sr 1, 1} =s^{\sr >}_{n-\sr 1, 1}
\phi_{\sr 1}(\dot y)$ and completes the proof.
\QED
\parno
{\bf Remark.}\ \ Let $y=[a_{\sr 1},\ldots,a_{n-\sr 1},n]\in \bS_{n-1}$ then by (i)
$\dot y=[a_{n-\sr 1},\ldots,a_{\sr 1},n]$ and by \ref{2.2.7}
$\phi_{\sr 1}({\dot y})=[1,a_{n-\sr 1}+1,\ldots,a_{\sr 1}+1].$
\subsubsection{}\label{3.4.3}
Recall the notation $(T\La c)$ and ${}^Tc$ from \ref{2.4.10}.
\par 
Given a cell $\Cscr$ in $\bS_{n- {\sr 1}}$ let  $\phi_{\sr 1}(\Cscr)$
be its displacement into $\bS\pr_{n- {\sr 1}},$ that is 
$\phi_{\sr 1}(\Cscr)=\{\phi_{\sr 1}(y)\ :\ y\in \Cscr\}.$
\begin{prop} Consider a standard Young tableau $T \in \bT_n$ with
$m$ corners $\{c_i\}_{i=1}^m,$ and set $q_i={}^T(c_i).$ Then
$${\Cscr}_T=\coprod_{i=1}^m s^{\sr >}_{q_i-\sr 1,1}\phi_{\sr 1}({\Cscr}_{(T\La c_i)}) 
.$$
\end{prop}
\Pf
Given $T\in \bT_n.$ By \ref{2.4.15}
$T^{\dagger}(y)=T(\ov y).$ Further by \ref{3.4.2}(i) $y=\ov y w_o.$ Thus we can represent ${\Cscr}_T$
as ${\Cscr}_T= {\Cscr}_{T^{\dagger}}w_o.$ Then by proposition \ref{3.1.1}:
$${\Cscr}_T=\coprod_{i=1}^m 
    s^{\sr <}_{q_i,n-\sr 1} {\Cscr}_{(T^{\dagger}\Ua c_i)}w_o.\eqno{(*)} $$ 
Note that
$${\Cscr}_{(T^{\dagger}\Ua c)}w_o=s^{\sr >}_{n-\sr 1,1}\phi_{\sr 1}({\Cscr}_{(T\La c)}).
                                                                     \eqno{(**)}$$
Indeed for any $y\in{\Cscr}_{(T^{\dagger}\Ua c)}$ one has 
by \ref{3.4.2}(iv) 
$$yw_o=yy^{\sr -1}s^{\sr >}_{n-\sr 1,1}\phi_1(\dot y)=s^{\sr >}_{n-\sr 1,1}\phi_1(\dot y).$$
Again by \ref{2.4.15} applied to $\bS_{n-1}$ one has that 
$y\in {\Cscr}_{(T^{\dagger}\Ua c)}$ iff $\dot y\in {\Cscr}_{(T^{\dagger}\Ua c)\dagger}$
and by \ref{2.4.10} $(*)\ \ {\Cscr}_{(T^{\dagger}\Ua c)\dagger}={\Cscr}_{(T\La c)}$ which completes the proof of $(**).$
\par
Finally $ s^{\sr <}_{q_i,n-\sr 1}\, s^{\sr >}_{n-\sr 1,1}=
s^{\sr >}_{q_i-\sr 1,1}$ which 
combined with $(*)$  proves the proposition.
\QED
\subsubsection{}\label{3.4.4}
Let us describe a few simple properties connected to insertion and deletion
by columns. Recall the notation from \ref{2.4.10}.
\par
Let $T$ be a non-standard Young tableau  and choose distinct 
$i,j\in {\Na}\ \setminus <T>.$  
\begin{lemma} Take $y$ such that $T(y)=T.$ Then
\begin{itemize}
\item[(i)]   $(i\Ra T)=T([i,y]).$ 
\item[(ii)]  $((i\Ra T)\La {}_Ti)=T$ and  $({}^T c\Ra(T\La c))=T$ for any corner $c$ of $T.$
\item[(iii)] $((i\Ra T)\Da j)=(i\Ra (T\Da j))=T([i,y,j]).$
\item[(iv)]  If $c,c\pr$ are distinct corners of $T$ then 
             $((T\La c)\Ua c\pr)=((T\Ua c\pr)\La c)$ and
             ${c\pr}^{(T\La c)}={c\pr}^T,\ {}^{(T\Ua c\pr)}c={}^Tc.$
\item[(v)]   If $c=c(r,s)$ is a corner of $T$ such that $|T^r|-|T^{r+1}|>1$ then
             $c\pr=c(r,s-1)$ is a corner of $(T\La c)$ with an entry 
$T^r_s$ and 
             $c\prpr=c(r,s-1)$ is a corner of $(T\Ua c)$. In that case one has
             $$((T\La c)\Ua c\pr)=((T\Ua c)\La c\prpr),\ c^T=
{c\pr}^{(T\La c)},\
             {}^Tc={}^{(T\Ua c)}c\prpr.$$
\end{itemize}
\end{lemma}
\Pf
\begin{itemize}
\item[(i)] By RS procedure one has $(T(y)\Da i)=T([y,i]).$
By the Schensted -Sch\" utzenberger theorem one has 
$(T(y)\Da i)^{\dagger}=T([i,\ov y]).$
Hence 
$$(i\Ra T(y))=(T^{\dagger}\Da i)^{\dagger}=
T^{\dagger}([\ov y,i])=T([i,y])$$
as required.
\item[(ii)] is obvious and is similar to \ref{2.4.9} $(*).$
\item[(iii)] Using twice (i) and RS insertion we get 
$((i\Ra T)\Da j)=(T([i,y])\Da j)=T([i,y,j])=
(i\Ra (T([y,j]))=(i\Ra (T\Da j)).$  
\item[(iv),] (v) are obtained  as follows. 
Set $p=c^T,\ p\pr={}^Tc\pr, 
T\pr=((T\Ua c)\La c\pr).$ 
One has
$$(p\pr\Ra (T\pr\Da p)){\buildrel {(iii)}\over =}
((p\pr \Ra T\pr)\Da p)=((T\Ua c)\Da p)
{\buildrel {\ref{2.4.9}(*)}\over =} T.$$
On the other hand 
$$T\prpr=((T\La c\pr)\Ua c)=
((\left(p\pr\Ra (T\pr\Da p)\right)\La c\pr)\Ua c)
{\buildrel {(ii)}\over =}((T\pr\Da p)\Ua c)
{\buildrel {\ref{2.4.9}}\over =} T\pr.$$ 
\end{itemize}
\QED
%
\subsubsection{}\label{3.4.5}
The knowledge of sets of offsprings for $\bS_{n- {\sr 1}}$ 
combined with  the above two corner deletions $\Ua,\ \La$ and two insertion
procedures $\Da,\ \Ra$ determines the set of offsprings of a cell in $\bS_n.$
Indeed by proposition \ref{3.4.3} 
$${\Cscr}_T=\coprod_{i=1}^m s^{\sr >}_{q_i-\sr 1,1}\phi_{\sr 1}
({\Cscr}_{(T\La c_i)})$$
and $\phi_{\sr 1}({\Cscr}_{(T\La c)})$ is a cell in $\bS\pr_{n- {\sr 1}}.$ So each $w\in \Cscr$
can be written as $w=s^{\sr >}_{q_i-\sr 1,1}y$ for some 
$y\in\phi_{\sr 1}({\Cscr}_{(T\La c_i)}).$  
Then right multiplication of $w$ by $s_{n-\sr 1}$ can be considered as
 multiplication of the corresponding 
$y\in\phi_{\sr 1}({\Cscr}_{(T\La c_i)})$
and then replacement of the resulting cell 
${\Cscr}_{ys_{n-\sr 1}}\in \bS\pr_{n-1}$ by the cell $s^{\sr >}_{q_i-\sr 1,1}{\Cscr}_{ys_{n-\sr 1}}.$ 
This means that offsprings obtained by action of $s_{n-\sr 1}$ can be read off
from sets of offsprings in $\bS\pr_{n- {\sr 1}}.$ Summarizing
\begin{theorem} Given a cell ${\Cscr}_T$ in 
$\bS_n$, with corners $c_i,\ i=1,2,\cdots,m.$
We can write 
$${\Cscr}_T=\coprod^m_{i=1} s^{\sr <}_{p_i,n-\sr 1}
{\Cscr}_{(T\Ua c_i)}=
        \coprod^m_{i=1}s^{\sr >}_{q_i-\sr 1,1}\phi_{\sr 1}({\Cscr}_{(T\La c_i)})$$
where $p_i=c^{\sr T},\ q_i={}^Tc_i.$ Then
 $$\Dscr({\Cscr}_T)=
\left\{
        \bigcup_{i=1}^m\{s^{\sr <}_{p_i,n-\sr 1}
\Dscr({\Cscr}_{(T\Ua c_i)})\}
                            \right\}
\bigcup\, 
\left\{
        \bigcup_{i=1}^m \{s^{\sr >}_{q_i-\sr 1,1}\phi_{\sr 1}
(\Dscr({\Cscr}_{(T\La c_i)}))\}
\right\}.$$
\end{theorem}
\subsubsection{}\label{3.4.6}
For the completion we add the description of the set of offsprings for a given 
cell ${\Cscr}_T$ provided  by the decomposition
$${\Cscr}_T=\coprod^m_{i=1}s^{\sr >}_{q_i-\sr 1,1}\phi_{\sr 1}({\Cscr}_{(T\La c_i)})$$
analogous to that in theorem \ref{3.3.3}. 
\par
We need the following 
very simple corollary of the  Schensted - Sch\" utzenberger theorem 
\begin{lemma} Given $T,S\in \bT_n.$ Then $T\dos S$ iff $S^{\dagger}\dos T^{\dagger}.$
\end{lemma}
\Pf
By \ref{2.4.15} $T^{\dagger}(w)=T(\ov w).$   
By \ref{2.1.5} for $w,y \in \bS_n$ one has
$y\dos w$ iff $S(y)\st S(w).$ It then suffices to note that by \ref{2.2.4} 
$S(\ov x)=R^+\setminus S(x),$ for all $x\in W.$
\QED
\subsubsection{}\label{3.4.7}
Let $T=T^{1,k}_{1,l}$ be a tableau with $m$ corners.
It is immediate from lemma \ref{3.4.6} that
$$\Dscr({\Cscr}_T)=\{{\Cscr}_{S^{\dagger}}\ |\ {\Cscr}_{T^{\dagger}}\in \Dscr({\Cscr}_S)\}$$ 
Set
$${\Dscr}_v({\Cscr}_T):=\bigcup_{i=1}^m\{s^{\sr >}_{q_i-\sr 1,1}
\phi_{\sr 1}(\Dscr({\Cscr}_{(T\La c_i)}))\}\quad {\rm and}\quad {\Dscr}_v(T):=\bigcup_{i=1}^m
\{(q_i\Ra S)\ :\ S\in \Dscr(T\La c_i)\} $$
For any corner $c_i\ne c(k,1)$ set 
$\dot T_{2,\infty}(c_i):=(T_{2,\infty}\La c_i)$ and $d_i:={}^{T_{2,\infty}}c_i.$
\begin{itemize}
\item[(i)] If $d_i>T^k_1$ then set
$${}_TS(c_i)=((T_1+d_i),\dot T_{2,\infty}(c_i))$$
\item [(ii)] If $ T^{k-1}_1<d_i<T^k_1$ then set 
$$S=:\left(d_i\Ra\left(T_1,\dot T_{2,\infty}(c_i)\right)\right).$$
\end{itemize}
Set
$${}_TS(c_i)=\begin{cases}S & {\rm if}\ \sh S>\sh T\cr
                    \emptyset & {\rm otherwise}\cr\end{cases}$$
\par 
As well set ${}_TS(c_i):=\emptyset$ if $d_i<T^{k-1}_1$ or $c_i=c(k,1).$ Set 
${\Dscr}_1(T)=\{{}_TS(c_i)\}_{i=1}^m.$
\begin{theorem}  $\Dscr(T)={\Dscr}_v(T)\bigcup {\Dscr}_1(T).$
\end{theorem}
\Pf
\ref{3.4.3} together with \ref{3.3.1} and \ref{3.4.6} show that
${\Dscr}_v(T)\st \Dscr(T).$ To get all the set $\Dscr(T)$ 
we must add to ${\Dscr}_v(T)$ the offsprings, obtained by the right multiplication
of the elements of $\Cscr$ by $s_{\sr 1}.$ Again by  \ref{3.4.6} ${\Cscr}_T\dor {\Cscr}_S$ 
if $S^{\dagger}\dor T^{\dagger}.$ Note that $\ov{ys_{n-1}}=s_{\sr 1}\ov y.$ 
Hence to complete the proof we must show the bijection $(T,{}_TS(c_i))^\dagger=(S_{T\pr}(c_j),T\pr).$ 
\par
Note that (i) corresponds \ref{3.3.3}(i) or (ii). Indeed, in this case 
$$ T\pr:={}_TS(c_i)^{\dagger}=\left(\begin{array}{c}(T_1+d_i)\cr 
\dot T_{2,\infty}^{\dagger}\end{array}\right)
,\quad {\rm and}\quad  \om^{\sr 1}(T\pr)=d_i $$ 
so that condition \ref{3.3.3}(i) or (ii) holds for $T\pr$ and 
$$T^{\dagger}=\left(\begin{array}{c}
({T\pr}^1-\om^{\sr 1}(T\pr))\cr 
({T\pr}^{2,\infty}\Da \om^{\sr 1}(T\pr))\cr\end{array}\right)
=S_{T\pr}(c_{\sr 1}).$$
\par
A similar computation shows that $T$ satisfies (ii) iff $T\pr={}_TS(c_i)^{\dagger}$
is such that $T^\dagger=S_{T\pr}(c_j)$ as defined by \ref{3.3.3}$(**).$ 
\QED

\section {\bf Closures of orbital varieties in $\gs\gl_n$ and Induced Duflo Order}  

\subsection{\bf Projections and embeddings of Duflo offsprings}

\subsubsection{}\label{4.1.1}  
Recall the notation from \ref{2.1.8} and the notation of 
$\pi_{i,j}$ from \ref{2.4.16}. For $w\in W$ set $w_{\sr \Iscr}:=\pi_{\sr \Iscr}(w).$
This element can be regarded as an element of $W_\Iscr$ and as an element of $W.$
Let ${\Cscr}_{w_{\sr \Iscr}}$ denote its cell in $W$ and ${\Cscr}^{\Iscr}_{w_{\sr \Iscr}}$ denote
its cell in $W_\Iscr.$ Respectively let ${\Vscr}_{w_{\sr \Iscr}}$ be the corresponding orbital variety
in $\gog$ and ${\Vscr}^{\Iscr}_{w_{\sr \Iscr}}$ be the corresponding orbital variety
in $\gl_{\sr \Iscr}.$ All the projections  
are in correspondence on orbital varieties and cells, namely
\begin{theorem} Let $\gog$ be a reductive algebra. Let $\Iscr\st \Pi.$ 
\begin{itemize}
\item[(i)]  For every $w\in W$ one has  $\pi_{\sr \Iscr}(\Cscr(w))\st \pi_{\sr \Iscr}(\Cscr(w_{\sr \Iscr}))=
            {\Cscr}_{\Iscr}(w_{\sr \Iscr}).$
\item[(ii)] For every orbital variety ${\Vscr}_w\st \gog$ one has $\pi_{\sr \Iscr}(\ov{\Vscr}_w)=
\ov{{\Vscr}^{\Iscr}_{w_\Iscr}}.$ 
\item[(iii)] If $\gog=\gs\gl_n$ and $\Iscr=\Pi_{<i,j>}$ then one also has $\pi_{i,j}(T(w))=T(\pi_{i,j}(w))$
where $\pi_{i,j}(w)$ is considered as an element of ${\boS}_{<i,j>}.$
\end{itemize}
\end{theorem}
\Pf
We prove 
(i), (ii) simultaneously. Note that lemma \ref{2.1.8} gives
$$\gn_{\sr \Iscr}\cap^w\gn=\gn_{\sr \Iscr}\cap^{w_{\sr \Iscr}}\gn_{\sr \Iscr}.$$
This implies that the projections 
$\pi_{\sr \Iscr}:\bB\rar \bB_\Iscr,\ \pi_{\sr \Iscr}:
\gn\rar\gn_{\sr \Iscr}$ and  $\pi_{\sr \Iscr}:W\rar W_\Iscr$ 
satisfy the following condition
$$\pi_{\sr \Iscr}(\bB(\gn\cap^w\gn))=\pi_{\sr \Iscr}(\bB)(\pi_{\sr \Iscr}(\gn)\cap^{\pi_{\sr \Iscr}(w)}
\pi_{\sr \Iscr}(\gn))=
  \bB_\Iscr(\gn_{\sr \Iscr}\cap^{w_{\sr \Iscr}}\gn_{\sr \Iscr}).\eqno{(*)}$$
By the continuity of the map $\pi_{\sr \Iscr}: \gn 
\rar \gn_{\sr \Iscr}$  in
Zariski topology  and $(*)$ we have
$$ \pi_{\sr \Iscr}(\ov{\bB(\gn\cap^w\gn)})\st
\ov{\bB_\Iscr(\gn_{\sr \Iscr}\cap^{w_{\sr \Iscr}}\gn_{\sr \Iscr})} \eqno{(**)} $$
Let us show the opposite inclusion. One has
$$\ov{\bB(\gn \cap^w\gn)} \supset\ov{\bB(\gn_{\sr \Iscr} \cap^{w_{\sr \Iscr}}\gn_{\sr \Iscr})}
\supset \ov{\bB_\Iscr(\gn_{\sr \Iscr} \cap^{w_{\sr \Iscr}} \gn_{\sr \Iscr})}$$
Using the fact that $\pi_{\sr \Iscr}\vert_{\gn_{\sr \Iscr}}={\rm id}$ we get 
$$\pi_{\sr \Iscr}(\ov{\bB(\gn \cap^w\gn)}) \supset
\pi_{\sr \Iscr}(\ov{\bB_\Iscr(\gn_{\sr \Iscr} \cap^{w_{\sr \Iscr}}\gn_{\sr \Iscr})})=
\ov{\bB_\Iscr(\gn_{\sr \Iscr} \cap^{w_{\sr \Iscr}}\gn_{\sr \Iscr})},$$
which together with $(**)$ provides (ii).
\par
Applying (ii) we get that for any $y\in \Cscr(w),$ that is such that 
$ \ov{\bB(\gn\cap^w\gn)})=\ov{\bB(\gn\cap^y\gn)}),$ one has 
$$\ov{\bB_\Iscr(\gn_{\sr \Iscr}\cap^{w_{\sr \Iscr}}\gn_{\sr \Iscr})}=\ov{\bB_\Iscr(\gn_{\sr \Iscr}
\cap^{y_{\sr \Iscr}}\gn_{\sr \Iscr})}$$ 
so that 
$$\pi_{\sr \Iscr}(\Cscr(w))\st {\Cscr}_\Iscr(w_{\sr \Iscr}).\eqno{
\binom{**}{*}}$$ 
\par
On the other hand consider $y_{\sr \Iscr}\in {\Cscr}_\Iscr(w_{\sr \Iscr}).$ One has that
$$\ov{\bB_\Iscr(\gn_{\sr \Iscr}\cap^{w_{\sr \Iscr}}\gn_{\sr \Iscr})}=\ov{\bB_\Iscr(\gn_{\sr \Iscr}
\cap^{y_{\sr \Iscr}}\gn_{\sr \Iscr})}$$ 
and by ~\cite[2.2.2]{Ca}  one has that for any $u\in W_\Iscr$
$$\gn\cap^u\gn=(\gn_{\sr \Iscr}\cap^u\gn_{\sr \Iscr})\oplus \gm_{\sr \Iscr}.
\eqno{\binom{**}{**}}$$
Note also that $\gm_{\sr \Iscr}$ is $\bB_\Iscr$ stable. Thus 
$$\begin{array}{rcl}
\ov{\bB(\gn\cap^{w_{\sr \Iscr}}\gn)}&=
\ov{\bM_\Iscr(\ov{\bB_\Iscr((\gn_{\sr \Iscr}\cap^{w_{\sr \Iscr}}\gn_{\sr \Iscr})\oplus \gm_{\sr \Iscr})})} 
                                                                          &{\rm by}\ \binom{**}{**}\cr 
&=\ov{\bM_\Iscr(\ov{\bB_\Iscr(\gn_{\sr \Iscr}\cap^{w_{\sr \Iscr}}\gn_{\sr \Iscr})}+ \gm_{\sr \Iscr})}
                                                                          &{\rm by\ {\bB_\Iscr}\ stability}\cr 
&=\ov{\bM_\Iscr(\ov{\bB_\Iscr(\gn_{\sr \Iscr}\cap^{y_{\sr \Iscr}}\gn_{\sr \Iscr})}+ \gm_{\sr \Iscr})} 
                                                                   &{\rm by\ hypothesis}\cr
&=\ov{\bB(\gn\cap^{y_{\sr \Iscr}}\gn)}
                                                                   &{\rm as\ above.}\cr
\end{array}$$  
Hence ${\Cscr}_\Iscr(w_{\sr \Iscr})\st\Cscr(w_{\sr \Iscr})$ which provides 
together with $\binom{**}{*}$ the assertion (i).
\par
(iii) is a straightforward corollary of \ref{2.1.8} and \ref{2.4.13}.
\QED
\subsubsection{}\label{4.1.2}
As an immediate corollary of \ref{4.1.1} (ii) we get that $\pi_{\Iscr}$ respects the geometric order,
 namely
\begin{cor} Let $\gog$ be a reductive algebra. Let $\Iscr\st \Pi.$
\begin{itemize}
\item[(i)]  If $\Vscr,\Wscr$ are orbital varieties of $\gog$ such that
 $\Vscr\go\Wscr$ then $\pi_{\sr \Iscr}(\Vscr)\go\pi_{\sr \Iscr}(\Wscr).$
\item[(ii)] If $y,w\in W$ are such that $y\go w$ then $\pi_{\sr \Iscr}(y)\go \pi_{\sr \Iscr}(w).$
\end{itemize}
\end{cor}
\subsubsection{}\label{4.1.3}
Recall $T^{<i,j>}$ from \ref{2.4.16}. 
\par
Induced Duflo order is also preserved under the projection 
$\pi_{\sr \Iscr}$ which is obvious from \ref{2.1.8}. 
Moreover if $\gog=\gs\gl_n$ one has
\begin{prop} Let $T,S\in\bT_n.$ If  $S\in \Dscr(T)$ then
for any  $i,j\ :\ 1\leq i<j\leq n$ one has $S^{<i,j>}\in\Dscr(T^{<i,j>}).$ 
\end{prop}
\Pf
Since $S$ is an offspring of $T$  there exists $w\in {\Cscr}_T$ and $\al_m\in \Pi_{n-{\sr 1}}$
such that $ws_m\dg w$ and $T(ws_m)=S.$
It is enough to show the proposition for $\pi_{{\sr 1},n-\sr 1}$ and $\pi_{{\sr 2},n}.$ 
Applying induction to these two cases we get the result for any  $\pi_{i,j}.$
The proofs for $\pi_{{\sr 1},n-\sr 1}$ and $\pi_{{\sr 2},n}$ are exactly the same
so we will show the proposition only for $\pi_{{\sr 1},n-\sr 1}.$ 
By \ref{2.4.13} (i) $w=ys^{\sr >}_{n-{\sr 1},i}$
where $y\in {\Cscr}_{(T-n)}.$ One has
$$ws_m=\begin{cases}ys_ms^{\sr >}_{n-{\sr 1},i}& {\rm if}\ m<i-1\cr
              ys^{\sr >}_{n-{\sr 1},i-{\sr 1}}& {\rm if}\ m=i-1\cr
              ys_{m-\sr 1}s^{\sr >}_{n-{\sr 1},i}& {\rm if} m>i\cr
\end {cases}$$
which implies the result. Note that in case $m=i-1$ we get 
$S^{<1,n-1>}=T^{<1,n-1>}$ but by our convention
in \ref{2.5.2}  every tableau is an offspring of itself.
\QED
\subsubsection{}\label{4.1.4}
Offsprings are
also preserved under embeddings $\Da,\Ra:\bT_n\rar\bT_{n+1}.$ Formally  
\begin{prop}. Let $T,S\in\bT_n.$ If  $S\in \Dscr(T)$ then
$(\phi_i(S)\Da i)\in\Dscr(\phi_i(T)\Da i)$ and 
$(i\Ra \phi_i(S))\in\Dscr(i\Ra \phi_i(T)).$
\end{prop}
\Pf 
Since $S$ is a offspring of $T$  there exists $w\in {\Cscr}_T$ and $\al_m\in \Pi_{n- {\sr 1}}$
such that $ws_m\dg w$ and $T(ws_m)=S.$
By \ref{2.4.7} for any $y\in \bS_n$ one has 
$(\phi_i(T(y))\Da i)=T(s^{\sr <}_{i,n}y)$ hence the first statement is just the reformulation of
\ref{3.3.1}. 
\par
Using the duality $T^\dagger(y)=T(\ov y)$ we get by \ref{2.4.7} that for any $y\in \bS_n$ one has 
$(i\Ra \phi_i(T(y)))=T(s^{\sr >}_{1,i-1}y)$ hence the second statement is equivalent to the first
by \ref{3.4.6}.
\QED
\subsubsection{}\label{4.1.5}
Note that the last  proposition provides us the following
\begin{cor} Given $S,T\in \bT_n$ (resp. ${\boT}_n$.) 
\begin{itemize}
\item[(i)]  Suppose $T^1=S^1.$ If $S^{2,\infty}$ is an offspring of $T^{2,\infty}$ then 
            $S$ is an offspring of $T.$ 
\item[(ii)] Suppose $T_1=S_1.$ If $S_{2,\infty}$ is an offspring of $T_{2,\infty}$ then 
            $S$ is an offspring of $T.$ 
\end{itemize}
\end{cor}
\Pf
It is enough to show only (i) since (ii) is equivalent to (i) by duality $T^\dagger(y)=T(\ov y)$
and \ref{3.4.6}.
\par
If $S^{2,\infty}$ is an offspring of $T^{2,\infty}$ then 
\ref{3.2.1} (i) and inductive use of \ref{4.1.4} gives that
$$S=((\cdots(S^{2,\infty}\Da T^{\sr 1}_{\sr 1})\cdots)\Da \om^{\sr 1}(T)) $$ 
is an offspring of 
$T=((\cdots(T^{2,\infty}\Da T^{\sr 1}_{\sr 1})\cdots)\Da \om^{\sr 1}(T)).$
\QED
\subsubsection{}\label{4.1.6}
The converse to corollary \ref{4.1.5} is not true as it is shown by the following 
\par
\vbox{\vskip 1truecm
{\bf Example.}\ \ \ 
Regard $\bT_{10}$
$$T=\vcenter{
\halign{&\hfill#\hfill
\tabskip4pt\cr
\multispan{9}{\hrulefill}\cr
\ssa
\vb &1 & & 2  && 3 & & 9 & \ts\vb\cr
\vsa
&&&&&&\multispan{3}{\hrulefill}\cr
\ssa
\vb & 4 && 5 && 6  &  \ts\vb\cr
\vsa
&&&&&&\multispan{0}{\hrulefill}\cr
\ssa
\vb & 7 && 8 && 10 &  \ts\vb\cr
\vsa
\multispan{7}{\hrulefill}\cr}}$$}
One can easily check that $T=T(w)$ where $w=[4,7,8,5,6,10,9,1,2,3].$ By 
\ref{2.2.4} Remark 2  $ws_{\sr 4}\dg w$ and 
$$S=T(ws_{\sr 4})= \vcenter{
\halign{&\hfill#\hfill
\tabskip4pt\cr
\multispan{9}{\hrulefill}\cr
\ssa
\vb &1 & & 2  && 3 & & 9 & \ts\vb\cr
\vsa
&&&&&&\multispan{3}{\hrulefill}\cr
\ssa
\vb & 4 && 5 && 8  &  \ts\vb\cr
\vsa
&&&&\multispan{3}{\hrulefill}\cr
\ssa
\vb & 6 && 10 &  \ts\vb\cr
\vsa
&&\multispan{3}{\hrulefill}\cr
\ssa
\vb & 7 &  \ts\vb\cr
\vsa
\multispan{3}{\hrulefill}\cr}}$$
Set
$$\dot T:=T^{2,\infty}=\vcenter{
\halign{&\hfill#\hfill
\tabskip4pt\cr
\multispan{7}{\hrulefill}\cr
\ssa
\vb & 4 && 5 && 6  &  \ts\vb\cr
\vsa
&&&&&&\multispan{0}{\hrulefill}\cr
\ssa
\vb & 7 && 8 && 10 &  \ts\vb\cr
\vsa
\multispan{7}{\hrulefill}\cr}}
\qquad{\rm and}\qquad
\dot S:=S^{2,\infty}= \vcenter{
\halign{&\hfill#\hfill
\tabskip4pt\cr
\multispan{7}{\hrulefill}\cr
\ssa
\vb & 4 && 5 && 8  &  \ts\vb\cr
\vsa
&&&&\multispan{3}{\hrulefill}\cr
\ssa
\vb & 6 && 10 &  \ts\vb\cr
\vsa
&&\multispan{3}{\hrulefill}\cr
\ssa
\vb & 7 &  \ts\vb\cr
\vsa
\multispan{3}{\hrulefill}\cr}}.$$
Let us show that $\dot S\not \dg \dot T.$ Indeed $\sh(\dot S)$
is a descendant of $\sh(\dot T)$ by \ref{2.3.3}. 
So that if $\dot S\dg \dot T$ it is an offspring of $\dot T.$ 
The only corner of $\dot T$ is $c(2,3).$ Since $c^{\dot T}=6$ and
$6\not\in <\dot S^1>$ we must have $\dot S\in {\Dscr}_n(\dot T)=\{S_T(c_{\sr 1})\}.$
The straightforward computation gives that $S\ne S_T(c_{\sr 1})$ thus
$\dot S\not \dg \dot T.$ 
\parno
{\bf Remark.}\ \ This example is important since in Part III we show that 
$\phi(\dot T)\gos \phi(\dot S)$ 
which shows that induced Duflo order is strictly weaker than geometric order.  
\subsubsection{}\label{4.1.7}
The following lemma describes the possible first rows of offsprings of 
a given $T.$
\begin{lemma} If $S$ is an offspring of $T$ and $S^1\ne T^1$ then $S^1$ is one
of the follows
\begin{itemize}
\item[(i)]   $S^1=(T^1-T^{\sr 1}_j)$ for some $j.$
\item[(ii)]  $S^1=((T^1-T^{\sr 1}_j)+b)$ for some $j$ and 
             $b>T^{\sr 1}_{j+\sr 1}.$
\item[(iii)] $S^1=(T^1\uar b).$
\end{itemize}
In particular given $T,\ S$ then $T\dos S$ implies $T^{\sr 1}_j\leq S^{\sr 1}_j.$
\end{lemma}
\Pf
We prove this by induction on $n.$
It holds for $n=2.$ Suppose it is true for $\bT_{n- {\sr 1}}$ and take $T\in \bT_n.$ 
Set $l=|T^1|.$
\par
If $n$ occurs in the corner $c=c(i,j), i>1$ of $T$ then $(T-n)^1=T^1.$
By \ref{4.1.3} $(S-n)$ is an offspring of $(T-n)$ hence by the induction hypothesis
the pair $(T-n),\ (S-n)$ satisfies one of the conclusions. Suppose this pair satisfies (ii) or 
(iii). Then $|(S-n)^1|=l$ and 
since $|S^1|\leq l$ one has that $S^1=(S-n)^1$ hence (ii) or (iii) 
is satisfied for $T,S.$
Now suppose the pair $(T-n),\ (S-n)$ satisfies (i). Then either $S^1=(S-n)^1$ and $T,S$ satisfy (i) 
or $S^1=((S-n)^1+n)=
((T-T^{\sr 1}_j)+n)$ and (ii) for $b=n$ is 
satisfied for $T,S.$
\par
Suppose $n$ occurs in the corner $c(1,l)$ of $T.$ By theorem \ref{3.3.3} either
there exist a corner $c$ such that 
$S=(\dot S\Da c^T)$ where $\dot S\in\Dscr(T\Ua c)$
or T satisfies \ref{3.3.3} (i) and $S=S_T(c_{\sr 1}).$  If $S=S_T(c_{\sr 1})$ then 
$S^{\sr 1}=(T^1-n)$ so (i) is satisfied. Otherwise by
the induction hypothesis the pair $\dot T,\ \dot S$ satisfies one of 
the conclusions.
If $c=c(1,l)$ then $T^1=(\dot T^1+n),\ S^1=(\dot S^1 +n),$ hence $T,S$ satisfy the same
conclusion as $\dot T,\dot S.$
It remains to consider $c>c(1,l).$ In this case $|\dot T^1|=l$ and $c^T<n.$
 
Suppose the pair $\dot T,\dot S$ satisfies (i), this is $\dot S^1=
(\dot T^1-\dot{T}^{\sr 1}_j).$ Then
$$S^1=((\dot T^1-\dot{T}^{\sr 1}_j)\dar c^T)=
\begin{cases}(T^1-T^{\sr 1}_j) & {\rm if}\ c^T\ne T^{\sr 1}_j,\cr
             (T^1-T^{\sr 1}_{j+1}) & {\rm if}\ c^T=T^{\sr 1}_j.\cr
\end{cases}$$
In both cases the pair $T,S$ satisfies (i).
Finally suppose the pair $\dot T,\dot S$ satisfies (ii) or (iii). Then  
$|\dot S^1|=|\dot T^1|=l$ and furthermore
$\dot{S}^{\sr 1}_j\geq \dot{T}^{\sr 1}_j,\ \forall j.$ 
So in particular  $\om^{\sr 1}(\dot S)=n$ and $c(1,l)$ is a corner of $\dot S.$
Now $\dot S\in\Dscr(\dot T)$ implies by \ref{4.1.3} that
$(\dot S-n)\in\Dscr(\dot T-n).$ Further by \ref{4.1.4} 
$S\pr=((\dot S-n)\Da c^T)=(S-n)$ is an 
offspring of $T\pr=((\dot T-n)\Da c^T)=(T-n).$ Since 
$|{S\pr}^1|=|{T\pr}^1|=l-1$ one has the pair $T\pr,\ S\pr$ satisfies (ii) 
or (iii) by induction hypothesis. Note that $T^1=({T\pr}^1+n),\ S^1=({S\pr}^1+n),$ 
hence the pair $T,S$ satisfies the same conclusion as $T\pr,\ S\pr.$
\QED
For example let 
$$T=\vcenter{
\halign{&\hfill#\hfill
\tabskip4pt\cr
\multispan{15}{\hrulefill}\cr
\ssa
\vb &1 & & 2  && 3 & & 4 & & 5& & 6 & & 8 & \ts\vb\cr
\vsa
&&&&&&\multispan{9}{\hrulefill}\cr
\ssa
\vb & 7 && 9 && 10  &  \ts\vb\cr
\vsa
\multispan{7}{\hrulefill}\cr}}.$$
Then $S$ is an offspring described in (i), $U$ is an offspring
described in (ii) and $V$ is an offspring described in (iii)
of lemma \ref{4.1.7} where
$$S=\vcenter{
\halign{&\hfill#\hfill
\tabskip4pt\cr
\multispan{13}{\hrulefill}\cr
\ssa
\vb &1 && 3 & & 4 & & 5& & 6 & & 8 & \ts\vb\cr
\vsa
&&&&&&&&\multispan{5}{\hrulefill}\cr
\ssa
\vb & 2&& 7 && 9 && 10  &  \ts\vb\cr
\vsa
\multispan{9}{\hrulefill}\cr}},$$
$$U=\vcenter{
\halign{&\hfill#\hfill
\tabskip4pt\cr
\multispan{15}{\hrulefill}\cr
\ssa
\vb &1 & & 2  && 3 & & 4 & & 5& & 8 & & 10 & \ts\vb\cr
\vsa
&&&&\multispan{11}{\hrulefill}\cr
\ssa
\vb & 6 && 9 &\ts\vb\cr
\vsa
&&\multispan{3}{\hrulefill}\cr
\ssa
\vb & 7&\ts\vb\cr
\vsa
\multispan{3}{\hrulefill}\cr}},\quad
V=\vcenter{
\halign{&\hfill#\hfill
\tabskip4pt\cr
\multispan{15}{\hrulefill}\cr
\ssa
\vb &1 & & 2  && 3 & & 4 & & 5& & 6 & & 10 & \ts\vb\cr
\vsa
&&&&\multispan{11}{\hrulefill}\cr
\ssa
\vb & 7 && 8 &\ts\vb\cr
\vsa
&&\multispan{3}{\hrulefill}\cr
\ssa
\vb & 9&\ts\vb\cr
\vsa
\multispan{3}{\hrulefill}\cr}}$$
\subsubsection{}\label{4.1.8}
This lemma provides 
\begin{prop} Given two diagrams $D_1<D_2$ and a Young tableau $T$ such
that $\sh(T)=D_1$ there exists $S$ such that $T\dos S$ and $\sh(S)=D_2.$  
\end{prop}
\Pf
It is enough to show that for any  $D_2$ descendant of $D_1$ there exist $S$  of shape  $D_2$
such that $S\dg T$.  $\sh(S)$ being a descendant of $\sh(T)$ forces $S$ to be an offspring of $T.$
\par
We prove the statement by induction on $n.$ For $\bS_2,\bS_3$ the statement
is trivial. Suppose it holds for $n\pr<n$ and take $D_1,D_2\in \bD_n.$ 
\par
If $D_1=(n),\ D_2=(n-1,1)$ then the assertion is trivial.
Hence we can assume that $D_1=(\l_{\sr 1},\l_{\sr 2},\cdots,\l_k)$ where $\l_{\sr 1}<n$ and $D_2=(\mu_{\sr 1},\cdots).$ 
\par
If $\mu_{\sr 1}=\l_{\sr 1}$ then by induction hypothesis we can find an 
offspring $S\pr$ of 
$T^{2,\infty}$ such that $\sh(S\pr)=D_2^{2,\infty}.$ By \ref{4.1.7} one has 
$S^{\sr 2}_j\geq T^{\sr 2}_j>T^{\sr 1}_j$  hence 
$$S=\left(\begin{array}{c}
T^1\cr S\pr\cr\end{array}\right)$$ 
is a standard Young tableau
and by \ref{4.1.5} $S$ is an offspring of $T$ and $\sh(S)=D_2.$
\par
We are reduced to the case $\mu_{\sr 1}=\l_{\sr 1}-1.$ In that case we get a result with the help of
case by case analysis.
\parno
(i)\ \ \ If $\l_{\sr 1}-\l_{\sr 2}\geq 2$ then $\mu_{\sr 1}=\l_{\sr 1}-1,\ \mu_{\sr 2}=
\l_{\sr 2}+1,\ \mu_i=\l_i$ for $i\geq 3.$ We will show that in that case there exist $T^1_i$ such that
$$S=\left(\begin{array}{l} T^1-T^1_i\cr
           T^2+T^1_i\cr
           T^{3,\infty}\cr\end{array}\right).\eqno (*)$$
Indeed this is trivially true for $\bT_2,\ \bT_3.$ Assume that this is true for $ \bT_{n\pr}$ where $n\pr\leq n-1$ 
and show for $\bT_n.$ 
\par
If $D_1$ has more than 2 rows 
let $t:=n-\l_1-\l_2.$ Consider $T^{1,2}.$ By induction hypothesis there exist $y\in {\Cscr}_T$
and $s_m$ such that $ys_m\dg y$ and 
$$T(ys_m)=\left(\begin{array}{c} T^1-T^1_i\cr
           T^2+T^1_i\cr\end{array}\right):=\dot S$$
Note that for any $j\ \dot S^2_j\leq T^2_j$ thus $S=\left(\begin{array}{l}\dot S\cr
           T^{3,\infty}\cr\end{array}\right)$ 
is a standard tableau.
By \ref{3.2.3} (iv) $T=T([w_r(T^{3,\infty}),y])$ and $S=T([w_r(T^{3,\infty}),ys_m])=T([w_r(T^{3,\infty}),y]s_{m+t}).$
Hence $S\in\Dscr(T).$
\par
Now consider that $D_1$ has only two rows. 
\begin{itemize}
\item[a)] If $\om^{\sr 1}(T)=n$ then $T$ satisfies \ref{3.3.3} (i) and
$$S_T(c_{\sr 1})=\left(\begin{array}{c}T^1-n\cr
                       T^2+n\cr\end{array}\right)$$
satisfies $(*)$.
\item[b)] If $\om^{\sr 2}=n$ 
then consider $\dot T=((T\Ua c(2,\l_2))\Ua c(1,\l_1)).$ One has
$$\dot T=\left(\begin{array}{l}T^1-\om^{\sr 1}(T)\cr
              T^2-n\cr\end{array}\right)$$
so that $\sh(\dot T)=(\l_1-1,\l_2-1)$ thus by induction hypothesis $\dot T$ has
an offspring $\dot S$ of the form 
$$\dot S=\left(\begin{array}{l}T^1-\om^{\sr 1}-T^1_i\cr
               T^2-n+T^1_i\cr\end{array}\right)$$ 
Further by \ref{4.1.4} 
$$S=((\dot S\Ua n)\Ua \om^{\sr 1}(T))= \left(\begin{array}{l}T^1-T^1_i\cr
               T^2+T^1_i\cr\end{array}\right)$$
is an offspring of $T.$ This completes case (i).
\end{itemize}

(ii)\ \ \ The case when $\l_{\sr 2}=\cdots=\l_j=\l_{\sr 1}-1,\ 
\l_{j+1}=\l_{\sr 1}-2$ for $j\geq 2$ and 
$\mu_{\sr 1}=\cdots=\mu_{j+1}=\l_{\sr 1}-1,\ \mu_i=\l_i$ for $i>j+1$ is proved in the same manner
as (i).
\begin{itemize}  
\item[a)] First let us show that if $D_1$ has $j+1$ rows
then there exist an offspring $S$ of shape 
$\mu=(\l_{\sr 1}-1,\cdots,\l_{\sr 1}-1)$
such that 
$$S^{j+1}_i\leq T^{j+1}_i\quad{\rm for\ all}\  i<\l_{\sr 1}-1 \qquad{\rm and}\qquad S^m_{\l_{\sr 1}-1}\leq T^m_{\l_{sr 1}-1}
\quad{\rm for\ all}\ m\leq j.\eqno{(*)}$$
 Indeed this is true for $\bT_3$ we can assume that it is true for
${\boT}_E,\ \vert E\vert <n.$
\begin{itemize} 
\item[1.] If $\om^{\sr 1}(T)>\om^{\sr 2}(T)$ then consider 
$$T\pr=\left(\begin{array}{l}T^2+\om^{\sr 1}(T)\cr T^{3,\infty}\cr\end{array}\right).$$
By induction hypothesis there exists a required $S\pr.$ Then by 
lemma \ref{4.1.7} ${S\pr}^{\sr 1}_i\geq T^{\sr 2}_i$ for every $1\leq i\leq l-1$
so that
$$S=\left(\begin{array}{l}T^1-\om_{\sr 1}(T)\cr S\pr\end{array}\right)$$ 
is a standard Young tableau.
Using subsequently \ref{4.1.4} we get that 
$S=(\cdots(S\pr\Da T^{\sr 1}_{\sr 1})\cdots\Da T^{\sr 1}_m)$ 
is an offspring 
of $T=(\cdots(T\pr\Da T^{\sr 1}_{\sr 1})\cdots\Da T^{\sr 1}_m)$ satisfying $(*).$
\item[2.] If $\om^{\sr 1}(T)<\om^{\sr 2}(T)$ and $\om^{j}(T)<\om^{j+1}(T)$ then consider 
$$T\pr=(\cdots((T\Ua c(j+1,\l_{\sr 1}-2))\Ua 
c(j,\l_{\sr 1}-1))\cdots\Ua c(1,\l_{\sr 1}))$$ 
this is $T$ without segment
$s_{c(j+1,\l_{\sr 1}-2)}.$
By induction hypothesis there exist an offspring 
$S\pr$ of $T\pr$ holding $(*)$ for $\l_{\sr 1}\pr=\l_{\sr 1}-1.$ Then note that $\om^{j+1}(T)=\max<T>$ hence
$\om^{j+1}(T)>\om^{j+1}(S\pr)$ which together with $(*)$ provides that
$$S=(\cdots(S\pr\Da \om^{j+1}(T))\cdots\Da \om^{\sr 1}(T))=
\left( S\pr,\left(\begin{array}{c}\om^{\sr 1}(T)\cr\vdots\cr \om^{j+1}(T)\cr\end{array}\right)\right)$$
is an offspring of $T$ satisfying $(*).$
\item[3.] If $\om^{\sr 1}(T)<\om^{\sr 2}(T)$ and $\om^{j}(T)>\om^{j+1}(T)$ then considering 
$$
T\pr=(\cdots((T\La c(j,\l_{\sr 1}-1))\La c(j+1,\l_{\sr 1}-2))\cdots\La c(j+1,1))=
\left(\begin{array}{l}T^{1,j-1}\cr T^{j}-\om^{\sr j}(T)\cr\end{array}\right)$$ 

that is $T$ without segment of $T^{\dagger}\ :\ s_{c(\l_{\sr 1}-1,j)}$ 
we get the result in the same way as in 2.
\end{itemize}
\item[b)] If $D_1$ has more than $j+1$ rows then by induction hypothesis
there exist $S\pr$ an offspring of $T^{1,j+1}$ satisfying
$(*).$ In particular 
$$S=\left(\begin{array}{l}S\pr\cr T^{j+2,\infty}\cr\end{array}\right)$$ 
is a standard tableau
and the same considerations as in (i) provide that $S$ is an offspring
of $T.$
\end{itemize}
\QED
This proposition shows that for every two nilpotent orbits 
${\Oscr}_1,{\Oscr}_2\in \gs\gl_n$ such that
${\Oscr}_2\subsetneq \ov {\Oscr}_1$ and every orbital variety 
${\Vscr}_1$ attached to ${\Oscr}_1$
one can find ${\Vscr}_2$ attached to ${\Oscr}_2$ such that
 ${\Vscr}_1\dos {\Vscr}_2$ and in particular ${\Vscr}_1\gos {\Vscr}_2.$
\subsection
{\bf Offsprings and descendants}\label{4.2}

\subsubsection{}\label{4.2.1}
By \ref{4.1.3} and \ref{4.1.4} both projections and embeddings
preserve Duflo offsprings. By \ref{4.1.1} projections respects geometric
order as well. As we will show in Part III embeddings
also respect geometric order. But neither of them preserve Duflo 
(and geometric) descendants. 
First of all it is obvious from \ref{3.3.3} that there exist 
$T,S$ such that $S$ is a descendant  of $T$ (in both geometric and 
Duflo order) but $\pi_{i,j}(S)=\pi_{i,j}(T)$ for some $i,j$
(resp. $(\phi_i(S)\Da i)=(\phi_i(T)\Da i)$ for some $i$).
As we show by the examples in the subsections
below there exist $T,S$ such that $S$ is a descendant of $T$ (in both geometric and 
Duflo order) and  
$\pi_{i,j}(S)>\pi_{i,j}(T)$ (resp. 
$(\phi_i(S)\Da i)>(\phi_i(T)\Da i)$ for some  $i$)
but $\pi_{i,j}(S)$ (reps. $(\phi_i(S)\Da i)$)  is not a 
descendant of $\pi_{i,j}(T)$ (resp. 
$(\phi_i(T)\Da i)$).  
\subsubsection{}\label{4.2.2}
Let us first show that
$S$ being a descendant of $T$ does not imply
that $(\phi_i(S)\Da i)$ is a descendant of 
$(\phi_i(T)\Da i).$ 
The first such example occurs in embedding from $\bT_3$
into $\bT_4.$ Consider  
$$T=
\vcenter{
\halign{&\hfill#\hfill
\tabskip4pt\cr
\multispan{5}{\hrulefill}\cr
\ssa
\vb & 1 &  &  2 & \ts\vb\cr
\vsa
&&\multispan{3}{\hrulefill}\cr
\ssa
\vb & 3 & \ts\vb\cr
\vsa
\multispan{3}{\hrulefill}\cr}}\quad {\rm and}\quad 
S=
\vcenter{
\halign{&\hfill#\hfill
\tabskip4pt\cr
\multispan{3}{\hrulefill}\cr
\ssa
\vb & 1 &\ts\vb\cr
\vsa
&&\cr
\ssa
\vb & 2 & \ts\vb\cr
\vsa 
&&\cr
\ssa
\vb & 3 & \ts\vb\cr
\vsa 
\multispan{3}{\hrulefill}\cr}}.$$
Note that $S>T$ in both geometric and Duflo order
because ${\Vscr}_S=\{0\}$ and $S$ is a descendant
of $T$ since ${\Oscr}_S$ is a descendant of ${\Oscr}_T.$
Now consider
$$T\pr=(T\Da 4)=\vcenter{
\halign{&\hfill#\hfill
\tabskip4pt\cr
\multispan{7}{\hrulefill}\cr
\ssa
\vb & 1 &  &  2 && 4 &\ts\vb\cr
\vsa
&&\multispan{5}{\hrulefill}\cr
\ssa
\vb & 3 & \ts\vb\cr
\vsa
\multispan{3}{\hrulefill}\cr}}\ ,\quad
S\pr=(S\Da 4)=
\vcenter{
\halign{&\hfill#\hfill
\tabskip4pt\cr
\multispan{5}{\hrulefill}\cr
\ssa
\vb & 1 && 4 &\ts\vb\cr
\vsa
&&\multispan{3}{\hrulefill}\cr
\ssa
\vb & 2 & \ts\vb\cr
\vsa 
&&\cr
\ssa
\vb & 3 & \ts\vb\cr
\vsa 
\multispan{3}{\hrulefill}\cr}}\quad {\rm and}\quad
P\pr=
\vcenter{
\halign{&\hfill#\hfill
\tabskip4pt\cr
\multispan{5}{\hrulefill}\cr
\ssa
\vb & 1 && 2 &\ts\vb\cr
\vsa
&&&&\cr
\ssa
\vb & 3 && 4& \ts\vb\cr
\vsa 
\multispan{5}{\hrulefill}\cr}}.$$
By \ref{3.3.3} $P\pr=S_{T\pr}(c_1)$ and 
$S\pr=S_{P\pr}(c_{\sr 1}).$
Thus $T\pr<P\pr<S\pr$ (both in geometric and Duflo orders) 
so that $S\pr$ is not a descendant of $T\pr.$ 
\subsubsection{}\label{4.2.3}
Now let us show that $S$ being a descendant of 
$T$ does not imply $\pi_{{\sr 2},n}(S)$
is a descendant of $\pi_{{\sr 2},n}(T).$ The first such example
occurs in projection from $\bT_5$ onto $\bT_4.$ Consider
$$T=
\vcenter{
\halign{&\hfill#\hfill
\tabskip4pt\cr
\multispan{7}{\hrulefill}\cr
\ssa
\vb & 1 && 3 && 5 &\ts\vb\cr
\vsa
&&&&\multispan{3}{\hrulefill}\cr
\ssa
\vb & 2 && 4 &\ts\vb\cr
\vsa 
\multispan{5}{\hrulefill}\cr}}\quad {\rm and} \quad 
S=
\vcenter{
\halign{&\hfill#\hfill
\tabskip4pt\cr
\multispan{5}{\hrulefill}\cr
\ssa
\vb & 1 && 3 &\ts\vb\cr
\vsa
&&&&\cr
\ssa
\vb & 2 && 4 &\ts\vb\cr
\vsa 
&&\multispan{3}{\hrulefill}\cr
\ssa
\vb & 5 &\ts\vb\cr
\vsa 
\multispan{3}{\hrulefill}\cr}}
$$
Note that $T=T(s_{\sr 1}s_{\sr 3})$ and 
$S=T(s_{\sr 1}s_{\sr 3}s_{\sr 4}s_{\sr 3})$ so that $T\dos S.$ Note also that
if there exists $P$ such that $T\go P\go S$ then 
$\pi_{\sr 1,4}(T)\go \pi_{\sr 1,4}(P)\go \pi_{\sr 1,4}(S).$ Since  
$\pi_{\sr 1,4}(T)=\pi_{\sr 1,4}(S)$ one has
$\pi_{\sr 1,4}(P)=\pi_{\sr 1,4}(T)$ which implies just by the shape 
consideration that $P=T$ or $P=S.$ We get that $S$ is a descendant of $T$ (both
in geometric and Duflo orders).
Now consider
$$T\pr=\phi^{\sr -1}_{\sr 1}\pi_{\sr 2,5}(T)=
\vcenter{
\halign{&\hfill#\hfill
\tabskip4pt\cr
\multispan{7}{\hrulefill}\cr
\ssa
\vb & 1 && 2 && 4 &\ts\vb\cr
\vsa
&&\multispan{5}{\hrulefill}\cr
\ssa
\vb & 3 &\ts\vb\cr
\vsa 
\multispan{3}{\hrulefill}\cr}}\ ,\quad 
S\pr=\phi^{\sr -1}_{\sr 1}\pi_{\sr 2,5}(S)=
\vcenter{
\halign{&\hfill#\hfill
\tabskip4pt\cr
\multispan{5}{\hrulefill}\cr
\ssa
\vb & 1 && 2 &\ts\vb\cr
\vsa
&&\multispan{3}{\hrulefill}\cr
\ssa
\vb & 3 &\ts\vb\cr
\vsa 
&&\cr
\ssa
\vb & 4 &\ts\vb\cr
\vsa 
\multispan{3}{\hrulefill}\cr}}\quad {\rm and}\quad
P\pr=\vcenter{
\halign{&\hfill#\hfill
\tabskip4pt\cr
\multispan{5}{\hrulefill}\cr
\ssa
\vb & 1 && 2 &\ts\vb\cr
\vsa
&&&&\cr
\ssa
\vb & 3 && 4 &\ts\vb\cr
\vsa 
\multispan{5}{\hrulefill}\cr}}$$
By \ref{3.3.3} one has $T\pr\dos P\pr\dos S\pr.$
\parno
{\bf Remark.}\ \ $S$ is obtained from $T$ just by moving the box
with the maximal number down to the first possible place. As we
show in Part III such $S$ is always a descendant of $T$
(both in geometric and Duflo orders).
\subsection
{\bf Definition of Induced Duflo Order }

\subsubsection{}\label{4.3.1}
One of the first general questions that arises about induced Duflo order is the
following. Let ${\Cscr}_1,\ {\Cscr}_2$ be two cells such that
${\Cscr}_1\dor {\Cscr}_2,$ can we always find representatives 
$x\in {\Cscr}_1,\ y\in {\Cscr}_2$ such that $x\dor y?$ The answer is negative 
(which is natural enough) and we  show this in the corresponding
\parno
{\bf Example.}\ \ \ Regard $\bS_5.$ Consider 
$$T=
\vcenter{
\halign{&\hfill#\hfill
\tabskip4pt\cr
\multispan{7}{\hrulefill}\cr
\ssa
\vb & 1 &  &  2 & &  5 & \ts\vb\cr
\vsa
&&&&\multispan{3}{\hrulefill}\cr
\ssa
\vb & 3 &  &  4 & \ts\vb\cr
\vsa
\multispan{5}{\hrulefill}\cr}},\qquad
{\Cscr}_T=
\left \{\begin{array}{c}
[3,1,4,2,5],\ [3,4,1,2,5],\ [3,1,4,5,2],\cr
         [3,4,1,5,2],\ [3,4,5,1,2]\cr\end{array}
\right \}.$$
Consider $x=[3,4,1,2,5]\in {\Cscr}_T,\ s_{\sr 3}x=
[3,4,2,1,5]\dgs x.$ Hence
$S=T(s_{\sr 3}x)\dg T$ and
$$S=\vcenter{
\halign{&\hfill#\hfill
\tabskip4pt\cr
\multispan{7}{\hrulefill}\cr
\ssa
\vb & 1 &  &  4 & &  5 & \ts\vb\cr
\vsa
&&\multispan{5}{\hrulefill}\cr
\ssa
\vb & 2 &  \ts\vb\cr
\vsa
&&\cr
\ssa
\vb & 3 &  \ts\vb\cr
\vsa
\multispan{3}{\hrulefill}\cr}},\ 
{\Cscr}_S=
\left \{\begin{array}{ccc}
[3,2,1,4,5],& [3,2,4,1,5],& [3,2,4,5,1]\cr
[3,4,2,1,5],& [3,4,2,5,1],& [3,4,5,2,1]\cr
\end{array}\right \}.$$
Consider $y=[3,2,1,4,5]\in {\Cscr}_S,\ s_{\sr 4}y=
[3,2,1,5,4]s_{\sr 4}\dgs y.$ Hence 
$U=T(s_{\sr 4}y)\dg S\dgs T$ and 
$$U=
\vcenter{
\halign{&\hfill#\hfill
\tabskip4pt\cr
\multispan{5}{\hrulefill}\cr
\ssa
\vb & 1 &  &  4 & \ts\vb\cr
\vsa
&&&&\multispan{3}{\hrulefill}\cr
\ssa
\vb & 2 &  &  5 & \ts\vb\cr
\vsa
&&\multispan{3}{\hrulefill}\cr
\ssa
\vb & 3 & \ts\vb\cr
\vsa
\multispan{3}{\hrulefill}\cr}},\quad  
{\Cscr}_U=\left \{\begin{array}{c}
[3,2,1,5,4],\ [3,2,5,1,4],\ [3,5,2,1,4],\cr
[3,2,5,4,1],\ [3,5,2,4,1]\cr
\end{array}
\right\}.$$
Since $t=[3,1,4,2,5]$ is a minimal element of ${\Cscr}_T$ and 
$z=[3,5,2,4,1]$ is a maximal element of ${\Cscr}_U$ it is enough to show that
$t\not\dos z.$ Indeed using \ref{2.2.4} one gets
$$\begin{array}{rcl}
\gn\cap^t\gn&=&X_{1,2}\oplus X_{1,4}\oplus X_{1,5}\oplus X_{2,5}\oplus 
              X_{3,4}\oplus X_{3,5}\oplus X_{4,5}\cr 
\gn\cap^z\gn&=&X_{2,4}\oplus X_{3,4}\oplus X_{3,5}\cr
\end{array}$$
\subsubsection{}\label{4.3.2}
The previous example  shows that in some sense the original
definition of induced Duflo order is not good. Theorem \ref{3.4.5} shows 
that induced Duflo order is the minimal partial order generalized by embeddings
$\Da$ and $\Ra$ from the natural order on $\bT_2.$

\bigskip
\parno
\centerline{ INDEX OF NOTATION}

\parno
\begin{tabular}{llll}
\ref{1.1}& ${\bG},\ \gog\quad\quad\quad\quad$ &\ref{2.2.7}& ${\bS}\pr_{n-1}$\\
\ref{1.2}& $\gn,\quad \gn^-,\quad \gh,\quad \Oscr\quad\quad\quad\quad$
&\ref{2.3.1}& $\lambda,\ \lambda^*,
\ D_{\lambda},\ {\bD}_n,\ J(u),\ {\Oscr}_{\lambda},\ D_u$\\
\ref{1.3}& $\Vscr$&\ref{2.3.2}& $D_{\lambda}\geq D_{\mu}$\\
\ref{1.5}& ${\bB},\ W$&\ref{2.4.1}& ${\bT}_n,\ \sh(T)$\\
\end{tabular}
\parno
\begin{tabular}{llll}
\ref{1.6}& ${\bS}_n$&\ref{2.4.2}& $T^i_j,\ r_{\sr T}(u),\ c_{\sr T}(u),\ h(T^i_j),\ T^i,\ T_j,\ |T_i|,$\\
\ref{1.8}& $\go,\ \Oscr_\Vscr$&& $\omega^i(T),\ \omega_i(T),\ T^{i,j},\ T^{i,\infty},\ T_{i,j},\ T_{i,\infty},$\\ 
\ref{1.9}& $\ell(w),\ \dor$&& $T^{\dagger},\ {\boT}_E,\ {\boT}_n,\ <T>,\ [T^i],\ [T_j]$\\
\ref{2.1.1}& $\Nscr, \ {\Oscr}\geq {\Oscr}\pr$&\ref{2.4.3}& $(T,S), \left( \begin{array}{c} T\cr S\cr 
\end{array}\right)$\\
\ref{2.1.2}& ${\bV},\ R,\ R^+,\ \Pi,\ X_{\al}$&\ref{2.4.4}& $c=c(i,\l_i),\ c< c\pr$\\
\ref{2.1.3}& $\gn\cap^w\gn,\ {\bH}(\gn\cap^w\gn)$&\ref{2.4.5}& $(R\dar j),\ (T\Da j),\ j_{\sr T}$\\
\ref{2.1.4}& ${\Cscr}_w,\ {\Cscr}_y\go {\Cscr}_w $&\ref{2.4.6}& 
RS, $T(w),\ {}_iT(w),\ \Cscr_T,\ 
T_{\Cscr}$ \\
\ref{2.1.5}& $S(w)$&\ref{2.4.8}& $(R-a),\, 
(R+j),\, (R\uar j),\, j^{\sr R},\, (T\uar j),\, j^{\sr T}$ \\
\ref{2.1.6}& ${\Vscr}_y\dor{\Vscr}_w,\ {\Cscr}_y\dor{\Cscr}_w$&\ref{2.4.9}& $(T\Ua c),\ c^{\sr T},\ s_c(T), (T\Ua a)$\\
\ref{2.1.7}& ${\bP}_{\al},\ \gp_{\al},\ {\bP}_{\Vscr},\ \gm_{\sr \Vscr},\ \tau(w),$&
\ref{2.4.10}& $(C-a),\ (C+j),\ (j\rar C),\ {}_{\sr C} j,$\\
           & $\tau({\bP}),\ \tau(\gp),\ \tau(\Vscr),\ \tau(\Cscr)$&&$(C\lar j),
\ {}^{\sr C}j\ (j\Ra T),\ (T\La c),\ {}^{\sr T}c$\\
\ref{2.1.8}& ${\bP}_{\Iscr},\ {\bM}_{\Iscr},\ {\bL}_{\Iscr},\ \gp_{\sr \Iscr},$&\ref{2.4.11}& $T-T^i_j$\\
           & $\gm_{\sr \Iscr},\ \gl_{\sr \Iscr},\ \gn_{\sr \Iscr},\ \pi_{\sr \Iscr},\ W_{\Iscr},$&
                                         \ref{2.4.14}& $\tau(T)$\\
           & $D_{\Iscr},\ R_{\Iscr}^+,\ w_{\sr \Iscr},\ d_{\sr \Iscr},\ f_{\sr \Iscr}$&
\ref{2.4.16}& $<i,j>,\ \Pi_{i,j},\ \pi_{i,j},\ T^{<i,j>},\ D^{<i,j>}_T,\ {\bC}_n$\\
\ref{2.2.1}& $e_{i,j},\ \al_{i,j}$&\ref{2.5.2}& ${\Dscr}(T),\ \Dscr(\Cscr)$\\
\ref{2.2.2}& $[a_{\sr 1},\ldots,a_n]$&
\ref{3.1.1}& ${\Cscr}_{(T\Ua c_i)}$\\
\ref{2.2.3}& $p_w(i)$&\ref{3.2.2}& $w_r(T),\ w_c(T)$ \\
\ref{2.2.5}& $<w>,\,(w-m),\, w_{[i,j]},\, \ov w,\, [x,y],$ &\ref{3.3.2}
& ${\Dscr}_o({\Cscr}_T),\ {\Dscr}_o(T),
\ {\Dscr}_n({\Cscr}_T),\ {\Dscr}_n(T)$\\
           & ${\boS}_n,\ {\boS}_E,\ \phi_j,\ \phi_j^{\sr -1},\ \phi,\ s^{\sr <}_{i,j},\ s^{\sr >}_{i,j}$&
                 \ref{3.3.3}& $S_T(c_i),\ {\Dscr}\pr_n(T)$\\
\end{tabular}

\end{document}